\title{Ensemble Kalman Filter Implementations Based on Covariance Matrix Estimation}
\author{
        Elias D. Nino-Ruiz and Adrian Sandu  \\
        Computational Science Laboratory, \\
        Department of Computer Science, \\
        Virginia Polytechnic Institute and State University,\\
        Blacksburg, VA 24060, USA \\
	enino@vt.edu, sandu@cs.vt.edu
        }        
        
\date{\today}

\documentclass[12pt]{article}

\usepackage{caption}
\usepackage{amsmath}
\usepackage{amssymb}
\usepackage[square,sort,comma,numbers]{natbib}
\usepackage{bbm, dsfont}
\usepackage{algorithmic}
\usepackage{algorithm}
\usepackage{graphicx}
\usepackage{epstopdf}
\usepackage{tabularx}
\usepackage{rotating}
\usepackage{multirow}
\usepackage{color, colortbl}
\usepackage{algorithmic}
\usepackage{subfloat}
\usepackage{subfig}

\newcommand{\htemp}{{\bf h}}
\newcommand{\utemp}{{\bf u}}
\newcommand{\vtemp}{{\bf v}}

\newcommand{\qga}{QG_{33 \times 33}}
\newcommand{\qgb}{QG_{65 \times 65}}
\newcommand{\qgc}{QG_{129 \times 129}}

\newcommand{\NensK}{N_k}
\newcommand{\Nens}{N} 
\newcommand{\Nobs}{m} 
\newcommand{\Nstate}{n} 
\newcommand{\X}{{\bf X}} 
\newcommand{\x}{{\bf x}} 
\newcommand{\J}{\mathcal{J}} 
\newcommand{\DU}{\mathcal{D}} 
\newcommand{\lp}{\left (} 
\newcommand{\rp}{\right )} 
\newcommand{\lb}{\left [} 
\newcommand{\rb}{\right ]} 
\renewcommand{\ln}{\left \|} 
\newcommand{\rn}{\right \|}
\newcommand{\B}{{\bf B}} 
\newcommand{\R}{{\bf R}} 
\newcommand{\N}{M} 
\newcommand{\y}{{\bf y}} 
\renewcommand{\H}{{\bf H}} 
\newcommand{\xm}{{\overline{\bf x}}} 
\newcommand{\W}{{\boldsymbol \alpha}} 
\newcommand{\I}{{\bf I}} 
\newcommand{\M}{\mathcal{M}} 
\newcommand{\Nor}{\mathcal{N}} 
\newcommand{\xt}{{\bf x}^{\rm true}} 
\newcommand{\Ho}{{\mathcal{H}}} 
\newcommand{\lle}{\left \{ } 
\newcommand{\rle}{\right \}} 

\newcommand{\Q}{{\bf Q}} 
\renewcommand{\d}{{\bf d}} 

\newcommand{\BO}[1]{\mathcal{O}\lp #1\rp} 
\newcommand{\C}{{\bf C}} 
\renewcommand{\P}{{\bf P}} 
\renewcommand{\Re}{\mathbbm{R}}
\renewcommand{\S}{{\bf S}} 
\newcommand{\errobs}{{\boldsymbol \epsilon}} 
\newcommand{\errbac}{{\boldsymbol \xi}} 
\newcommand{\Prob}{{\mathcal{P}}} 

\newcommand{\D}{{\bf D}}

\newcommand{\WPOD}{{\boldsymbol \beta}}

\newcommand{\Z}{{\bf Z}} 
\newcolumntype{N}{>{\centering\arraybackslash} m{0.30\textwidth} }
\newcolumntype{V}{>{\centering\arraybackslash} m{0.02\textwidth} }

\newcommand{\basis}{{\boldsymbol \Psi}}

\newcommand{\Pobs}{p}

\newcommand{\zero}{{\bf 0}}

\newcommand{\BEST}{\widehat{\bf B}}

\newcommand{\expect}{\mathbb{E}} 
\newcommand{\Sig}{\widehat{\boldsymbol \Sigma}}
\newcommand{\Eigen}{{\boldsymbol \Sigma}}

\newcommand{\target}{{\bf T}}
\newcommand{\CEST}{{\bf \widehat{C}}}

\newcommand{\EX}[1]{\widetilde{#1}}

\newcommand{\WS}{{\boldsymbol \lambda}}
\newcommand{\U}{{\bf U}}
\newcommand{\V}{{\bf V}}

\newcommand{\E}{{\bf E}}
\newcommand{\GA}{{\boldsymbol \Gamma}}
\newcommand{\PI}{{\boldsymbol \Pi}}

\newcommand{\GammaB}{{\boldsymbol{\Gamma}}}

\newcommand{\ones}{{\bf 1}}

\makeindex

\begin{document}

\thispagestyle{empty}
\setcounter{page}{0}

\begin{Huge}
\begin{center}
Computational Science Laboratory Technical Report CSL-TR-16-2015\\
\today
\end{center}
\end{Huge}
\vfil
\begin{huge}
\begin{center}
Elias D. Nino and Adrian Sandu
\end{center}
\end{huge}

\vfil
\begin{huge}
\begin{it}
\begin{center}
``Ensemble Kalman Filter Implementations Based on Covariance Matrix Estimation''
\end{center}
\end{it}
\end{huge}
\vfil
\vfil

\begin{large}
\begin{center}
Computational Science Laboratory \\
Computer Science Department \\
Virginia Polytechnic Institute and State University \\
Blacksburg, VA 24060 \\
Phone: (540)-231-2193 \\
Fax: (540)-231-6075 \\ 
Email: {\rm enino@vt.edu},{\rm sandu@cs.vt.edu} \\
Web: {\rm http://csl.cs.vt.edu}
\end{center}
\end{large}

\vspace*{1cm}

\begin{center}
\begin{tabular}{c}
\includegraphics[width=1.in]{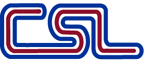}\\ 
.  \\
\includegraphics[width=1.in]{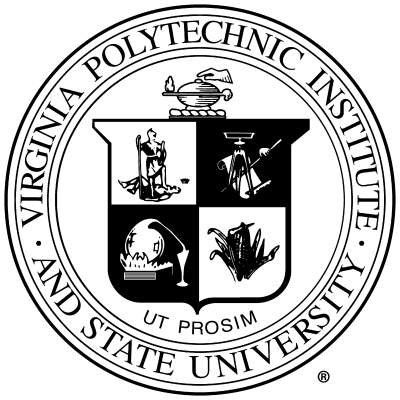} \\
\end{tabular}
\end{center}

\newpage

\maketitle
\tableofcontents{}

\begin{abstract}

This paper develops efficient ensemble Kalman filter (EnKF) implementations based on shrinkage covariance estimation. The forecast ensemble members at each step are used to estimate the background error covariance matrix via the Rao-Blackwell Ledoit and Wolf estimator, which has been developed specifically developed to approximate high-dimensional covariance matrices using a small number of samples. Additional samples are taken from the normal distribution described by the background ensemble mean and the estimated background covariance matrix in order to increase the size of the ensemble and reduce the sampling error of the filter. This increase in the size of the ensemble is obtained without running the forward model. After the assimilation step, the additional samples are discarded and only the initial members are propagated. Two implementations are considered. In the EnKF Full-Space (EnKF-FS) approach the assimilation process is performed in the model space, while the EnKF Reduce-Space (EnKF-RS) formulation performs the analysis in the subspace spanned by the ensemble members. Numerical experiments carried out with a quasi-geostrophic model show that the proposed implementations outperform current methods such as the traditional EnKF formulation, square root filters, and inflation-free EnKF implementations. 
The proposed implementations provide good results with small ensemble sizes ($\sim 10$) and small percentages of observed components from the vector state. These results are similar (and in some cases better) to traditional methods using large ensemble sizes ($\sim 80$) and large percentages of observed components. 
The computational times of the new implementations remain reasonably low.
\end{abstract}


\section{Introduction}\label{sec:introduction}
%
Sequential data assimilation estimates the current unknown state $\xt \in \Re^{\Nstate \times 1}$ of a physical system as follows. The background (prior) state $\x^b \in \Re^{\Nstate \times 1}$ is given by a physical model initialized
with the best estimate of the state at a previous time:
\begin{eqnarray}
\label{eq:model-operator}
\x^b \equiv \x_{\rm current} = \M_{t_{\rm previous}\rightarrow t_{\rm current}} \lp \x_{\rm previous} \rp \,.
\end{eqnarray}
The background errors are assumed to be unbiased and normally distributed
\begin{eqnarray}
\label{eq:prior-estimate}
\errbac = \x^b - \xt \in  \Nor \lp \zero_{\Nstate} ,\, \B\rp
\end{eqnarray}
where $\Nstate$ is the dimension of the model state, $\zero_{q}$ is the vector in the $q$-dimensional space whose components are all zeros, and $\B \in \Re^{\Nstate \times \Nstate}$ is the unknown background error covariance matrix. 
Observations of the true state are taken
\begin{eqnarray}
\label{eq:observations}
\displaystyle 
\y = \Ho \lp \xt \rp + \errobs \in \Re^{\Nobs \times 1} \,,
\end{eqnarray}
where $\Ho: \Re^{\Nstate \times 1} \rightarrow \Re^{\Nobs \times 1}$ is the observation operator, and $\Nobs$ is the number of observed components.
The observational errors are assumed to be normally distributed $\errobs \sim \Nor \lp \zero_{\Nobs},\, \R \rp$,  where $\R \in \Re^{\Nobs \times \Nobs}$ is the observation error covariance matrix that  is assumed to possess a simple structure (e.g., block diagonal) and therefore its inverse can be easily computed.


Under the Gaussian assumption on data and model errors the negative logarithms of the a posteriori probability density is the 3D-Var cost function \cite{DVAR3}:
\begin{eqnarray}
\label{eq:cost-function-3DVAR}
\displaystyle 
\J \lp \x \rp =  \frac{1}{2}\cdot \ln \x-\x^b \rn_{\B^{-1}}^2 + \frac{1}{2} \cdot \ln \y-\Ho \lp \x \rp\rn_{\R^{-1}}^{2}\,.
\end{eqnarray}
The maximum likelihood estimate of the state is obtained by minimizing the cost function \eqref{eq:cost-function-3DVAR}, i.e., the analysis state $\x^a \in \Re^{\Nstate \times 1}$ is the solution of the following optimization problem:
\begin{eqnarray}
\label{eq:solution-model-space}
\displaystyle
\x^a = \underset{\x}{\arg \min} \,\, \J \lp \x \rp \,.
\end{eqnarray}
The solution of \eqref{eq:solution-model-space} over the subspace spanned by the ensemble members is:
\begin{eqnarray}
\label{eq:analysis-state-general}
\x^a = \x^b+\basis \cdot \W_{\B} \,,
\end{eqnarray}
where $\basis = \B^{1/2}$ is a set of basis vectors satisfying $\basis\, \basis^T = \B$, and the vector of weights is given by 
\begin{eqnarray}
\label{eq:optimal-weights-general}
\displaystyle 
\W_{\B} = \V_{\B}^T \cdot \lp \R + \V_{\B} \cdot \V_{\B}^T \rp^{-1} \cdot \lp \y - \Ho (\x^b) \rp \,,
\end{eqnarray}
where $\V_{\B} = \H \cdot \basis \in \Re^{\Nobs \times \Nstate}$ and $\H = \Ho'$. By propagating in time the analysis \eqref{eq:analysis-state-general} a new prior (background) state is obtained for the next assimilation cycle.

Since $\B$ is unknown the direct use of \eqref{eq:analysis-state-general} is infeasible. Methods for estimating the background error covariance matrix have been developed \cite{Furrer2007227,Cheng2010}. However, in the absence of prior information about the true structure of $\B$ biased estimators are often obtained. Buhener \cite{Buhener2005}  discusses the impact of having biased estimators for the background errors in the assimilation process. AIn the context of sequential data assimilation an ensemble of model realizations \cite{Evensen:2006:DAE:1206873} is built in order to represent the background error statistics:
\begin{eqnarray}
\label{eq:initial-ensemble}
\displaystyle \X^b = \lb \x^b_1,\, \x^b_2,\, \ldots \,, \x^b_{\Nens} \rb \in \Re^{\Nstate \times \Nens}\,,
\end{eqnarray}
where $\x^b_i \in \Re^{\Nstate \times 1}$ is the i-th ensemble member (model run) in the \textit{background stage}. Estimates of $\x^b$ and $\B$ are obtained via the empirical moments of the ensemble \eqref{eq:initial-ensemble}
\begin{eqnarray}
\label{eq:ensemble-mean}
\displaystyle 
\x^b \approx \xm^b = \frac{1}{\Nens} \cdot \sum_{i=1}^{\Nens} \x^b_{i} \in \Re^{\Nstate \times 1} \,,
\end{eqnarray}
and
\begin{eqnarray}
\label{eq:ensemble-covariance-matrix}
\displaystyle
\B \approx \P^b = \S \cdot \S^T \in \Re^{\Nstate \times \Nstate} \,,
\end{eqnarray}
where the matrix of scaled member deviations is
\begin{eqnarray}
\label{eq:matrix-of-member-deviations}
\displaystyle
\S = \frac{1}{\sqrt{\Nens-1}} \cdot \lb \X^b - \xm^b \otimes \ones_{\Nens}^T \rb \in \Re^{\Nstate \times \Nens},
\end{eqnarray}
where $\ones_{q}$ is the vector in the $q$-th dimensional space whose components are all ones.
The set of basis vectors \eqref{eq:matrix-of-member-deviations} does not span the full space of model errors. By replacing the estimators \eqref{eq:ensemble-mean}  and \eqref{eq:ensemble-covariance-matrix} in \eqref{eq:analysis-state-general}, the analysis can be approximated as follows:
\begin{eqnarray}
\label{eq:analysis-estimators-general}
\displaystyle
\x^a \approx \xm^a = \xm^b + \S \cdot \W \in \Re^{\Nstate \times 1} \,,
\end{eqnarray}
where 
\begin{eqnarray}
\label{eq:optimal-weights-basis-approximation}
\displaystyle
\W = \V \cdot \lp \R + \V \cdot {\V}^T \rp^{-1} \cdot \lp \y - \Ho(\xm^b) \rp \in \Re^{\Nens \times 1} \,,
\end{eqnarray}
and $\V = \H \cdot \S \in \Re^{\Nobs \times \Nens}$. The analysis state \eqref{eq:analysis-estimators-general} is the optimal solution in the space spanned by the ensemble members \eqref{eq:initial-ensemble},
but in general it is not optimal in the model space. Since $\Nstate \gg \Nens$ information contained in \eqref{eq:optimal-weights-general} is not represented by the set of basis vectors $\S$. From another point of view, the ensemble covariance matrix $\P^b$ is rank deficient and therefore there are insufficient degrees of freedom to explain the full error. This problem can be alleviate by making use of localization techniques \cite{JGRD:JGRD50654,citeulike:6937506,TELA:TELA460}. However, since the structure of $\B$ remains unknown, the use of localization generally increases the bias in the background error estimation. 

There is an opportunity to avoid the intrinsic need of inflation in ensemble based methods by replacing the covariance matrix \eqref{eq:ensemble-covariance-matrix} with a more accurate and well-conditioned estimate of $\B$. We do not want to impose any kind of structure on $\B$ since it will make our approach sensitive to problems faced by current implementations. Instead, we seek to capture most of the information contained in the ensemble. Shrinkage covariance estimators developed to estimate high-dimensional covariance matrices from a small number of samples \cite{Ledoit2004365} fit very well in the context of sequential data assimilation.

The remaining part of the paper is organized as follows. Section \ref{sec:preliminaries} reviews ensemble based data assimilation and shrinkage covariance estimation. In section \ref{sec:proposed-implementation} the two novel implementations of the ensemble Kalman filter based on shrinkage covariance estimation are proposed. Experimental results making use of a quasi-geostrophic model are given in section \ref{sec:experimental-results}. Section \ref{sec:conclusions} summarizes the conclusions of this work.

\section{Background}
\label{sec:preliminaries}

In this section we review relevant concepts with regard to shrinkage covariance estimation and ensemble based methods in sequential data assimilation.

\subsection{Covariance estimation}
\label{subsec:covariance-estimation}

Many problems in science and engineering require an estimate of a covariance matrix and/or its inverse, where the matrix dimension $\Nstate$ is large compared to the sample size $\Nens$. Different applications ranging from variational \cite{citeulike:6937506,citeulike:975315} to sequential \cite{JGRD:JGRD50654,citeulike:5092716} data assimilation rely on accurately estimated covariance matrices. 

Let $\{s_1, \, s_2, \, \ldots, \, s_{\Nens}\}$ be a sample of independent identical distributed $\Nstate$-dimensional Gaussian vectors 
\begin{eqnarray*}
\displaystyle 
s_i \sim \Nor( \zero_{\Nstate} ,\, \Q) \in \Re^{\Nstate \times 1}.
\end{eqnarray*}
A common approach is to estimate $\Q \in \Re^{\Nstate \times \Nstate}$ by the sample covariance matrix $\C_s$
\begin{eqnarray}
\label{eq:empirical-moments-sample-covariance}
\C_s &=& \frac{1}{\Nens-1} \cdot \sum_{i=1}^{\Nens} s_i \otimes s_i^T \in \Re^{\Nstate \times \Nstate} \,.
\end{eqnarray}
$\C_s$ is the maximum likelihood estimator when it is invertible \cite{Ledoit2004365}. However, under the condition $\Nstate \gg \Nens$, this is not the case. The simpler thing to do in order to deal with the rank-deficiency of $\C_s$ is to impose some  structure (i.e., localization in ensemble based methods). However, in the absence of prior information about the true structure of $\Q$, $\C_s$ will poorly describe the correlations between different components of the samples $\{s_i\}_{1 \le i \le \Nens}$.
In order to improve estimation of covariance matrices many methods have been proposed in the literature  based on tapering procedures \cite{cai2010,Xia6252067}, minimizing the log-determinant divergence  \cite{ravikumar2011}, and greedy methods \cite{GreedyEstimation12}. Another class of well-conditioned estimators is based on shrinkage approximations \cite{Farebrother1978,DeMiguel20133018,Couillet201499,Park20141,Bickel2014204,Fisher20111909}. These approximations express the estimated covariance matrix as a weighted average of some target matrix $\target \in \Re^{\Nstate \times \Nstate}$ and the empirical covariance matrix \eqref{eq:empirical-moments-sample-covariance}. To better understand this assume that the components of $s_i$ are uncorrelated for $1 \le i \le \Nens$. A simple estimate of $\Q$ is given by
\begin{eqnarray*}
\target = \frac{{\rm tr} \lp \C_s \rp}{\Nstate} \cdot \I_{\Nstate \times \Nstate} \,,
\end{eqnarray*}
where $\I_{\Nstate \times \Nstate}$ is the identity matrix in the $\Nstate$-dimensional space. Note that this structure will reduce the variance but will increase the bias when the diagonal assumption is not fulfilled. A reasonable trade-off is achieved by the shrinkage of $\C_s$ towards $\target$ and provides the followng class of estimators
\begin{eqnarray}
\label{eq:shinkrage-estimator}
\displaystyle
\CEST = \gamma \cdot \target + \lp 1 - \gamma \rp \cdot \C_s  \in \Re^{\Nstate \times \Nstate}\,,
\end{eqnarray}
where $\gamma \in [0,\,1]$. The problem is then reduced to find an optimal value for $\gamma$ in which the squared loss
\begin{eqnarray}
\label{eq:squared-loss}
\displaystyle
\expect \lb \ln \CEST - \Q \rn_{F}^2 \rb
\end{eqnarray}
is minimized, where $\ln \bullet \rn_{F}$ denotes the Frobenius norm. There are many shrinkage based estimators derived from the minimization of \eqref{eq:squared-loss} subject to \eqref{eq:shinkrage-estimator}. We restrict our exploration to three well-accepted methods: the Ledoit and Wolf estimator \cite{Ledoit2004365}, the Rao-Blackwell Ledoit and Wolf estimator \cite{Yilun5484583} and the oracle approximating shrinkage estimator \cite{Yilun5484583,Yilun5743027}.

The distribution-free Ledoit and Wolf (LW) estimator \cite{Ledoit2004365} has been proven more accurate than the sample covariance matrix and some estimators proposed in finite sample decision theory. Moreover, it is better conditioned than the true covariance matrix \cite{Ledoit2004365}. The optimal $\gamma$ value proposed by this estimator is 
\begin{eqnarray}
\label{eq:LW-estimator}
\displaystyle
\gamma_{LW} = \min \lp \frac{\sum_{i=1}^{\Nens} \ln \C_s - s_i \otimes s_i^T \rn_{F}^2}{\Nens^2 \cdot \lb {\rm tr} \lp \C_s^2 \rp - \frac{{\rm tr}^2 \lp \C_s \rp}{\Nstate} \rb} ,\, 1\rp
\end{eqnarray}
and the LW estimator $\CEST_{LW}$ is obtained by using $\gamma_{LW}$ in \eqref{eq:shinkrage-estimator}.

The Rao-Blackwell Ledoit and Wolf (RBLW) estimator \cite{Yilun5743027,Yilun5484583}  provably improves the LW method under Gaussian assumptions. The motivation of this estimator is  that, under Gaussian assumptions, all the information required in order to get a well-conditioned estimate of $\Q$ is contained in $\C_s$. The proposed value for $\gamma$ is
\begin{eqnarray}
\label{eq:RBLW-estimator}
\displaystyle 
\gamma_{RBLW} = \min \lp \frac{\frac{\Nens-2}{\Nstate} \cdot {\rm tr}\lp \C_s^2\rp + {\rm tr}^2\lp \C_s \rp }{(\Nens+2) \cdot \lb {\rm tr}\lp \C_s^2\rp - \frac{{\rm tr}^2 \lp \C_s\rp}{\Nstate}\rb},\, 1\rp
\end{eqnarray}
and the corresponding estimator $\CEST_{RBLW}$ is obtained by replacing \eqref{eq:RBLW-estimator} in the equation \eqref{eq:shinkrage-estimator}. In addition, in \cite[Theorem 2]{Yilun5743027}, it is proven that
\begin{eqnarray*}
\displaystyle \expect \lb \ln \CEST_{RBLW} - \Q \rn_{F}^2  \rb \le  \expect \lb \ln \CEST_{LW} - \Q \rn_{F}^2  \rb \,,
\end{eqnarray*}
which rigorously shows the RBLW estimator to be a better approximation of $\Q$ than the LW estimator under the Gaussian assumption.

The oracle approximating shrinkage (OAS) estimator \cite{Yilun5484583} is an iterative approximation of the unimplementable oracle method \cite[Section 3]{Yilun5484583}. The optimal $\gamma$ at each iteration $j$ is given by
\begin{subequations}
\begin{eqnarray}
\displaystyle
\label{eq:OAS-estimator-matrix-j}
\CEST_j &=& \gamma_j \cdot \target + (1-\gamma_j) \cdot \C_s \,,\\
\label{eq:OAS-estimator-gamma-value}
\gamma_{j+1} &=&  \frac{\lp 1-\frac{2}{\Nstate}\rp \cdot {\rm tr}\lp \CEST_j \cdot \C_s \rp + {\rm tr}^2 \lp \CEST_j \rp}{\lp \Nens+1-\frac{2}{\Nstate} \rp \cdot {\rm tr} \lp \CEST_j \cdot \C_s \rp + \lp 1-\frac{\Nens}{\Nstate} \rp \cdot {\rm tr}^2 \lp \CEST_j \rp}
\end{eqnarray}
\end{subequations}
where the initial estimator $\CEST_0$ can be any estimator (i.e., the LW, the RBLW, or even the sample covariance matrix).

All the estimators presented in this section provide well-conditioned approximations to the unknown covariance matrix $\Q$. We center our attention on the RBLW estimator since in high dimensional problems, such those found in data assimilation, this estimator can be implemented easily, and under the Gaussian assumption it provides better approximations than the LW estimator.

\subsection{Sequential data assimilation methods}
\label{subsec:data-assimilation-methods}

Ensemble data assimilation methods are widely used in applications to weather, oceanography, and climatology \cite{DAIssues}. These methods represent the background error statistics by the empirical moments of the ensemble \eqref{eq:ensemble-mean} and \eqref{eq:ensemble-covariance-matrix}\cite{EnKF1657419}. The trajectory of each ensemble member and the dispersion of the overall ensemble around the background state provide meaningful information about the background error distribution. One of the most important advantages of ensemble DA is the flow-dependent forecast error covariance matrix \cite{Buhener2005}. When observations are available the analysis state \eqref{eq:analysis-estimators-general} is computed. The next step is to generate an ensemble which describes the analysis uncertainty around this optimal state. We briefly discuss three  implementations that achieve this, the ensemble square root filter, the ensemble transform Kalman filter, and the (basic) ensemble Kalman filter. 

In ensemble square root filters \cite{Tippett2003}, the new analysis ensemble in the analysis is built as follows:
\begin{eqnarray}
\label{eq:ensrf-analysis}
\displaystyle
\X^a = \xm^a \otimes \ones_{\Nens}^T + \S \cdot \lb \I_{\Nens \times \Nens} - {\V}^T \cdot \lp \R + \V \cdot \V^T \rp^{-1} \cdot \V \rb^{1/2} \in \Re^{\Nstate \times \Nens} \,.
\end{eqnarray}
As expected the analysis ensemble members live in the subspace space spanned by the background ensemble, this is, the space spanned by the columns of $\S$. All possible information that can be obtained from the model states is contained in this set of basis vectors.

The covariance matrix in the observation space is
\begin{eqnarray}
\label{eq:weighted-covariance-matrix}
{\bf W}_{\rm obs} = \R+\V \cdot \V^T \in \Re^{\Nobs \times \Nobs} \,,
\end{eqnarray}
the linear system ${\bf W}_{\rm obs} \cdot \Z_{\bf V} = \V \in \Re^{\Nobs \times \Nens}$ can be solved via the iterative Sherman Morrison formula (ISMF) \cite{Nino2014}
\begin{subequations}
\label{eq:ISMF}
\begin{eqnarray}
\displaystyle
\label{eq:ISMF-stage-1}
\htemp^{(k)} &=& \lp 1+ \vtemp_k^T \cdot \utemp_{k}^{(k-1)}\rp^{-1} \cdot \utemp_{k}^{(k-1)} \in \Re^{\Nobs \times 1} \\
\label{eq:ISMF-stage-2}
\Z^{(k)} &=& \Z^{(k-1)} - \htemp^{(k)} \cdot \lb \vtemp_k^T \cdot \Z^{(k-1)} \rb \in \Re^{\Nobs \times \Nens} \\
\label{eq:ISMF-stage-3}
\U^{(k)} &=& \U^{(k-1)} - \htemp^{(k)} \cdot \lb \vtemp_k^T \cdot \U^{(k-1)} \rb \in \Re^{\Nobs \times \Nens}
\end{eqnarray}
\end{subequations}
for $1 \le k \le \Nens$,  where $\Z^{(0)} = \R^{-1} \cdot \lp \y - \H \cdot \xm^b \rp$ and $\U^{(0)} = \R^{-1} \cdot \V$. 
The implementation \eqref{eq:ISMF} requires no more than $\BO{\Nens^2 \cdot \Nobs}$ long computations, where $\vtemp_k$ and $\utemp_k^{(k-1)}$ are the $k$-th column of the matrices $\V$ and $\U^{(k-1)}_k$, respectively. Moreover, by applying the singular value decomposition to
\begin{eqnarray}
\label{eq:square-root-filter}
\displaystyle
\V^T \cdot \Z_{\V} = \U_{\Z} \cdot \Eigen_{\Z} \cdot \U_{\Z}^T \in \Re^{\Nens \times \Nens} \in \Re^{\Nens \times \Nens}\,,
\end{eqnarray}
the square root in \eqref{eq:ensrf-analysis} 
\begin{eqnarray*}
\GammaB = \U_{\Z} \cdot \lb \I_{\Nens \times \Nens} - \Eigen_{\Z}\rb^{1/2} \cdot \U_{\Z}^T \,,
\end{eqnarray*}
can be efficiently computed with no more than $\BO{\Nens^3}$ long computations. 

Another widely used square root filter implementation is the ensemble transform Kalman filter (EnTKF) \cite{Nuru2013}. By making use of the matrix identity
\begin{eqnarray}
\label{eq:EnTKF-matrix-identity}
\displaystyle \I - \V \cdot \lb \R + \V \cdot \V^T \rb^{-1} \cdot \V^T = \lb \I+\V \cdot \R \cdot \V^T \rb^{-1} \in \Re^{\Nobs \times \Nobs} \,,
\end{eqnarray} 
and the singular value decomposition
\begin{eqnarray}
\label{eq:entkf-svd}
\displaystyle
\V^T = \U_{\V} \cdot \Sig_{\V} \cdot \V_{\V}^T \in \Re^{\Nens \times \Nobs} \,,
\end{eqnarray}
the analysis state \eqref{eq:analysis-estimators-general} can be written as follows
\begin{eqnarray}
\label{eq:entkf-analysis}
\displaystyle
\x^a = \xm^b + \S \cdot \WPOD
\end{eqnarray}
with the optimal weights $\WPOD \in \Re^{\Nens \times \Nens}$ given by
\begin{eqnarray}
\label{eq:entkf-weights}
\displaystyle 
\WPOD = \U_{\V} \cdot \Sig_{\V} \cdot \lp \I_{\Nens \times \Nens}+\Sig_{\V}^T \cdot \Sig_{\V} \rp^{-1} \cdot \V_{\V}^T \cdot \sqrt{\R} \cdot \lp \y - \H \cdot \xm^b \rp.
\end{eqnarray}
The new ensemble members are built as follows
\begin{eqnarray}
\label{eq:entkf-new-samples}
\displaystyle
\X^a = \xm^a \otimes \ones_{\Nens}^T  + \S \cdot \U_{\V} \cdot  \lb \I_{\Nens \times \Nens}+\Sig_{\V} \cdot \Sig_{\V}^T\rb^{1/2} \cdot \U_{\V}^T.
\end{eqnarray}

In the ensemble Kalman filter (EnKF) \cite{EnKFEvensen} each ensemble member  \eqref{eq:initial-ensemble} and the observations \eqref{eq:observations} are treated as random variables, and the  $i$-th ensemble member is updated as follows:
\begin{eqnarray}
\label{eq:analysis-ensemble-members}
\displaystyle
\x^a_i = \x^b_i + \S \cdot \WPOD_i \in \Re^{\Nstate \times 1} \,, \quad i = 1,\dots,\Nens\,,
\end{eqnarray}
where the optimal weights for $i$-th ensemble member are given by
\begin{eqnarray}
\label{eq:optimal-weights-ith-member}
\WPOD_i = \cdot \V^T \cdot \lp \R+\V \cdot \V^T \rp^{-1} \cdot \lb \y^s_i - \H \cdot \x^b_i\rb \in \Re^{\Nens \times 1}
\end{eqnarray}
and
\begin{eqnarray}
\label{eq:observation-sample}
\displaystyle
\y^s_i \sim \Nor \lp \y,\, \R \rp \,,
\end{eqnarray}
The addition of the perturbed observations \eqref{eq:observation-sample} in the analysis provides asymptotically correct analysis-error covariance estimates for large ensemble sizes and makes the formulation of the EnKF statistical consistent \cite{Thomas2002}. However, it also has been proven that the inclusion of perturbed observations introduces sampling errors in the assimilation \cite{SamplingErrors1,SamplingErrors2}.

One of the important problems faced by current ensemble based methods is filter divergence due to the insufficient degrees of freedom ($\Nens \ll \Nstate$). To alleviate this deficiency localization is used to impose  structure on the sample covariance matrix \eqref{eq:ensemble-covariance-matrix} according to the physics of the model. Intuitively, the correlations between individual model variable errors decays with distance, e.g., exponentially:
\begin{eqnarray}
\label{eq:localization}
{\boldsymbol \rho}_{ij} =  \exp \lp \frac{-d \lp i,j\rp^2}{2 \cdot L^2}\rp \,,
\end{eqnarray}
where $d \lp i,j\rp$ is the physical distance between the locations of the $i$-th and $j$-th model components, and $L$ is the localization radius. The localized covariance matrix is obtained as
\begin{eqnarray}
\label{eq:forcing-to-have-structure}
\displaystyle
\P^b_{\rm Loc} = {\boldsymbol \rho} \circ \P^b \in \Re^{\Nstate \times \Nstate}.
\end{eqnarray}
where $\circ$ is the Schur product. This method can be impractical since it relies on the explicit computation of $\P^b$. Moreover, there are no guaranties that this method captures the true structure of $\B$. More sophisticated covariance estimation methods have been proposed in the context of data assimilation. A classic approximation is the Hollingworth and Lonnberg method \cite{TELA:TELA460} in which the difference between observations and background states are treated as a combination of background and observations errors. However, this method provides statistics of background errors in observation spaces, requires uniform observing network (not the case in practice), and the resulting statistics are biased towards data-dense areas. Another method has been proposed by Benedetti and Fisher \cite{QJ:QJ37}  based on forecast differences in which the spatial correlations of background errors are assumed to be similar at 24 and 48 hours forecasts. This method can be efficiently implemented in practice, however, it does not perform well in data-sparse regions, and  the statistics provided are a mixture of analysis and background errors. Since the structure of $\B$ remains unknown, assumptions made about its structure may increase the bias in the estimate. Furthermore, the balance of variables in the model with some physical meaning can be disturbed when localization is utilized \cite{TELA:TELA076}. 

A different approach is based on the 3D-Var cost function in the ensemble space. Any vector $\x \in \Re^{\Nstate \times 1}$ in the ensemble subspace can be written as
\begin{eqnarray}
\label{eq:vector-in-ensemble-space}
\displaystyle
\x = \xm^b + \U \cdot \W 
\end{eqnarray}
where $\U$ is the matrix of anomalies
\begin{eqnarray}
\label{eq:matrix-of-anomalies}
\displaystyle
\U = \lb \x^b_1 - \xm^b,\, \x^b_2 - \xm^b, \, \ldots, \, \x_{\Nens}^b - \xm^b  \rb \in \Re^{\Nstate \times \Nens} \,,
\end{eqnarray}
and $\W \in \Re^{\Nens \times 1}$ is a vector  to be determined. The columns of $\U$ and $\S$ span the same space. Using \eqref{eq:vector-in-ensemble-space} the 3D-Var cost function \eqref{eq:cost-function-3DVAR} in the ensemble space reads
\begin{eqnarray}
\label{eq:3Dvar-ensemble-space}
\displaystyle
\J_{\rm ens} \lp \W \rp = \frac{1}{2} \cdot \ln \U \cdot \W \rn_{\B^{-1}}^2 + \frac{1}{2} \cdot \ln \d - \Q \cdot \W \rn_{\R^{-1}}^2 \,,
\end{eqnarray}
where the optimal value of the control variable $\W$ 
\begin{eqnarray}
\label{eq:solution-ensemble-space}
\displaystyle
\W^{*} = \underset{\W}{\arg \min} \, \J_{\rm ens} \lp \W \rp  \,,
\end{eqnarray}
provides the analysis state in \eqref{eq:vector-in-ensemble-space}
\begin{eqnarray}
\label{eq:analysis-state-in-the-ensemble-space}
\xm^a = \xm^b + \U \cdot \W^{*} \,.
\end{eqnarray}
Two recent formulations based on this approximation are the finite size anddual ensemble Kalman filters \cite{Bocquet2012}. These formulation avoid the intrinsic needed of inflation by choosing Jeffrey's prior for background errors:
\begin{eqnarray*}
\Prob \lp \x^b , \B \rp = \Prob_J \lp  \x^b \rp \cdot \Prob_J \lp \B \rp \,,
\end{eqnarray*}
where the parameters $\x^b$ and $\B$ are assumed to be independent. 

In the case of the finite size ensemble Kalman filter (EnKF-N) the cost function in the ensemble space reads
\begin{eqnarray}
\label{eq:cost-function-EnKF-N}
\displaystyle
\J_{\rm ens}^\textsc{fn} \lp \W \rp = \frac{1}{2} \ln \y - \Ho \lp \xm^b + \U \cdot \W \rp \rn_{\R^{-1}}^2 + \frac{\Nens}{2} \cdot \log \lp 1+\frac{1}{\Nens} + \ln \W \rn^2 \rp.
\end{eqnarray}
Minimization of this cost function provides the optimal weights in the ensemble space
\begin{eqnarray}
\label{eq:optimization-problem-EnKF-FN}
\displaystyle
\W^{*} = \underset{\W}{\arg \min} \, \J_{\rm ens}^\textsc{fn} \lp \W \rp  \,,
\end{eqnarray}
and the analysis is computed via \eqref{eq:analysis-state-in-the-ensemble-space}. The projection of the analysis covariance matrix on the ensemble space is approximated by the inverse of the Hessian of \eqref{eq:cost-function-EnKF-N} at the optimal value \eqref{eq:optimization-problem-EnKF-FN}. The Hessian reads:
\begin{eqnarray}
\label{eq:Hessian-EnKF-FN}
\nabla^2_{\W,\W} \J_{\rm ens}^\textsc{fn} \lp \W \rp  &=& \lb \H \cdot \U \rb^T \cdot \R^{-1}  \cdot \H \cdot \U \\ 
\nonumber
&+&  \Nens \cdot \frac{\lp 1 + \frac{1}{\Nens} + \ln \W\rn^2 \rp \cdot \I_{\Nens \times \Nens} - 2 \cdot \W \cdot \W^T }{\lp 1+\frac{1}{\Nens} + \ln \W \rn^2 \rp^2} \in \Re^{\Nens \times \Nens}.
\end{eqnarray}
The analysis ensemble is generated as follows:
\begin{eqnarray}
\label{eq:analysis-EnKF-N}
\displaystyle
\X^a = \xm^a \otimes \ones_{\Nens}^T + \U \cdot \lle \lp \Nens-1\rp \cdot \lb \nabla^2_{\W,\W} \J_{\rm ens}^\textsc{fn} \lp \W^{*} \rp  \rb^{-1}  \rle^{1/2} \cdot {\boldsymbol \Phi} \in \Re^{\Nstate \times \Nens} \,,
\end{eqnarray}
where $ {\boldsymbol \Phi}\in \Re^{\Nens \times \Nens}$ is an arbitrary orthogonal matrix which preserves the ensemble mean \eqref{eq:ensemble-mean}. 

Another approach is based on the dual formulation of the cost function \eqref{eq:Hessian-EnKF-FN} \cite{Bocquet2012} %
\begin{eqnarray}
\label{eq:cost-function-EnKF-DU}
\displaystyle 
\DU_{\rm ens}^\textsc{du} \lp \zeta \rp = \ln \y - \H \cdot \x^b \rn_{{{\bf W}_{\rm \zeta}}^{-1}}^2 + \zeta \cdot \lp 1 + \frac{1}{\Nens} \rp + \Nens \cdot \log \lp \frac{\Nens}{\zeta}\rp -\Nens \,,
\end{eqnarray}
where the weighted covariance matrix ${\bf W}_{\rm \zeta}$ reads
\begin{eqnarray}
\label{eq:weighted-covariance-W-Dual}
\displaystyle
{\bf W}_{\rm \zeta} = \R+\frac{1}{\zeta} \cdot \V \cdot \V^T \in \Re^{\Nobs \times \Nobs}.
\end{eqnarray}
The dual optimization problem is one-dimansional
\begin{eqnarray}
\label{eq:optimization-problem-EnKF-DU}
\displaystyle
\zeta^{*} = \underset{ \zeta \in  \lp 0, \frac{\Nens}{1+\frac{1}{\Nens}} \rb}{\arg \min} \, \DU_{\rm ens}^\textsc{du} \lp \zeta \rp 
\qquad \textnormal{subject to} \quad \zeta \in  \lp 0, \frac{\Nens}{1+\frac{1}{\Nens}} \rb.
\end{eqnarray}
The optimal state is computed as follows:
\begin{eqnarray}
\label{eq:optimal-analysis}
\displaystyle 
\xm^a = \xm^b + \U \cdot \lb \V^T \cdot \R^{-1} \cdot \V + \zeta^{*} \cdot \I_{\Nens \times \Nens} \rb^{-1} \cdot \V^T \cdot \R^{-1} \cdot \lb \y - \H \cdot \xm^b \rb \,,
\end{eqnarray}
and the following analysis ensemble is built:
\begin{eqnarray}
\label{eq:analsysis-ensemble-EnKF-DU}
\displaystyle \X^a = \xm^a \otimes \ones_{\Nens}^T + \U \cdot \lle \lp \Nens-1 \rp \cdot \lb \V^T \cdot \R^{-1} \cdot \V + \zeta^{*} \cdot \I_{\Nens \times \Nens} \rb^{-1} \rle^{1/2} \cdot  {\boldsymbol \Phi}.
\end{eqnarray}

In this paper we consider a different representation of the background error statistics by making use of shrinkage covariance estimation. The idea is not to impose any structure on $\P^b$ but to obtain a well-conditioned estimator $\BEST$ of the background error covariance matrix $\B$ wherein using all the possible information brought from the ensemble members. Samples from the distribution $\Nor \lp \xm^b ,\, \BEST \rp$ are taken in order to better represent the error statistics and to increase the number of degrees of freedom. Two novel EnKF implementations based on the Rao-Blackwell Ledoit and Wolf estimator are presented in the next section. 

\section{Ensemble Filters Based on Shrinkage Covariance Estimators}
\label{sec:proposed-implementation}

In this section, we propose two efficient implementations of the EnKF based on the RBLW estimator \eqref{eq:RBLW-estimator}. As mentioned before, we do not impose any kind of structure on $\P^b$ since the information brought by the ensemble members is more than only background errors
\begin{eqnarray*}
\P^b = \B + \Q + \C \in \Re^{\Nstate \times \Nstate} \,,
\end{eqnarray*}
where $\Q$ is the covariance of model errors and $\C \in \Re^{\Nstate \times \Nstate}$ is the covariance matrix of additional errors whose sources are unknown for us. Errors coming from different sources are assumed to be uncorrelated.  We seek to exploit the information brought by ensemble members and use the RBLW covariance estimator \eqref{eq:RBLW-estimator} to build a covariance matrix that captures all error correlations. The standard form of this estimatordepends on the explicit representation of $\P^b$. The efficient implementation for high-dimensional covariance matrices presented in section \ref{subsec:RBLW-EnKF} avoids the explicit computation of $\P^b$. Section \ref{subsec:enkf-rblw} discusses two EnKF implementations based on the RBLW estimator. Section \ref{subsec:sampling-high-dimensions} develops an efficient sampling method in high dimensions for drawing samples from the prior error distribution based on the RBLW estimate. Finally, section \ref{subsec:comparison-of-the-EnKF-FS-RS} discusses the similarities and differences between the two proposed implementations.

\subsection{RBLW estimator for covariance matrices in high-dimensions}
\label{subsec:RBLW-EnKF}
Consider the sample covariance matrix \eqref{eq:ensemble-covariance-matrix}.
In the context of data assimilation the RBLW  estimator \eqref{eq:shinkrage-estimator},\eqref{eq:RBLW-estimator} reads
\begin{subequations}
\label{eq:EnkF-optimal-values}
\begin{eqnarray}
\label{eq:RBLW-EnKF}
\displaystyle
\BEST = \gamma_{\BEST} \cdot \mu_{\B} \cdot \I + (1 - \gamma_{\BEST}) \cdot \P^b \in \Re^{\Nstate \times \Nstate} \,,
\end{eqnarray}
where
\begin{eqnarray}
\label{eq:EnKF-mu-value}
\displaystyle 
\mu_{\BEST} &=& \frac{{\rm tr} \lp \P^b \rp}{\Nstate} \,, \\
\label{eq:EnKF-gamma-value}
\gamma_{\B} &=& \min \lp \frac{\frac{\Nens-2}{ \Nstate} \cdot {\rm tr} \lp \lb\P^b \rb^2 \rp + {\rm tr}^2 \lp  \P^b \rp}{\lp \Nens+2 \rp \cdot \lb {\rm tr} \lp \lb \P^b \rb^2 \rp - \frac{{\rm tr}^2 \lp \P^b \rp}{\Nstate} \rb} , \, 1 \rp \,.
\end{eqnarray}
\end{subequations}
Since the dimension of the model state is high ($\Nstate \sim \BO{10^7}$), the direct computation of  \eqref{eq:RBLW-estimator-EnKF} is impractical as it requires the explicit representation of the sample covariance matrix $\P^b$. An alternative manner to compute ${\rm tr} \lp \P^b \rp$ and ${\rm tr} \lp [\P^b]^2 \rp$ is proposed. Consider the eigenvalue decomposition of $\P^b$
\begin{eqnarray}
\label{eq:Pb-SVD}
\displaystyle
\P^b = \U_{\P^b} \cdot \Eigen_{\P^b} \cdot \U^T_{\P^b} \in \Re^{\Nstate \times \Nstate} \,,
\end{eqnarray}
where $\Eigen_{\P^b} \in \Re^{\Nstate \times \Nstate}$ is a diagonal matrix whose diagonal components $\sigma_i$, for $1 \le i \le \Nstate$, are the eigenvalues of $\P^b$ and $\U_{\P^b} \in \Re^{\Nstate \times \Nstate}$ is a set of orthogonal basis vectors spanning the ensemble space (since $\P^b$ is rank deficient). By definitiont ${\rm tr} \lp \P^b \rp = \sum_{i=1}^{\Nstate} \sigma_i$ and ${\rm tr} \lp \lb \P^b \rb^2 \rp = \sum_{i=1}^{\Nstate} \sigma_i^2$. Since there are only $\Nens-1$ eigenvalues different from zero we obtain:
\begin{eqnarray*}
\displaystyle
{\rm tr} \lp \P^b \rp = \sum_{i=1}^{\Nens-1} \sigma_i \,, \qquad
{\rm tr} \lp \lb \P^b \rb^2 \rp = \sum_{i=1}^{\Nens-1} \sigma_i^2 \,,
\end{eqnarray*}
and the computations in the set of equations \eqref{eq:EnkF-optimal-values} can be efficiently performed whenever the first $\Nens-1$ eigenvalues of $\P^b$ can be easily obtained. Consider the singular value decomposition (SVD) for the set of basis vectors \eqref{eq:matrix-of-member-deviations}
\begin{eqnarray}
\label{eq:SVD-over-S}
\S = \U_{\S} \cdot \Sig_{\S} \cdot \V_{\S}^T \in \Re^{\Nstate \times \Nens} \,,
\end{eqnarray}
where $\Sig_{\S} \in \Re^{\Nstate \times \Nens}$ is a diagonal matrix holding the singular values $\widehat{\sigma_i}$ of $\S$, for $1\le i \le \Nens$. Likewise, $\U_{\S} \in \Re^{\Nstate \times \Nstate}$ and $\V_{\S} \in \Re^{\Nens \times \Nens}$ are the left and right singular vectors, respectively. Since $\P^b = \S \cdot \S^T $ we have $\Eigen_{\P^b} = \Sig_{\S} \cdot \Sig_{\S}^T$ and
\begin{eqnarray*}
\displaystyle
{\rm tr} \lp \P^b \rp &=& \sum_{i=1}^{\Nens-1} \sigma_i  = \sum_{i=1}^{\Nens-1} \widehat{\sigma_i}^2, \\
{\rm tr} \lp \lb \P^b \rb^2 \rp &=& \sum_{i=1}^{\Nens-1} \sigma_i^2 =\sum_{i=1}^{\Nens-1} \widehat{\sigma_i}^4. 
\end{eqnarray*}
The computational effort of the SVD decomposition \eqref{eq:SVD-over-S} is $\BO{\Nens^2 \cdot \Nstate}$.
The traces in \eqref{eq:EnkF-optimal-values} can be computing without calculating the sample covariance matrix $\P^b$ by making use of the inexpensive SVD decomposition of $\S$. None of the singular vector of $\S$ are required, but only the singular values $\widehat{\sigma_i}$, for $1 \le i \le \Nens-1$. The parameter values in \eqref{eq:EnkF-optimal-values} are computed as follows:
\begin{subequations}
\label{eq:EnkF-optimal-values-efficient}
\begin{eqnarray}
\label{eq:EnKF-mu-value-efficient}
\displaystyle 
\mu_{\BEST} &=& \frac{\sum_{i=1}^{\Nens-1} \widehat{\sigma_i}^2}{\Nstate} \,, \\
\label{eq:EnKF-gamma-value-efficient}
\gamma_{\B} &=& \min \lp \frac{\frac{\Nens-2}{ \Nstate} \cdot \sum_{i=1}^{\Nens-1} \widehat{\sigma_i}^4 + \lb  \sum_{i=1}^{\Nens-1} \widehat{\sigma_i}^2 \rb^2}{\lp \Nens+2 \rp \cdot \lb \sum_{i=1}^{\Nens-1} \widehat{\sigma_i}^4  - \frac{\lb \sum_{i=1}^{\Nens-1} \widehat{\sigma_i}^2  \rb^2}{\Nstate} \rb} , \, 1 \rp \,.
\end{eqnarray}
\end{subequations}
With $\varphi = \mu_{\BEST} \cdot \gamma_{\BEST}$ and $\delta = 1 - \gamma_{\BEST}$ the estimated covariance matrix \eqref{eq:EnkF-optimal-values} is
\begin{eqnarray}
\label{eq:RBLW-estimator-EnKF}
\displaystyle
\BEST = \varphi \cdot \I_{\Nstate \times \Nstate} + \delta  \cdot \S \cdot \S^T \in \Re^{\Nstate \times \Nstate}\,.
\end{eqnarray}
%
%
\subsection{EnKF implementations based on the RBLW estimator}
\label{subsec:enkf-rblw}

By replacing the estimated error covariance matrix \eqref{eq:RBLW-estimator-EnKF} in \eqref{eq:analysis-ensemble-members}, the EnKF analysis in matrix form becomes
\begin{eqnarray}
\label{eq:EnKF-cov-est}
\displaystyle
\X^a = \X^b + \BEST \cdot \H^T \cdot \lp \R + \H \cdot \BEST \cdot \H^T \rp^{-1} \cdot \D \in \Re^{\Nstate \times \Nens} \,,
\end{eqnarray}
where the matrix of innovations $\D \in \Re^{\Nobs \times \Nens}$ is
\begin{eqnarray}
\label{eq:EnKF-matrix-of-innovations}
\displaystyle
\D = \lb \y_1^s - \H \cdot \x_1^b ,\, \y_2^s - \H \cdot \x_2^b, \, \ldots,\, \y_{\Nens}^s - \H \cdot \x_{\Nens}^b \rb
\end{eqnarray}
and the data $\y^s_i$ for $1 \le i \le \Nens$ is drawn from the distribution \eqref{eq:observation-sample}. We have
\begin{eqnarray}
\nonumber
\displaystyle
\X^a &=& \X^b + \lp \varphi \cdot \I_{\Nstate \times \Nstate} + \delta  \cdot \S \cdot \S^T  \rp \cdot \H^T  \\ \nonumber
& \cdot & \lp \R + \H \cdot \lp \varphi \cdot \I_{\Nstate \times \Nstate} + \delta  \cdot \S \cdot \S^T \rp \cdot \H^T \rp^{-1} \cdot \D \\
\label{eq:EnKF-cov-summ}
\displaystyle \X^a &=& \X^b + \E \cdot \PI \cdot \Z_{\BEST} + \varphi \cdot \H^T \cdot \Z_{\BEST}, 
\end{eqnarray}
where $\E = \sqrt{\delta} \cdot \S \in \Re^{\Nstate \times \Nens}$, $\PI = \H \cdot \E \in \Re^{\Nobs \times \Nens}$, and $\Z_{\BEST} \in \Re^{\Nobs \times \Nens}$ is given by the solution of the linear system
\begin{eqnarray}
\label{eq:solution-linear-system}
\displaystyle \lp \GA+ \PI \cdot \PI^T \rp \cdot \Z_{\B} &=&  \D \,, \\
\nonumber
\GA &=& \R+\varphi \cdot \H \cdot \H^T \in \Re^{\Nobs \times \Nobs}.
\end{eqnarray}
When $\H$ possesses a simple structure (e.g., indexes to observed components from vector states) the matrix $\GA$ also has a simple structure (since in practice $\R$  is block diagonal). By letting $\U^{(0)} = \GA^{-1} \cdot \PI \in \Re^{\Nobs \times \Nens}$ and $\Z_{\B}^{(0)} = \GA^{-1} \cdot \D$, the linear system \eqref{eq:solution-linear-system} can be efficiently solved via the ISMF with no more than $\BO{\Nens^2 \cdot \Nobs}$ long computations. When $\GA$ has no special structure its inverse can be calculated off-line. 

In order to obtain a better representation of the background error statistics (uncertainty) about the background state \eqref{eq:ensemble-mean}  additional samples can be taken from the distribution
\begin{eqnarray}
\label{eq:artificial-member}
\EX{\x}^{b}_{i} \sim \Nor \lp \xm^b ,\, \BEST \rp \in \Re^{\Nstate \times 1} \,,\quad 1 \le i \le K.
\end{eqnarray}
This yields to a new ensemble formed of two kinds of members, real and synthetic. The \textit{real members} $\{\x^b_i\}_{i=1}^{\Nens}$ are obtained by model propagation of the previous analysis ensemble. The \textit{synthetic members} $\{\EX{\x}^b_i\}_{i=1}^{K}$ are artificially built by taking samples from the distribution \eqref{eq:artificial-member} and do not require additional model runs. 

\begin{figure}[H]
     \centering
     \subfloat[$K = 0$]{\includegraphics[width=0.5\textwidth]{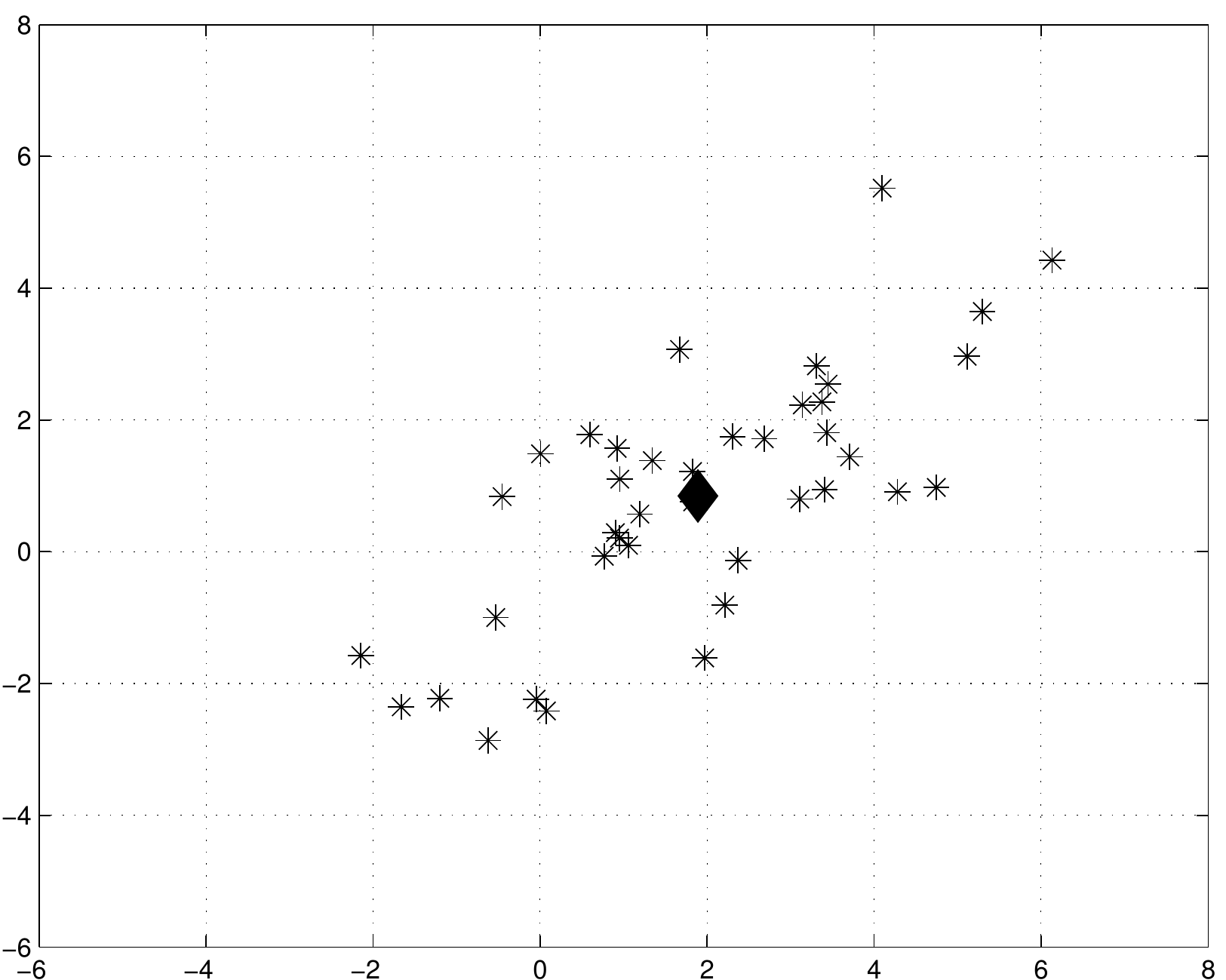}\label{K0}}
%
     \subfloat[$K = 120$]{\includegraphics[width=0.5\textwidth]{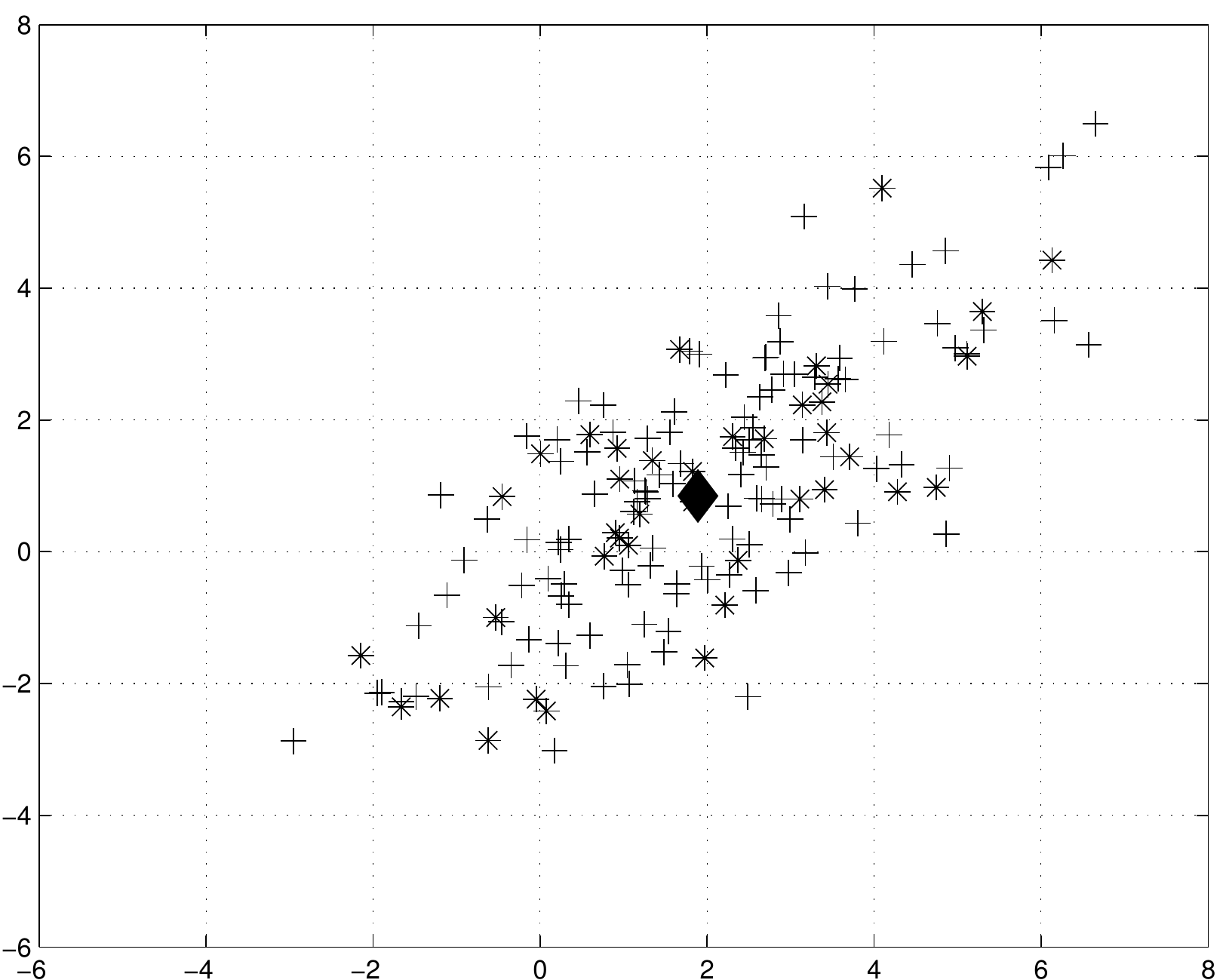}\label{K3}}
         \caption{Error distribution for different values of $K$ using two-dimensional projections of ensemble members from the Lorenz 96 model. The number of real members is $\Nens = 40$. In the plots $*$ represents real ensemble members, $+$ artificial members, and $\diamond$ is the ensemble mean.}
     \label{fig:effect-of-different-K-values}
\end{figure}

The artificial increase in the size of the ensemble is therefore a relatively inexpensive modality to bring in additional degrees of freedom in the solution of the optimization problem \eqref{eq:solution-model-space}. Figure \ref{fig:effect-of-different-K-values} exemplifies the effect of additional members using two-dimensional projections of ensemble member states from the Lorenz 96 model. Figure \ref{K0} shows the spread of the real ensemble members (for the background uncertainty around $\xm^b$). Figures \ref{K3} shows the distribution when artificial members are added to the background ensemble, resulting in a better representation of the background error and therefore a decrease in the sampling error. 

The analysis state is now computed in the subspace spanned by both real and artificial members:
\begin{eqnarray}
\label{eq:extended-basis}
\xm^a \in {\bf span} \lle \x^b_1, \, \x^b_2, \, \ldots, \, \x^b_{\Nens} ,\, \EX{\x}^{b}_1, \,\EX{\x}^{b}_2, \, \ldots, \, \EX{\x}^{b}_{K}  \rle.
\end{eqnarray}
The \textit{extended ensemble} reads:
\begin{eqnarray}
\label{eq:ensemble-extended-background}
\EX{\X}^{b} = \lb \x^b_1, \, \x^b_2, \, \ldots, \, \x^b_{\Nens} ,\, \EX{\x}^{b}_1, \,\EX{\x}^{b}_2, \, \ldots, \, \EX{\x}^{b}_{K} \rb \in \Re^{\Nstate \times \NensK} \,.
\end{eqnarray}
where $\NensK = \Nens+K$. Similar to \eqref{eq:EnKF-cov-summ}, the background ensemble \eqref{eq:ensemble-extended-background} is updated using
\begin{eqnarray}
\label{eq:ensemble-update-extended}
\displaystyle
\X^a = \X^b + \EX{\E} \cdot \EX{\PI} \cdot \EX{\Z}_{\BEST}+\varphi \cdot \H^T \cdot \EX{\Z}_{\B} \in \Re^{\Nstate \times \Nens} \,,
\end{eqnarray}
where
\begin{subequations}
\label{eq:extended-matrices-analysis}
\begin{eqnarray}
\label{eq:extend-E}
\EX{\E} &=& \sqrt{\delta} \cdot \EX{\S}  \in \Re^{\Nstate \times \NensK } \,, \\ 
\label{eq:extend-PI}
\EX{\PI} &=& \H \cdot \EX{\E} \in \Re^{\Nobs \times \NensK} \,, \\ 
\label{eq:extend-Z}
\EX{\Z}_{\BEST} &=& \lp \GA+\EX{\PI} \cdot {\EX{\PI}}^T \rp^{-1} \cdot \D  \in \Re^{\Nobs \times \Nens} \,, \\
\label{eq:extend-S}
\EX{\S} &=& \EX{\X} - \xm^b \otimes \ones_{\NensK}^T \in \Re^{\Nstate \times \NensK} \,.
\end{eqnarray}
\end{subequations}
The solution of \eqref{eq:extend-Z} can be efficiently obtained via the ISMF. Even if the artificial members \eqref{eq:artificial-member} are used in the covariance approximation for the analysis step, according to \eqref{eq:ensemble-update-extended} only the states of real members are adjusted. After the assimilation step the artificial members are discarded and only the real members form the analysis ensemble. This strategy does not increase the number of ensemble members to be propagated by model runs (and therefore, the computational effort). When more computational resources become available, or when some real members are lost due to hardware failures, selected artificial members can be updated and propagated as well. In this case they become real members during the next assimilation step. Moreover, some real members can be replaced by artificial members such as to refresh the ensemble directions, e.g., in order to prevent filter divergence.

To summarize the above discussion, the implementation of the ensemble Kalman filter based on the RBWL estimator with weighted covariance matrix in the observation space (EnKF-FS) consists of the following steps:
\begin{enumerate}
\item Estimate the background error covariance matrix \eqref{eq:RBLW-estimator-EnKF} based on the samples $\lle \x^{b[i]} \rle_{i=1}^{\Nens}$, for $1 \le i \le \Nens$.
\item Compute the innovation matrix $\D$ according to \eqref{eq:EnKF-matrix-of-innovations}.
\item Draw $K$ artificial members from the distribution \eqref{eq:artificial-member}.
\item Compute the set of matrices \eqref{eq:extended-matrices-analysis}.
\item Perform the assimilation \eqref{eq:ensemble-update-extended}.
\item Propagate the ensemble members $$\x^{b}_{\rm next} =  \M_{t_{\rm current} \rightarrow t_{\rm next}} \lp \x^{a}_{\rm current}\rp \,,$$%
until the next assimilation step, for $1 \le i \le \Nens$. 
\end{enumerate}

Another efficient implementation of the filter based on the RBLW estimator can be obtained via the 3D-Var cost function in the ensemble space \eqref{eq:3Dvar-ensemble-space}. This involves projecting the weighted covariance matrix \eqref{eq:weighted-covariance-matrix} onto the ensemble space  \eqref{eq:extended-basis} rather than onto the observation space. 

After taking the samples \eqref{eq:artificial-member} and building the ensemble \eqref{eq:ensemble-extended-background} a new set of basis vectors can be built as follows 
\begin{eqnarray}
\label{eq:basis-vectors-EnKF}
\displaystyle
\EX{\U} = \lb \x_1^b - \xm^b,\, \ldots,\, \x_{\Nens}^b- \xm^b,\EX{\x}_1^{b} - \xm^b,\, \ldots,\, \EX{\x}^b_{K}- \xm^b\rb \in \Re^{\Nstate \times \NensK} \,,
\end{eqnarray}
where $\xm^b$ is given in \eqref{eq:ensemble-mean}. By replacing \eqref{eq:basis-vectors-EnKF} and the RBLW estimator of the background error covariance matrix \eqref{eq:RBLW-estimator-EnKF} in the 3D-Var cost function \eqref{eq:3Dvar-ensemble-space} we obtain
\begin{eqnarray}
\label{eq:extended-3d-var}
\displaystyle 
\J_{\rm ens} \lp \WS\rp =  \frac{1}{2} \cdot \left\Vert \EX{\U} \cdot \WS \right\Vert^2_{\BEST^{-1}} + \frac{1}{2} \cdot \ln \D-\EX{\Q} \cdot \WS \rn_{\R^{-1}}^2
\end{eqnarray}
where $\WS \in \Re^{\NensK \times \Nens}$ is the matrix of weights whose i-th column represent the coordinates of the i-th ensemble member in the space \eqref{eq:extended-basis}, for $1 \le i \le \NensK$, and $\EX{\Q} = \H \cdot \EX{\U} \in \Re^{\Nobs \times \NensK}$. The resulting 3D-Var optimization problem is
\begin{eqnarray}
\label{eq:optimization-problem-EnKF-RBLW-reduce-space}
\displaystyle
\WS^{*} = \underset{\WS}{\arg \min} \,\, \J_{\rm ens} \lp \WS \rp \,,
\end{eqnarray}
and has the solution
\begin{eqnarray}
\label{eq:optimal-weights-extended-ensemble}
\displaystyle
\WS^{*} = \lb \EX{\U}^T \cdot \Z_{\BEST \U} + \EX{\Q}^T \cdot \R^{-1} \cdot \EX{\Q} \rb^{-1} \cdot \EX{\Q}^T \cdot \R^{-1} \cdot \D \in \Re^{\NensK \times \Nens} \,,
\end{eqnarray}
where $\Z_{\BEST \U} = \BEST^{-1} \cdot \U \in \Re^{\Nstate \times \NensK}$. The resulting analysis ensemble is
\begin{eqnarray}
\label{eq:analysis-reduce-space}
\displaystyle
\X^a = \X^b + \EX{\U} \cdot \WS^{*} \in \Re^{\Nstate \times \Nens} \,.
\end{eqnarray}

To summarize, the implementation of the ensemble Kalman filter based on covariance estimation with weighted covariance matrix in the ensemble space (EnKF-FS) consists of the following steps:
\begin{enumerate}
\item Estimate the background error covariance matrix \eqref{eq:RBLW-estimator-EnKF} based on the samples $\lle \x^{b[i]} \rle_{i=1}^{\Nens}$, for $1 \le i \le \Nens$.
\item Compute the innovation matrix $\D$ according to \eqref{eq:EnKF-matrix-of-innovations}.
\item Draw $K$ artificial members according to \eqref{eq:artificial-member}.
\item Compute the matrix of optimal weights \eqref{eq:optimal-weights-extended-ensemble}.
\item Perform the assimilation \eqref{eq:analysis-reduce-space}.
\item Propagate the ensemble members $$\x^{b}_{\rm next} =  \M_{t_{\rm current} \rightarrow t_{\rm next}} \lp \x^{a}_{\rm current}\rp \,,$$%
until the next assimilation step, for $1 \le i \le \Nens$. 
\end{enumerate}
\subsection{Sampling in high-dimensions based on the RBLW estimator}
\label{subsec:sampling-high-dimensions}

Both implementations discussed in Section \ref{subsec:enkf-rblw} use samples from the distribution \eqref{eq:artificial-member}. Such samples can be generated as follows:
\begin{eqnarray}
\label{eq:straight-forward-computation}
\EX{\x}^{b}_i = \xm^b + \BEST^{1/2} \cdot \errbac_i =  \xm^b + \left( \varphi \cdot \I_{\Nstate \times \Nstate} + \delta \cdot \S \cdot \S^T\right)^{1/2} \cdot \errbac_i
\end{eqnarray}
for $1 \le i \le K$, where $\errbac_i \sim \Nor \lp \zero_{\Nstate} ,\, \I_{\Nstate \times \Nstate }\rp$. However, this computation requires the explicit representation in memory of the estimated error covariance matrix $\BEST$, which is prohibitive for high-resolution models. Moreover, the square root matrix $\BEST^{1/2}$ is required making the use of \eqref{eq:straight-forward-computation}  impractical. 

We need an equivalent strategy to obtain the samples \eqref{eq:artificial-member} that requires a reasonable computational effort and does not use a full representation of the covariance matrix $\BEST$. Toward this end consider the random vectors
\begin{eqnarray*}
\displaystyle
\errbac^1_i &\sim & \Nor \lp \zero_{\Nstate},\, \I_{\Nstate \times \Nstate }\rp  \in \Re^{\Nstate \times 1} \,, \\
\errbac^2_i &\sim & \Nor \lp \zero_{\Nens},\, \I_{\Nens \times \Nens }\rp  \in \Re^{\Nens \times 1} \,,
\end{eqnarray*}
and let 
\begin{eqnarray*}
\text{Cov} \lp \errbac^1_i ,\, \errbac^2_i \rp &=& \errbac^1_i \otimes {\errbac^2_i}^T = \zero_{\Nstate \times \Nens} \,, \\
\text{Cov} \lp \errbac_2 ,\, \errbac_1 \rp &=& \errbac^2_i \otimes {\errbac^1_i}^T = \zero_{\Nens \times \Nstate} \,. 
\end{eqnarray*}
We make the following substitution in  \eqref{eq:straight-forward-computation}
\[
{\BEST}^{1/2} \cdot \errbac_i \sim \sqrt{\varphi} \cdot \errbac^1_i + \sqrt{\delta} \cdot \S \cdot \errbac^2_i.
\] 
This does not change the statistics since
\begin{eqnarray*}
&& \expect \lb \lp \sqrt{\varphi} \cdot \I_{\Nstate \times \Nstate} \cdot \errbac^1_i + \sqrt{\delta} \cdot \S \cdot \errbac^2_i \rp \cdot \lp \sqrt{\varphi} \cdot \I_{\Nstate \times \Nstate} \cdot \errbac^1_i + \sqrt{\delta} \cdot \S \cdot \errbac^2_i \rp^T \rb  \\
&&= \varphi \cdot \underbrace{ \errbac^1_i \otimes {\errbac^1_i}^T}_{\text{Cov}\lp \errbac^1_i,\errbac^1_i\rp = \I_{\Nstate \times \Nstate}} + \sqrt{\varphi \cdot \delta} \cdot \underbrace{ \errbac^1_i \otimes {\errbac^2_i}^T}_{\text{Cov}\lp \errbac^1_i,\errbac^2_i\rp = \zero_{\Nstate \times \Nens}} + \sqrt{\varphi \cdot \delta} \cdot \underbrace{ \errbac^2_i \otimes {\errbac^1_i}^T}_{\text{Cov}\lp \errbac^2_i,\errbac^1_i\rp = \zero_{\Nens \times \Nstate}} \\
&&\quad + \delta \cdot \S \cdot \underbrace{ \errbac^2_i \otimes {\errbac^2_i}^T}_{\text{Cov}\lp \errbac^2_i,\errbac^2_i\rp = \I_{\Nens \times \Nens}} \cdot \S^T = \varphi \cdot \I_{\Nstate \times \Nstate} + \delta \cdot \S \cdot \S^T \\
&&= \BEST \,.
\end{eqnarray*}
The artificial ensemble members are obtained as follows:
\begin{eqnarray}
\label{eq:efficient-compute-artificial-members}
\displaystyle
\EX{\x}^{b}_{i} = \xm^b +  \sqrt{\varphi} \cdot \I_{\Nstate \times \Nstate} \cdot \errbac^1_i + \sqrt{\delta} \cdot \S \cdot \errbac^2_i,
\quad i=1,\dots,K.
\end{eqnarray}

The components of the random variables $\errbac^1_i$ and $\errbac^2_i$ are drawn independently from the standard normal distribution $\Nor \lp 0,\,1 \rp$. For large model resolutions the components of the random vectors can be prepared independently taking advantage of parallel computations. Moreover, the random vectors $\errbac^1_i$ and $\errbac^2_i$ can be sampled prior the assimilation process in an off-line computation. 

The estimated error covariance matrix is never represented explicitly in memory. Instead, the estimator $\BEST$ is represented via the triplet
\begin{eqnarray*}
\BEST \equiv \lb \varphi,\, \mu,\, \S \rb  \,.
\end{eqnarray*}
which contains two scalars and one matrix of dimension $\Nstate \times \Nens$. In addition, the scalars $\varphi$ and $\mu$ are computed making use only of the matrix $\S$.  This data is sufficient for correct sampling from the distribution \eqref{eq:artificial-member}.

\subsection{Comparison of EnKF-FS and EnKF-RS versions of the filter}
\label{subsec:comparison-of-the-EnKF-FS-RS}

Although both EnKF-FS and EnKF-RS methods are based on the EnKF equations and RBLW estimator, their underlying theoretical properties are slightly different. To facilitate the comparison of the two proposed implementations we bring the EnKF-FS analysis equation \eqref{eq:ensemble-update-extended} to the form \eqref{eq:analysis-state-general}:
\begin{eqnarray}
\label{eq:equivalent-to-analysis-solution}
\displaystyle
\X^a = \X^b +\BEST^{1/2} \cdot \W_{\BEST} \in \Re^{\Nstate \times \Nens} \,,
\end{eqnarray}
where the weights $\W_{\BEST} \in \Re^{\Nstate \times \Nens}$ are given by
\begin{eqnarray*}
\displaystyle
\W_{\BEST} = \BEST^{1/2} \cdot \H^T \lp \R+\H \cdot \BEST \cdot \H^T \rp^{-1} \cdot \D \in \Re^{\Nstate \times \Nens} \,.
\end{eqnarray*}
It is readily apparent from equations \eqref{eq:equivalent-to-analysis-solution} and \eqref{eq:analysis-reduce-space} that EnKF-FS and EnKF-RS implementations differ in the number of degrees of freedom used in the assimilation process. In the EnKF-FS approach the columns of $\BEST^{1/2}$ serve as the basis set for generating an ensemble of background deviations. Since the estimated background error covariance matrix is full-rank the optimal solution \eqref{eq:optimal-analysis} is searched for in the full space.
The matrix identity
\begin{eqnarray*}
\lb \BEST^{-1} + \H^T \cdot \R^{-1} \cdot \H \rb^{-1} \cdot \H^T \cdot \R^{-1} \equiv \BEST \cdot \H^T \cdot \lb \R+ \H \cdot \BEST \cdot \H^T \rb^{-1}
\end{eqnarray*}
together with \eqref{eq:equivalent-to-analysis-solution} reveal that the weighted covariance matrix of the EnKF-FS implementation
\begin{eqnarray}
\label{eq:weighted-covariance-matrix-for-EnKFFS}
\displaystyle 
{\bf W} = \BEST^{-1} + \H^T \cdot \R^{-1} \cdot \H  \in \Re^{\Nstate \times \Nstate} \,,
\end{eqnarray}
is related to the weighted covariance matrix of the EnKF-RS method by the relation
\begin{eqnarray}
\label{eq:weighted-covariance-matrix-for-EnKFRS}
\displaystyle
{\bf W}_{\rm ens} = \EX{\U}^T \cdot {\bf W} \cdot \EX{\U} \in \Re^{\NensK \times \NensK}\,.
\end{eqnarray}
Threfore when the size of the ensemble $\NensK$ is increased (by adding real or artificial members) more information from the matrix ${\bf W}$ is captured by its projection onto the $\NensK$-dimensional space. Note that when $\NensK \rightarrow \Nstate$ and $\EX{\U}$ is orthonormal we have that ${\bf W}_{\rm ens} \to {\bf W}$. Consequently, the number of artificial members will play an important role in the performance of the EnKF-RS implementation.

\section{Experimental Results}
\label{sec:experimental-results}

This section tests the new EnKF implementations on a data assimilation problem using the quasi-geostrophic model presented in \cite{TELA:TELA299}. A comparison is done in two steps: first, the proposed implementations are compared against the well-known EnKF implementations presented in section \ref{sec:preliminaries}, and next the quality of the results for the EnKF-FS and EnKF-RS based on different values of $\Nens$ and $K$ are assessed. 

The oceans form a complex flow system influenced by the rotation of the Earth, the density stratification due to temperature and salinity, as well as other factors. The quasi-geostrophic (QG) model is a simple approximation of the
real behavior of the ocean. It is defined by the following partial differential equation:
\begin{eqnarray}
\label{eq:QG-model}
\displaystyle { \omega}_t + r \cdot J \lp { \omega},\, { \psi} \rp + \beta \cdot { \psi}_{x} - v \cdot \nabla^4 { \psi} = - \mu \cdot \nabla^2 { \psi} + \tau \cdot \sin \lp \frac{2 \cdot \pi \cdot y }{L_y} \rp 
\end{eqnarray}
in $\Omega \in [0,\, L_x] \times [0,\,L_y]$, where $x$ and $y$ represent the horizontal and vertical space components, ${\omega}$ is the vorticity, $\psi$ is the stream function, $J \lp { \psi}, { \omega}\rp $ is the Jacobian of two fields
\begin{eqnarray}
\label{eq:QG-Jacobian}
\displaystyle
J \lp { \psi}, { \omega}\rp  = { \psi}_x \cdot { \omega}_y - { \psi}_y \cdot { \omega}_x \,,
\end{eqnarray}
and $\nabla^2$ is the Laplacian operator. The coefficients $\beta$, $v$, $\mu$ and $\tau$ are associated with the horizontal vorticity, the horizontal friction, the biharmonic horizontal friction, and the horizontal wind stress at the surface of the ocean, respectively. Moreover, the vorticity is related to the stream function by the elliptical equation:
\begin{eqnarray}
\displaystyle 
\nabla^2 { \psi} = {\omega} \,.
\end{eqnarray}
The spatial domain for our experiments is $\Omega = [0,\,1] \times [0,\,1]$. The interior is covered by computational grids of different resolutions, and we denote by $D_1$ and $D_2$ the number of horizontal and vertical grid points, respectively. The different model resolutions are presented in the table \ref{tab:qg-instances}.
\begin{table}[H]
\centering
\begin{tabular}{|c|c|c|c|c|c|c|} \hline
Instance & $D_1$ & $D_2$ & $\Nstate =  D_1 \cdot D_2$ \\ \hline
$\qga$& 31 & 31 & 961 \\ \hline
$\qgb$& 63 & 63 & 3,969 \\ \hline
$\qgc$& 127 & 127 & 16,129 \\ \hline
\end{tabular}
\caption{Quasi-geostrophic instances for the computational tests in terms of the number of horizontal $D_1$ and vertical $D_2$ grid points. $D_1$ and $D_2$ do not consider boundary points.}
\label{tab:qg-instances}
\end{table}

The numerical data assimilation experiments are characterized by the following settings:
\begin{itemize}
\item Initial vorticities have the form 
\begin{eqnarray*}
\displaystyle
{\omega_0} = \sin(4 \, x_i \, y_j) \cdot \cos(2 \, x_i \, y_j)+\sin(2 \, x_i \, y_j)+\cos(4 \, x_i \, y_j)\,,
\end{eqnarray*}
for $1 \le i \le D_1$ and $1 \le j \le D_2$.
\item The initial background error covariance is
\begin{eqnarray*}
\displaystyle
\B = \lb \sigma^{\B} \rb^2 \cdot \I_{\Nstate \times \Nstate} ,\,
\end{eqnarray*}
where the standard deviation  $\sigma^{\B}$ is chosen to be $0.05$ or $0.15$ (times the true vorticity) in different experiments . 
\item The observational errors are uncorrelated with variances $\lb 0.01 \rb^2$.
\item The number of observed components from the vector state is given by
\begin{eqnarray}
\Nobs = \Pobs \cdot \Nstate \,,
\end{eqnarray}
where $\Pobs$ is the percentage of observed components. We consider two values for $\Pobs$, 0.7 and 0.9, corresponding to a sparser and a denser network, respectively.
\end{itemize}
Other aspects of the numerical simulation are described below:
\begin{itemize}
\item The EnKF methods are implemented in C++ while the forward model (QG) \cite{QGPavelSakov} is implemented making use of FORTRAN. 
\item The partial derivatives are discretized by central finite differences.
\item The matrix and vector computations are efficiently carried out using the BLAS library \cite{BLASLIB}.
\item Matrix decompositions as well as eigenvalue computations are performed using the LAPACK library \cite{LAPACKLIB}.
\item The Arakawa method \cite{Jespersen1974383} is utilized in order to compute the Jacobian \eqref{eq:QG-Jacobian}.
\item The time discretization of the model \eqref{eq:QG-model} uses of a fourth order Runge-Kutta method. The time step size $1.27$ (units) which represent one hour in the ocean. The integration is performed for 1000 hours.
\item The NLOPT library \cite{NLOPT} is utilized to numerically solve the optimization problems \eqref{eq:cost-function-EnKF-N} and \eqref{eq:cost-function-EnKF-DU}.
\item The GSL-GNU scientific library \cite{GSL-GNU} is utilized  to generate the synthetic background and data errors.
\end{itemize}

\subsection{Comparison with current EnKF implementations}

We compare the EnKF-RS and EnKF-FS methods against current well-known EnKF implementations. The root mean square errors (RMSE) and CPU times for different methods making use of the $\qga$, $\qgb$ and $\qgc$ model instances are reported inTables \ref{tab:rmse-time-nrefin-5}, \ref{tab:rmse-time-nrefin-6}, and \ref{tab:rmse-time-nrefin-7}, respectively. {For a given number of snapshots $\lle \x^a_i \rle_{i=1}^{\N}$ and a reference trajectory $\lle \x^{\rm true}_i \rle_{i=1}^{\N}$, where $\N$ is the number of assimilation times, the RMSE is defined as follows
\begin{eqnarray*}
\text{RMSE} = \sqrt{\frac{1}{\N} \cdot \sum_{i=1}^\N \ln \x_i^a - \x^{\rm true}_i\rn^2 }.
\end{eqnarray*}

We vary the size of the ensemble $\Nens$, the percentage of observed components $\Pobs$, and the initial background error $\sigma^{\B}$. As expected the traditional ensemble implementations  EnKF, EnSRF, and EnTKF provide accurate analyses for  different model resolutions and number of observed components. Moreover it can be seen in figures \ref{fig:Comparison-Nens-40} and \ref{fig:Comparison-Nens-80} that the RMSE decreases as the simulation progresses. The traditional EnKF implementations provide the lowest elapsed time among the compared methods (i.e., EnKF and EnSRF implementations) which explains why they are attractive for  use in real applications. The RMSE values for the EnSRF and EnTKF implementations are identical since both filters are deterministic and  EnTKF is just an efficient implementation of the EnSRF. The inflation-free methods such as the EnKF-FN and EnK-DU implementations provide slightly better accurate results than other methods. For example, for the largest instance $\qgc$, table \ref{tab:rmse-time-nrefin-7} shows the EnKF-DU to perform better than current implementations. In most of the cases, strong duality holds: $\J^\textsc{fn}_{ens}(\W^{*}) = \DU_{ens}^\textsc{du}(\zeta^{*})$ and the slight differences between the optimal cost function values \eqref{eq:cost-function-EnKF-N} and \eqref{eq:cost-function-EnKF-DU} are consistent with the numerical approximation errors in the solution of the optimization problems \eqref{eq:optimization-problem-EnKF-FN} and \eqref{eq:optimization-problem-EnKF-DU}, respectively. 

Since the analysis state in the EnKF-DU formulation is obtained via the solution of the one-dimensional optimization problem \eqref{eq:optimization-problem-EnKF-DU}, we expect this method to be faster than the EnKF-FN implementation where the analysis requires the solution of the $\Nens$-dimensional optimization problem \eqref{eq:cost-function-EnKF-N}. This fact is also pointed out by Boquet in \cite[Section 2.2]{Bocquet2012}, where the cost of computing the inverse \eqref{eq:weighted-covariance-W-Dual} is assumed to be negligible in the dual formulation. However, this statement seems to be true for small model resolutions (i.e., Bocquet makes use of the Lorenz 96 model with 40 variables) and it holds for the smallest QG instance $\qga$. Nevertheless, for the $\qgb$ case, the difference between the CPU times for the primal and dual implementations is almost negligible and even more, for the largest instance $\qgc$, the EnKF-FN performs better than its dual approach for $\sigma^{\B} = 0.15$. Although the cost function \eqref{eq:cost-function-EnKF-DU} depends only on $\zeta$, every step in the optimization process requires the solution of the linear system \eqref{eq:weighted-covariance-W-Dual} whose computational cost is not negligible in practice. Furthermore, when the initial background error is large the EnKF-DU computes many times the inverse of \eqref{eq:weighted-covariance-W-Dual} and therefore its performance decreases considerably. 

Even if the traditional implementations perform very well in terms of RMSE and elapsed time, the most accurate results are obtained by the proposed new EnKF implementations. The results presented in figures \ref{fig-40-07-001-005-5} and \ref{fig-80-07-001-005-5} show that the EnKF-RS method performs much better than traditional and inflation-free methods for the small instance $\qga$, and are slightly better in the larger instances. The most accurate results among all the compared methods are the ones obtained by the EnKF-FS implementation. The results reported in figures \ref{fig:Comparison-Nens-40} and \ref{fig:Comparison-Nens-80} show that, for all the instances and configurations, the RMSE obtained via the EnKF-FS outperforms the other methods by at least 60\%. The CPU times for the proposed implementations are just slightly larger than those from the compared EnKF implementations. Since most of the computational time is spent in propagating the ensemble members, a modest increase in analysis time retains the potential of the new methods to perform well in practical applications.
\begin{table}[H]
\centering
\begin{tabular}{|c|c|l|c|c|c|c|} \cline{4-7}
\multicolumn{1}{c}{ } &  \multicolumn{1}{c}{ } &  \multicolumn{1}{c}{ } &\multicolumn{2}{|c|}{RMSE} & \multicolumn{2}{|c|}{CPU Time}   \\ \cline{4-7}
 \multicolumn{1}{c}{ }&  \multicolumn{1}{c}{ } &  \multicolumn{1}{c}{ } &\multicolumn{2}{|c|}{$\sigma^{\B}$} & \multicolumn{2}{|c|}{$\sigma^{\B}$}   \\ \hline
$\Nens$ & $\Pobs$ & Method & 0.05 & 0.15 & 0.05 & 0.15  \\ \hline
\multirow{14}{*}{40} & \multirow{7}{*}{0.7} & EnKF & 1.287 & 3.858& 0.013 & 0.009 \\
 &  & EnSRF & 1.289 & 3.864& 0.013 & 0.013 \\
 &  & EnTKF & 1.289 & 3.864& 0.026 & 0.028 \\
 &  & EnKF-FN & 1.286 & 3.855& 0.125 & 0.116 \\
 &  & EnKF-DU & 1.285 & 3.854& 0.045 & 0.036 \\
 &  & EnKF-FS & 0.659 & 1.974& 0.024 & 0.036 \\
 &  & EnKF-RS & 1.16 & 3.49& 0.421 & 0.252 \\ \cline{2-7}
 &  \multirow{7}{*}{0.9} & EnKF & 1.279 & 3.836& 0.019 & 0.021 \\
 &  & EnSRF & 1.282 & 3.841& 0.01 & 0.017 \\
 &  & EnTKF & 1.282 & 3.841& 0.028 & 0.017 \\
 &  & EnKF-FN & 1.279 & 3.831& 0.115 & 0.115 \\
 &  & EnKF-DU & 1.276 & 3.828& 0.072 & 0.087 \\
 &  & EnKF-FS & 0.371 & 1.116& 0.048 & 0.028 \\
 &  & EnKF-RS & 1.07 & 3.209& 0.395 & 0.222 \\ \hline
\multirow{14}{*}{80} &  \multirow{7}{*}{0.7} & EnKF & 1.268 & 3.803& 0.065 & 0.061 \\
 &  & EnSRF & 1.275 & 3.82& 0.034 & 0.052 \\
 &  & EnTKF & 1.275 & 3.82& 0.143 & 0.076 \\
 &  & EnKF-FN & 1.264 & 3.79& 0.608 & 0.681 \\
 &  & EnKF-DU & 1.263 & 3.79& 0.162 & 0.181 \\
 &  & EnKF-FS & 0.646 & 1.927& 0.119 & 0.105 \\
 &  & EnKF-RS & 1.069 & 3.144& 1.303 & 1.285 \\ \cline{2-7}
 &  \multirow{7}{*}{0.9} & EnKF & 1.252 & 3.756& 0.074 & 0.077 \\
 &  & EnSRF & 1.26 & 3.773& 0.058 & 0.051 \\
 &  & EnTKF & 1.26 & 3.773& 0.168 & 0.12 \\
 &  & EnKF-FN & 1.249 & 3.737& 0.894 & 0.766 \\
 &  & EnKF-DU & 1.246 & 3.737& 0.155 & 0.181 \\
 &  & EnKF-FS & 0.358 & 1.074& 0.192 & 0.14 \\
 &  & EnKF-RS & 0.879 & 2.588& 1.269 & 1.26 \\ \hline
\end{tabular}
\caption{RMSE and CPU-time (TIME) for the EnKF, EnSRF, EnTKF, EnKF-FN, EnKF-DU, EnKF-RS and EnKF-FS implementations applied to the $\qga$ instance ($\Nstate = 961$).}
\label{tab:rmse-time-nrefin-5}
\end{table}
\begin{table}[H]
\centering
\begin{tabular}{|c|c|l|c|c|c|c|} \cline{4-7}
\multicolumn{1}{c}{ } &  \multicolumn{1}{c}{ } &  \multicolumn{1}{c}{ } &\multicolumn{2}{|c|}{RMSE} & \multicolumn{2}{|c|}{CPU Time}   \\ \cline{4-7}
 \multicolumn{1}{c}{ }&  \multicolumn{1}{c}{ } &  \multicolumn{1}{c}{ } &\multicolumn{2}{|c|}{$\sigma^{\B}$} & \multicolumn{2}{|c|}{$\sigma^{\B}$}   \\ \hline
$\Nens$ & $\Pobs$ & Method & 0.05 & 0.15 & 0.05 & 0.15  \\ \hline
\multirow{14}{*}{40} & \multirow{7}{*}{0.7} & EnKF & 1.67 & 5.011& 0.072 & 0.076 \\
 &  & EnSRF & 1.675 & 5.025& 0.059 & 0.047 \\
 &  & EnTKF & 1.675 & 5.025& 0.056 & 0.098 \\
 &  & EnKF-FN & 1.671 & 4.987& 0.185 & 0.261 \\
 &  & EnKF-DU & 1.662 & 4.987& 0.17 & 0.245 \\
 &  & EnKF-FS & 0.893 & 2.668& 0.108 & 0.205 \\
 &  & EnKF-RS & 1.61 & 4.818& 1.869 & 1.801 \\  \cline{2-7}
 & \multirow{7}{*}{0.9} & EnKF & 1.667 & 5.001& 0.095 & 0.061 \\
 &  & EnSRF & 1.671 & 5.015& 0.06 & 0.035 \\
 &  & EnTKF & 1.671 & 5.015& 0.067 & 0.117 \\
 &  & EnKF-FN & 1.667 & 4.975& 0.262 & 0.269 \\
 &  & EnKF-DU & 1.658 & 4.975& 0.273 & 0.244 \\
 &  & EnKF-FS & 0.52 & 1.533& 0.209 & 0.139 \\
 &  & EnKF-RS & 1.589 & 4.761& 1.975 & 2.011 \\ \hline
\multirow{14}{*}{80} & \multirow{7}{*}{0.7} & EnKF & 1.653 & 4.958& 0.189 & 0.189 \\
 &  & EnSRF & 1.661 & 4.982& 0.174 & 0.183 \\
 &  & EnTKF & 1.661 & 4.982& 0.265 & 0.377 \\
 &  & EnKF-FN & 1.644 & 4.915& 0.934 & 0.919 \\
 &  & EnKF-DU & 1.638 & 4.916& 0.567 & 0.787 \\
 &  & EnKF-FS & 0.866 & 2.592& 0.634 & 0.495 \\
 &  & EnKF-RS & 1.546 & 4.612& 7.627 & 7.542 \\ \cline{2-7}
 & \multirow{7}{*}{0.9} & EnKF & 1.65 & 4.949& 0.271 & 0.302 \\ 
 &  & EnSRF & 1.657 & 4.97& 0.227 & 0.181 \\
 &  & EnTKF & 1.657 & 4.97& 0.428 & 0.429 \\
 &  & EnKF-FN & 1.639 & 4.888& 1.099 & 1.01 \\
 &  & EnKF-DU & 1.63 & 4.888& 0.85 & 0.936 \\
 &  & EnKF-FS & 0.494 & 1.486& 0.607 & 0.502 \\
 &  & EnKF-RS & 1.507 & 4.509& 7.845 & 8.012 \\ \hline
\end{tabular}
\caption{RMSE and CPU-time (TIME) for the EnKF, EnSRF, EnTKF, EnKF-FN, EnKF-DU, EnKF-RS and EnKF-FS implementations applied to the $\qgb$ instance ($\Nstate = 3969$).}
\label{tab:rmse-time-nrefin-6}
\end{table}
\begin{table}[H]
\centering
\begin{tabular}{|c|c|l|c|c|c|c|} \cline{4-7}
\multicolumn{1}{c}{ } &  \multicolumn{1}{c}{ } &  \multicolumn{1}{c}{ } &\multicolumn{2}{|c|}{RMSE} & \multicolumn{2}{|c|}{CPU Time}   \\ \cline{4-7}
 \multicolumn{1}{c}{ }&  \multicolumn{1}{c}{ } &  \multicolumn{1}{c}{ } &\multicolumn{2}{|c|}{$\sigma^{\B}$} & \multicolumn{2}{|c|}{$\sigma^{\B}$}   \\ \hline
$\Nens$ & $\Pobs$ & Method & 0.05 & 0.15 & 0.05 & 0.15  \\ \hline
\multirow{14}{*}{40}  & \multirow{7}{*}{0.7}  & EnKF & 1.708 & 5.125& 0.307 & 0.186 \\
 &  & EnSRF & 1.712 & 5.136& 0.174 & 0.153 \\
 &  & EnTKF & 1.712 & 5.136& 0.237 & 0.223 \\
 &  & EnKF-FN & 1.707 & 5.1& 0.793 & 0.589 \\
 &  & EnKF-DU & 1.7 & 5.1& 0.788 & 0.746 \\
 &  & EnKF-FS & 0.963 & 2.847& 0.572 & 0.611 \\
 &  & EnKF-RS & 1.682 & 5.042& 5.607 & 5.832 \\ \cline{2-7}
 & \multirow{7}{*}{0.9}  & EnKF & 1.709 & 5.127& 0.263 & 0.277 \\
 &  & EnSRF & 1.712 & 5.134& 0.166 & 0.192 \\
 &  & EnTKF & 1.712 & 5.134& 0.294 & 0.296 \\
 &  & EnKF-FN & 1.704 & 5.089& 0.888 & 0.849 \\
 &  & EnKF-DU & 1.696 & 5.089& 1.047 & 1.041 \\
 &  & EnKF-FS & 0.582 & 1.661& 0.627 & 0.745 \\
 &  & EnKF-RS & 1.676 & 5.02& 6.368 & 6.292 \\ \hline
\multirow{14}{*}{80}  & \multirow{7}{*}{0.7}  & EnKF & 1.702 & 5.104& 0.721 & 0.693 \\
 &  & EnSRF & 1.709 & 5.125& 0.554 & 0.576 \\
 &  & EnTKF & 1.709 & 5.125& 0.838 & 0.734 \\
 &  & EnKF-FN & 1.69 & 5.063& 3.085 & 2.596 \\
 &  & EnKF-DU & 1.687 & 5.063& 2.826 & 2.931 \\
 &  & EnKF-FS & 0.941 & 2.81& 2.191 & 2.373 \\
 &  & EnKF-RS & 1.654 & 4.94& 26.384 & 26.355 \\ \cline{2-7}
 & \multirow{7}{*}{0.9}  & EnKF & 1.7 & 5.099& 0.834 & 0.991 \\
 &  & EnSRF & 1.708 & 5.123& 0.667 & 0.729 \\
 &  & EnTKF & 1.708 & 5.123& 0.924 & 0.963 \\
 &  & EnKF-FN & 1.683 & 5.034& 4.07 & 3.485 \\
 &  & EnKF-DU & 1.678 & 5.034& 3.565 & 3.501 \\
 &  & EnKF-FS & 0.558 & 1.622& 2.745 & 2.442 \\
 &  & EnKF-RS & 1.634 & 4.908& 27.672 & 26.208 \\ \hline
\end{tabular}
\caption{RMSE and CPU-time for for the EnKF, EnSRF, EnTKF, EnKF-FN, EnKF-DU, EnKF-RS and EnKF-FS implementations applied to the $\qgc$ instance ($\Nstate = 16129$).}
\label{tab:rmse-time-nrefin-7}
\end{table}%
\begin{figure}
     \centering
     \subfloat[$\Nstate = 961$ and $\sigma^B = 0.05$]{\includegraphics[width=0.4\textwidth]{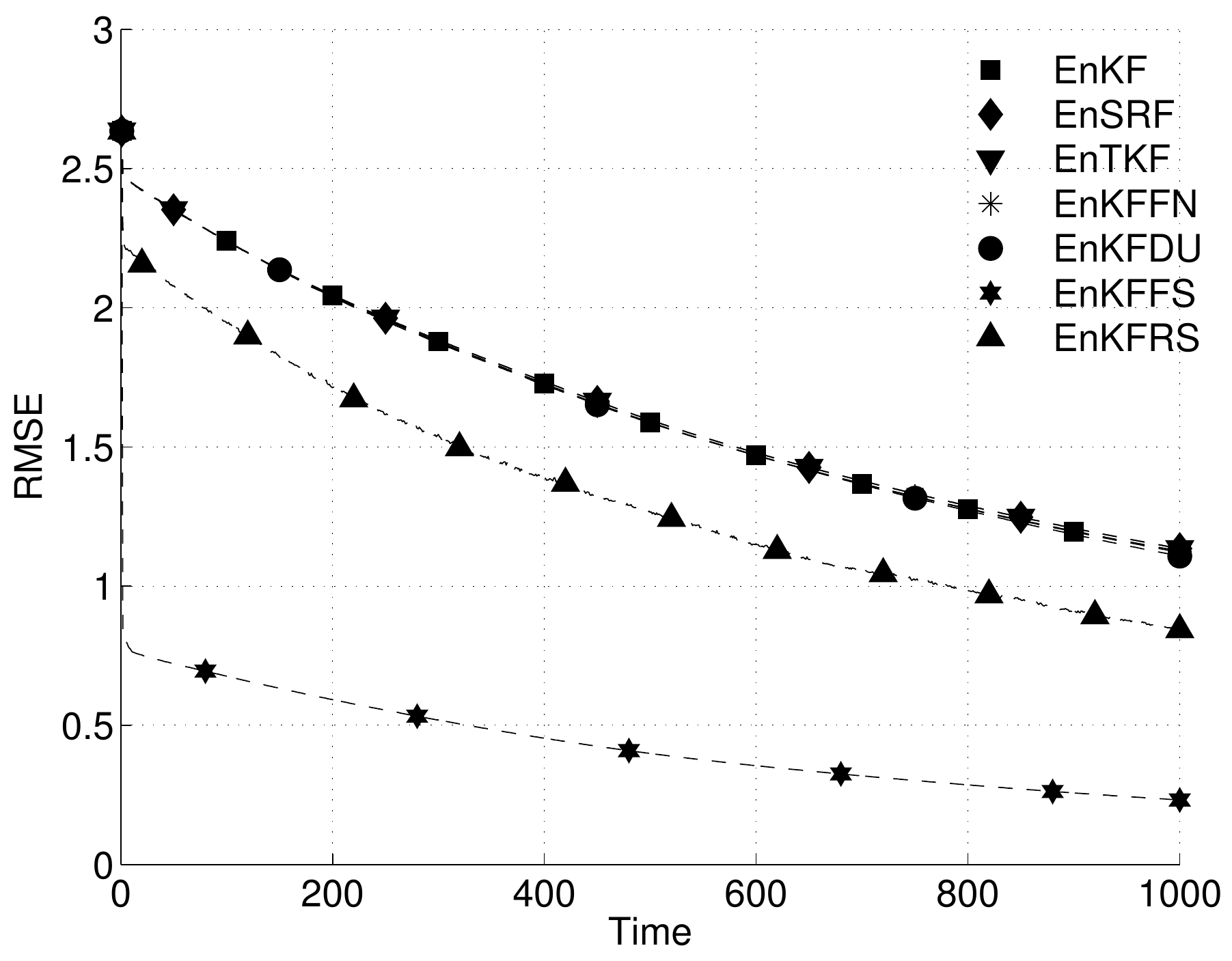}\label{fig-40-07-001-005-5}} %
          \subfloat[$\Nstate = 961$ and $\sigma^B = 0.15$]{\includegraphics[width=0.4\textwidth]{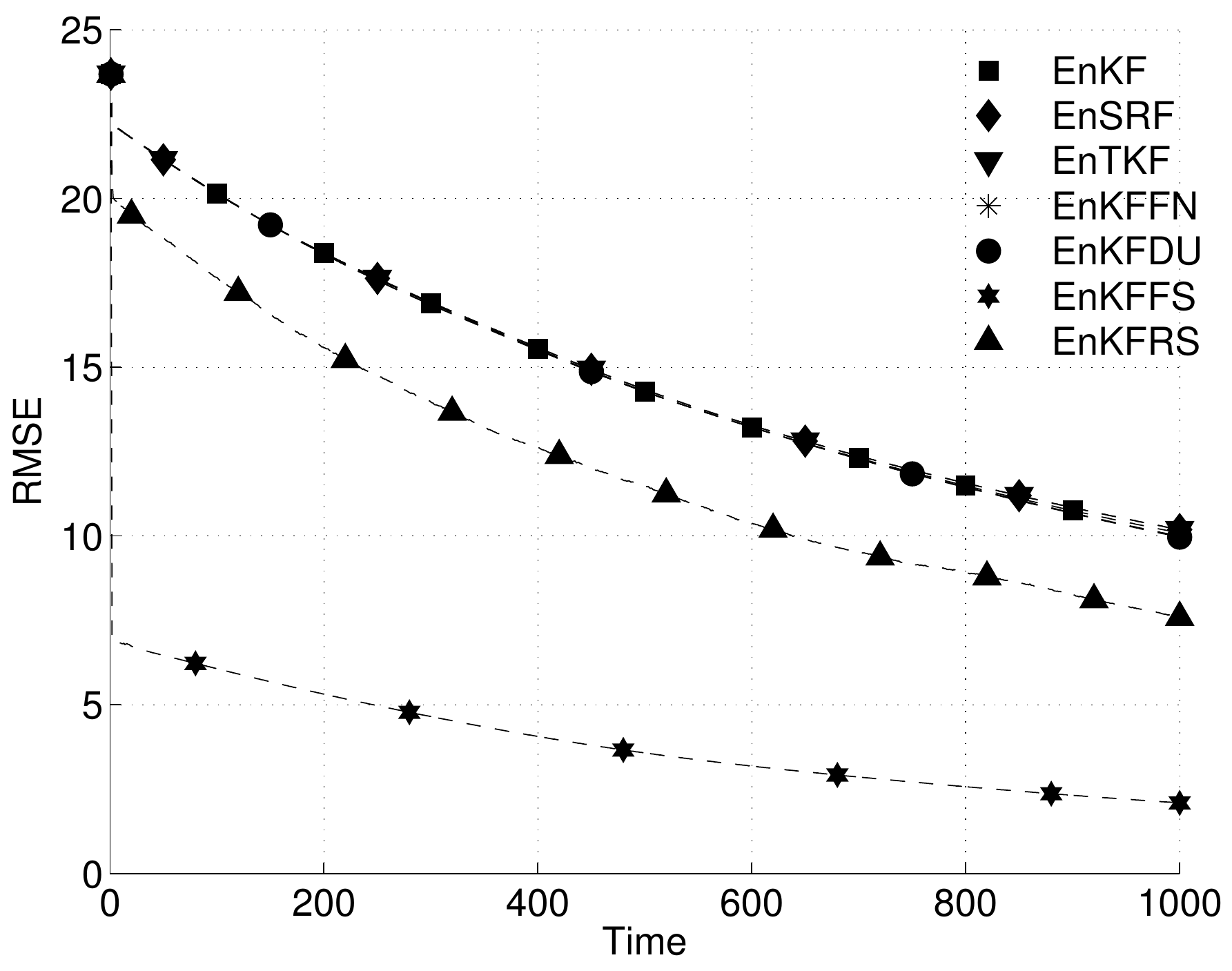}\label{fig-40-07-001-015-5}}
          
     \subfloat[$\Nstate =  3969$ and $\sigma^B = 0.05$]{\includegraphics[width=0.4\textwidth]{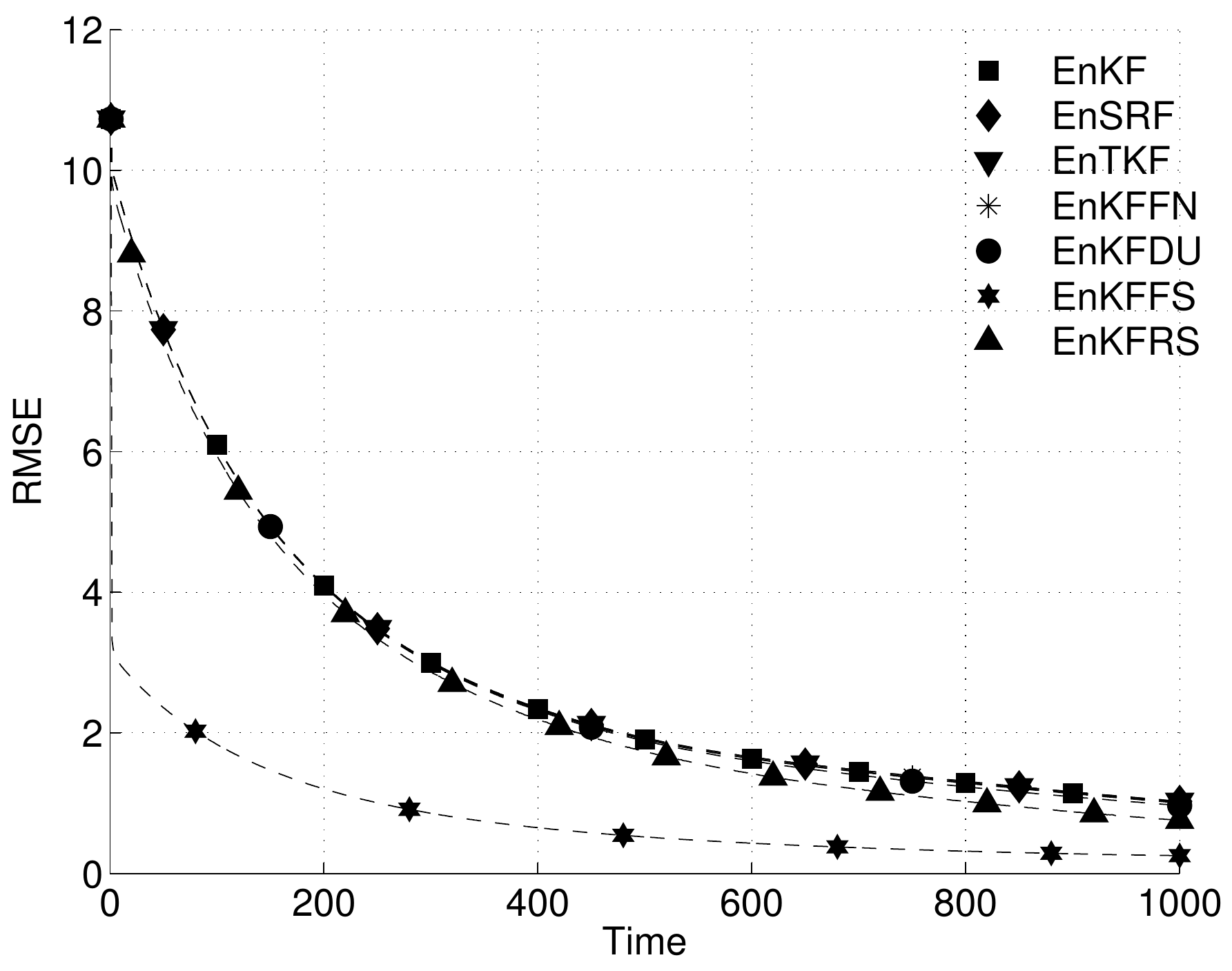}\label{fig-40-07-001-005-6}} %
          \subfloat[$\Nstate =  3969$ and $\sigma^B = 0.15$]{\includegraphics[width=0.4\textwidth]{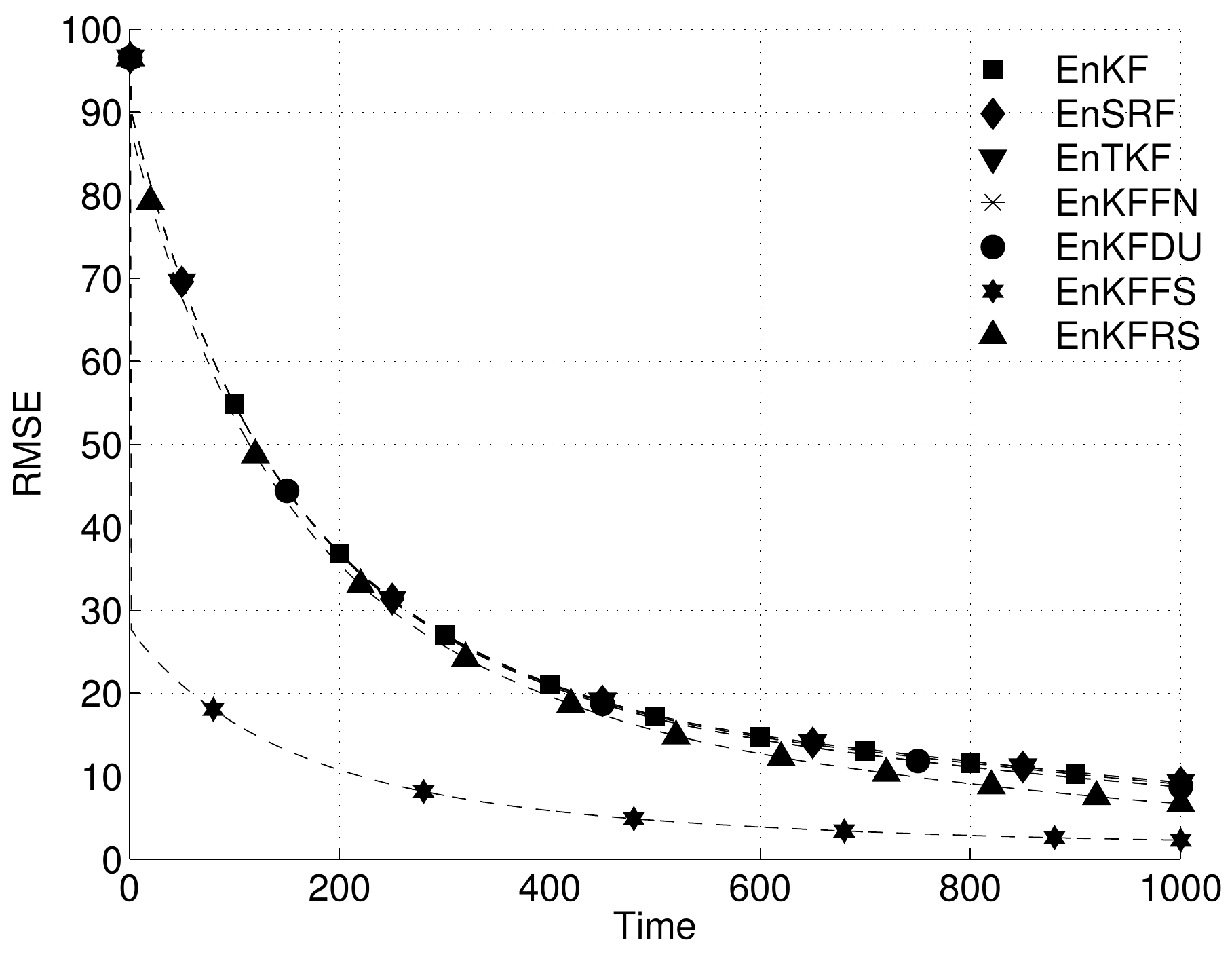}\label{fig-40-07-001-015-6}}
          
     \subfloat[$\Nstate =  16129$ and $\sigma^B = 0.05$]{\includegraphics[width=0.4\textwidth]{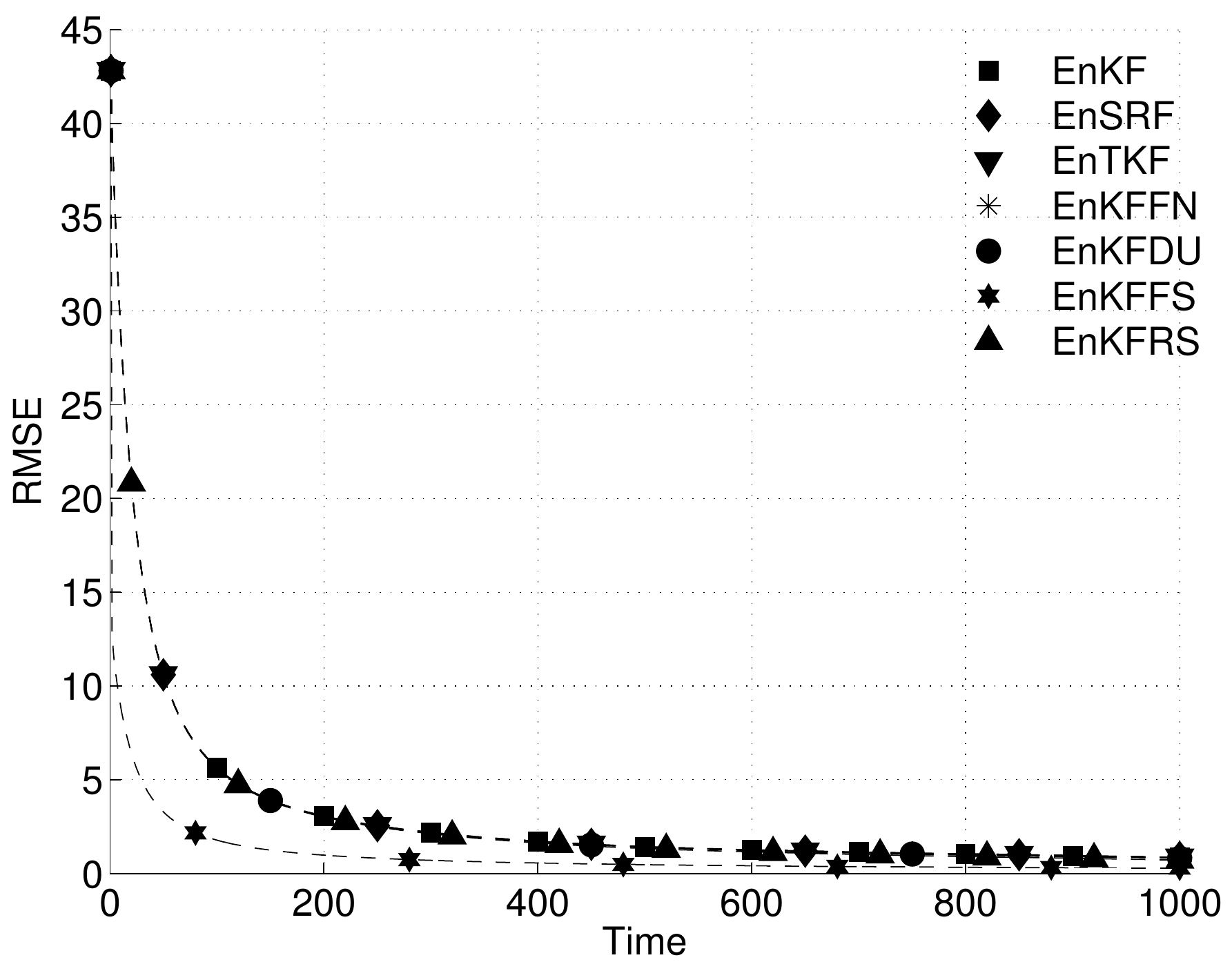}\label{fig-40-07-001-005-7}} %
          \subfloat[$\Nstate =  16129$ and $\sigma^B = 0.15$]{\includegraphics[width=0.4\textwidth]{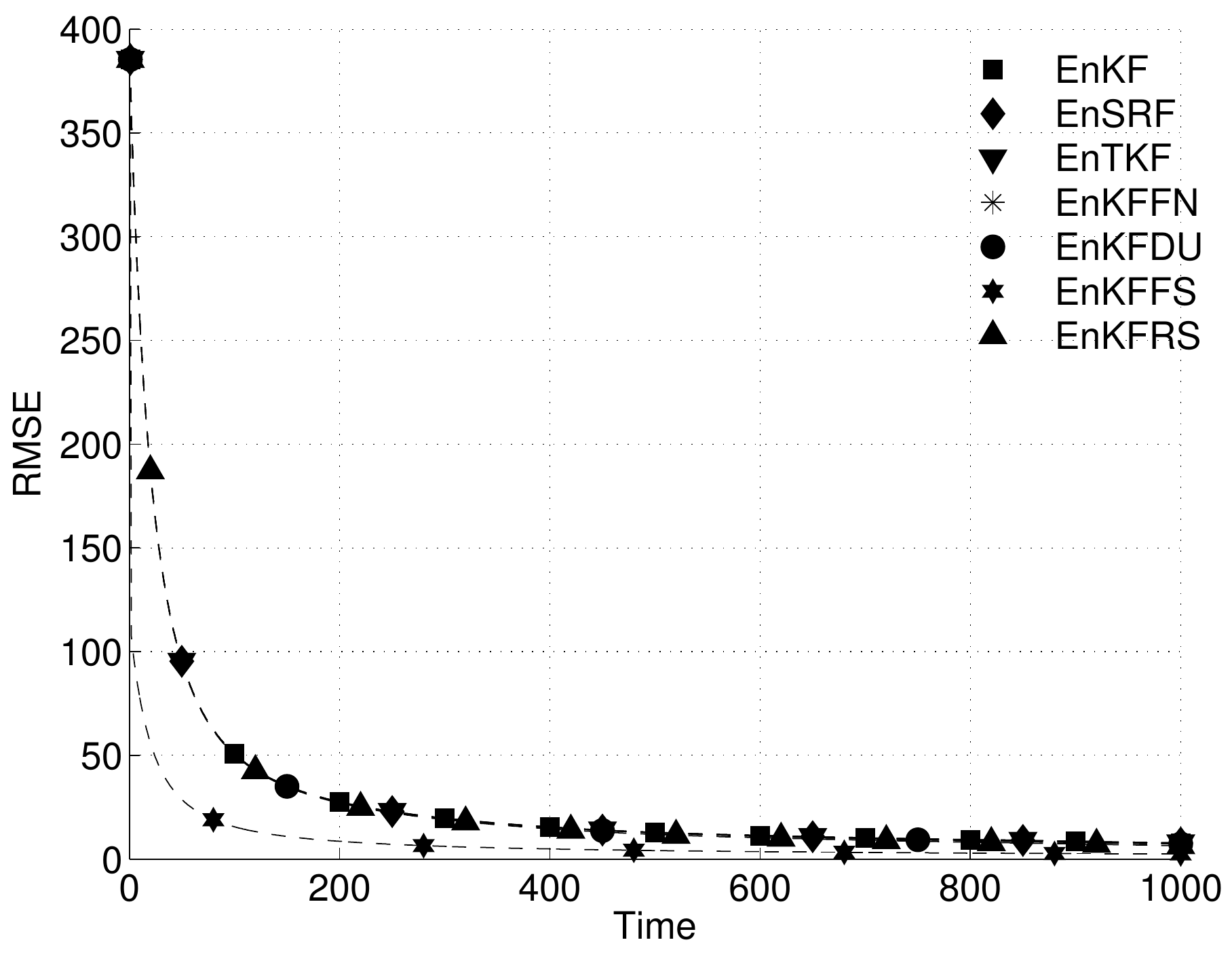}\label{fig-40-07-001-015-7}}
     \caption{Plots of RMSE values for the EnKF, EnSRF, EnTKF, EnKF-FN, EnKF-DU, EnKF-FS and EnKF-RS for $\Nens = 40$ and $\Pobs = 0.7$.}
     \label{fig:Comparison-Nens-40}
\end{figure}
\begin{figure}
     \centering
     \subfloat[$\Nstate = 961$ and $\sigma^B = 0.05$]{\includegraphics[width=0.4\textwidth]{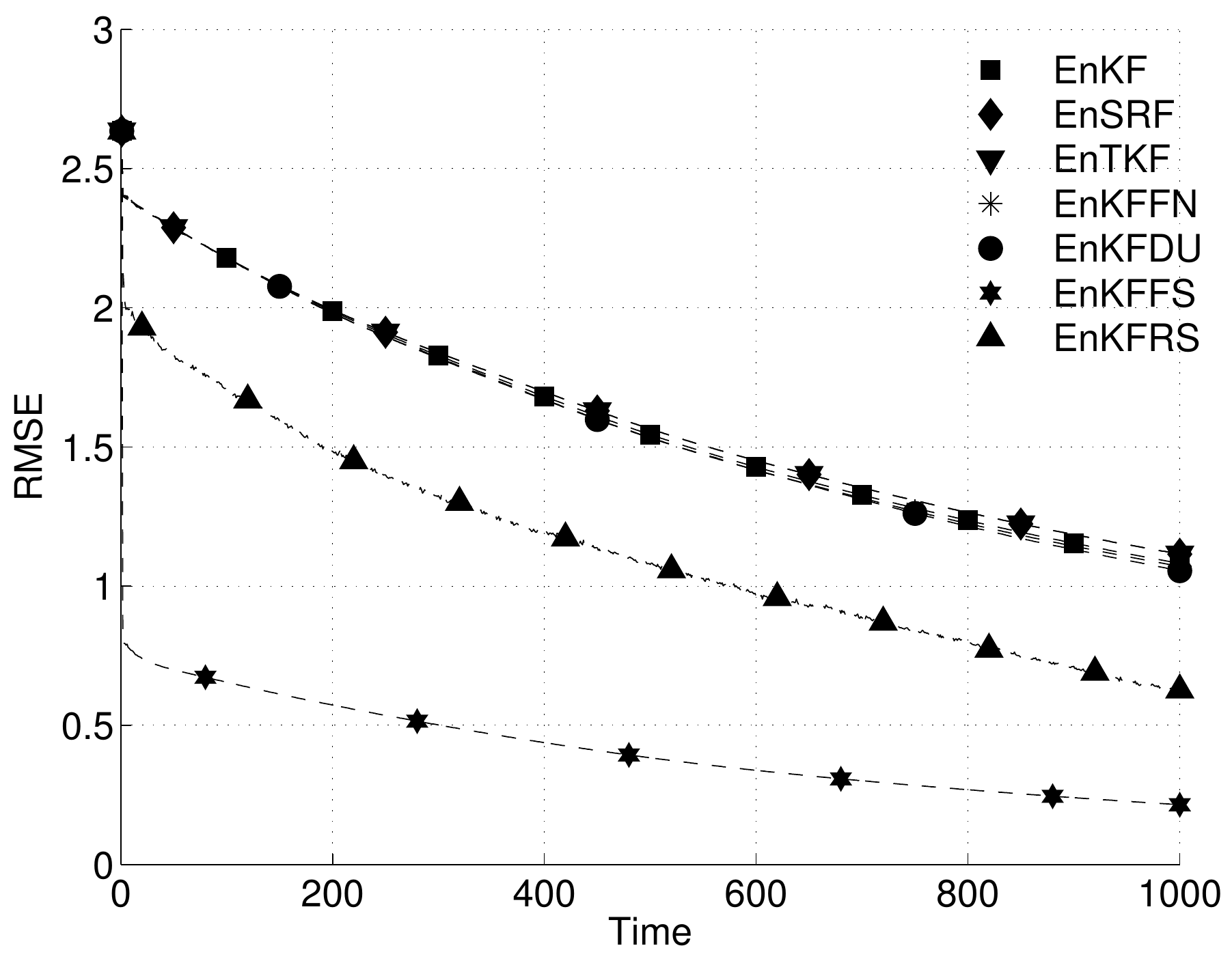}\label{fig-80-07-001-005-5}} %
          \subfloat[$\Nstate = 961$ and $\sigma^B = 0.15$]{\includegraphics[width=0.4\textwidth]{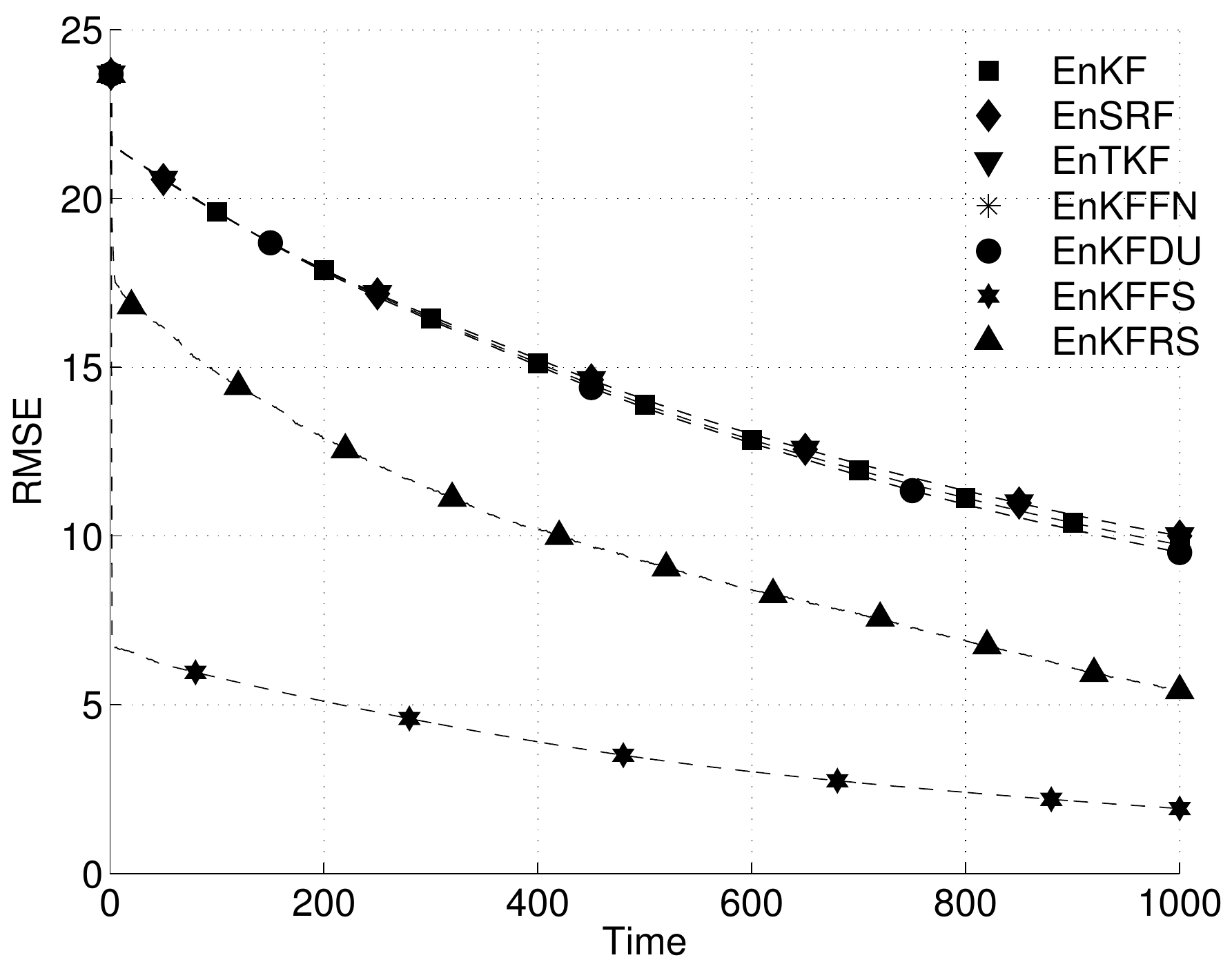}\label{fig-80-07-001-015-5}}
          
     \subfloat[$\Nstate =  3969$ and $\sigma^B = 0.05$]{\includegraphics[width=0.4\textwidth]{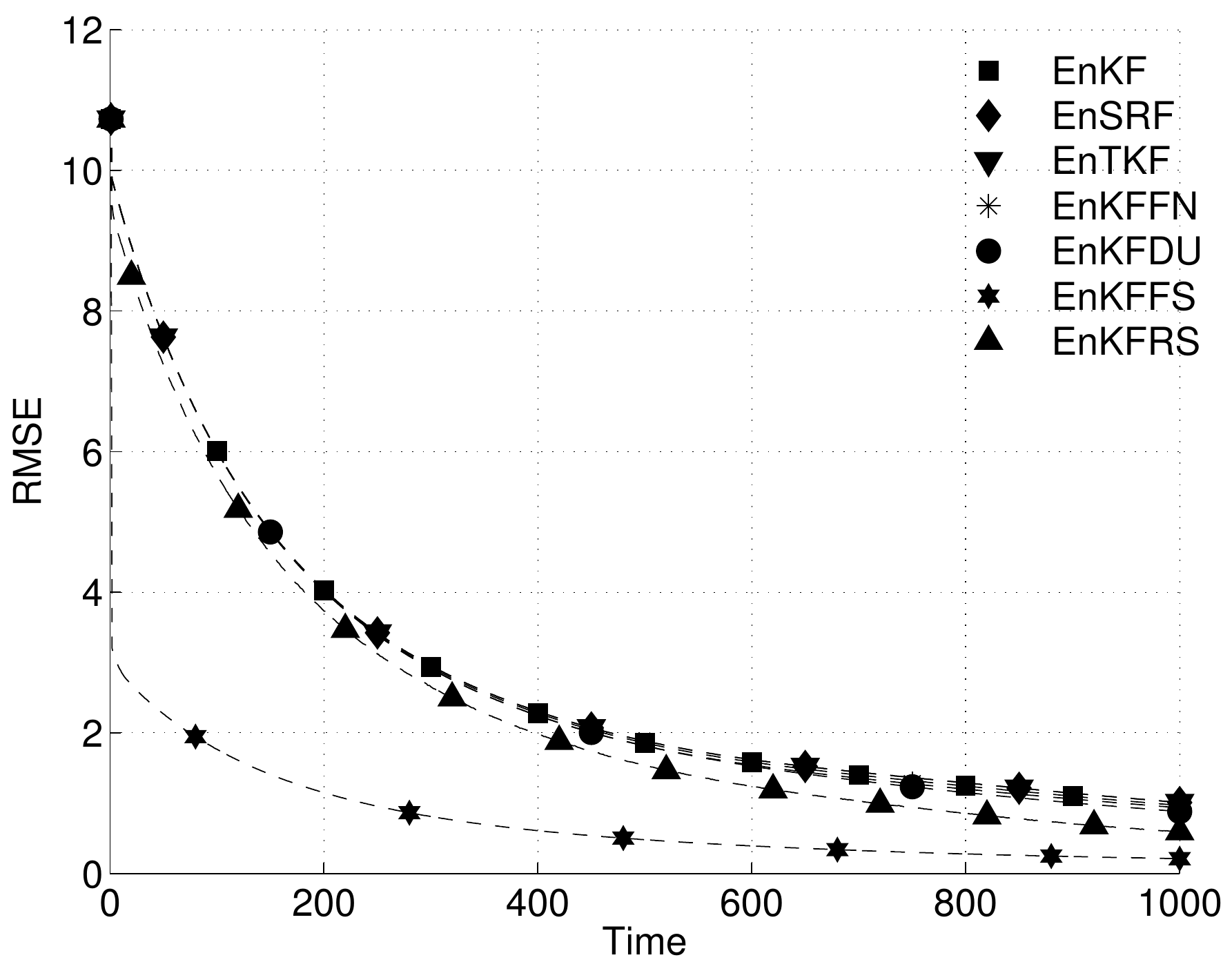}\label{fig-80-07-001-005-6}} %
          \subfloat[$\Nstate =  3969$ and $\sigma^B = 0.15$]{\includegraphics[width=0.4\textwidth]{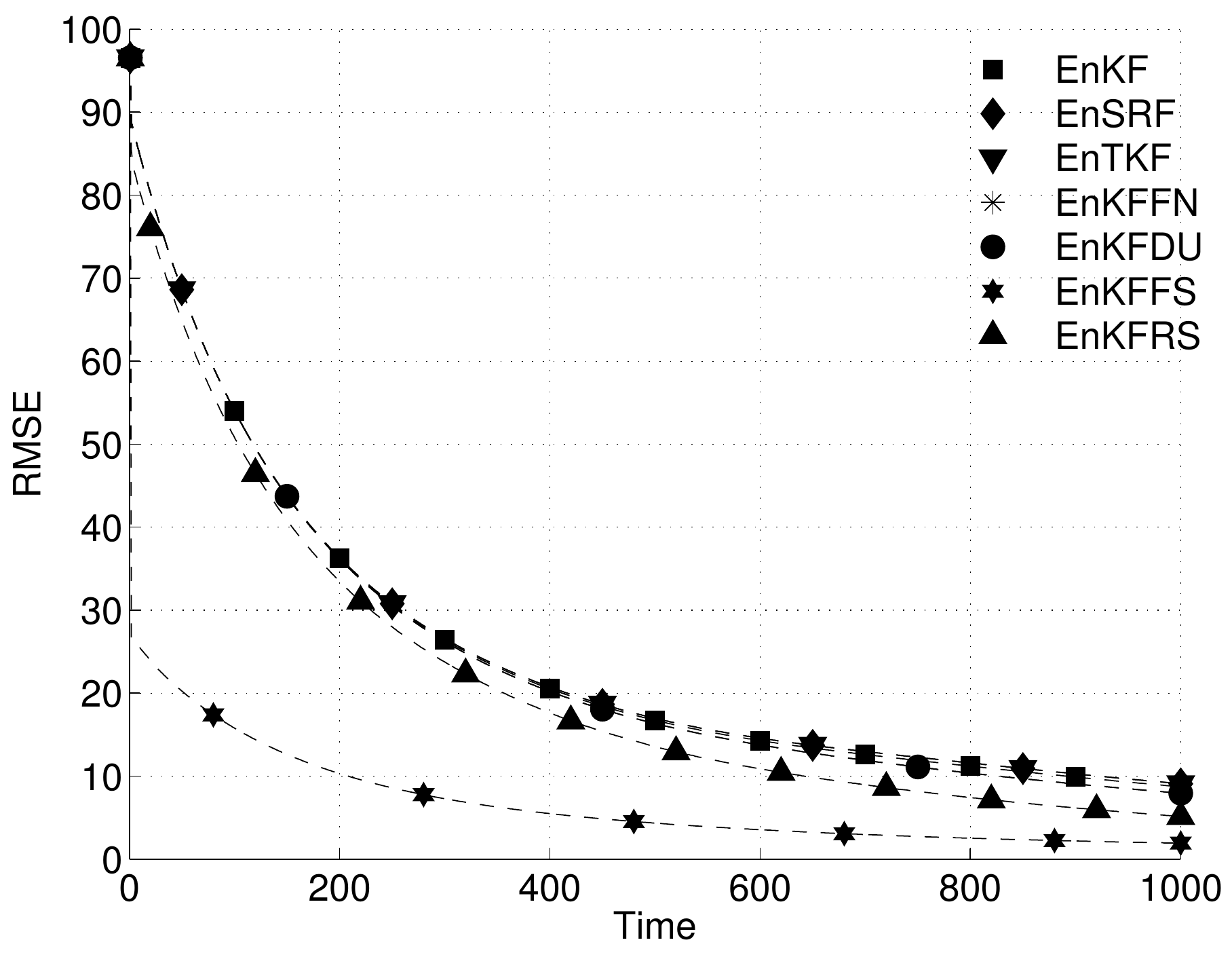}\label{fig-80-07-001-015-6}}
          
     \subfloat[$\Nstate =  16129$ and $\sigma^B = 0.05$]{\includegraphics[width=0.4\textwidth]{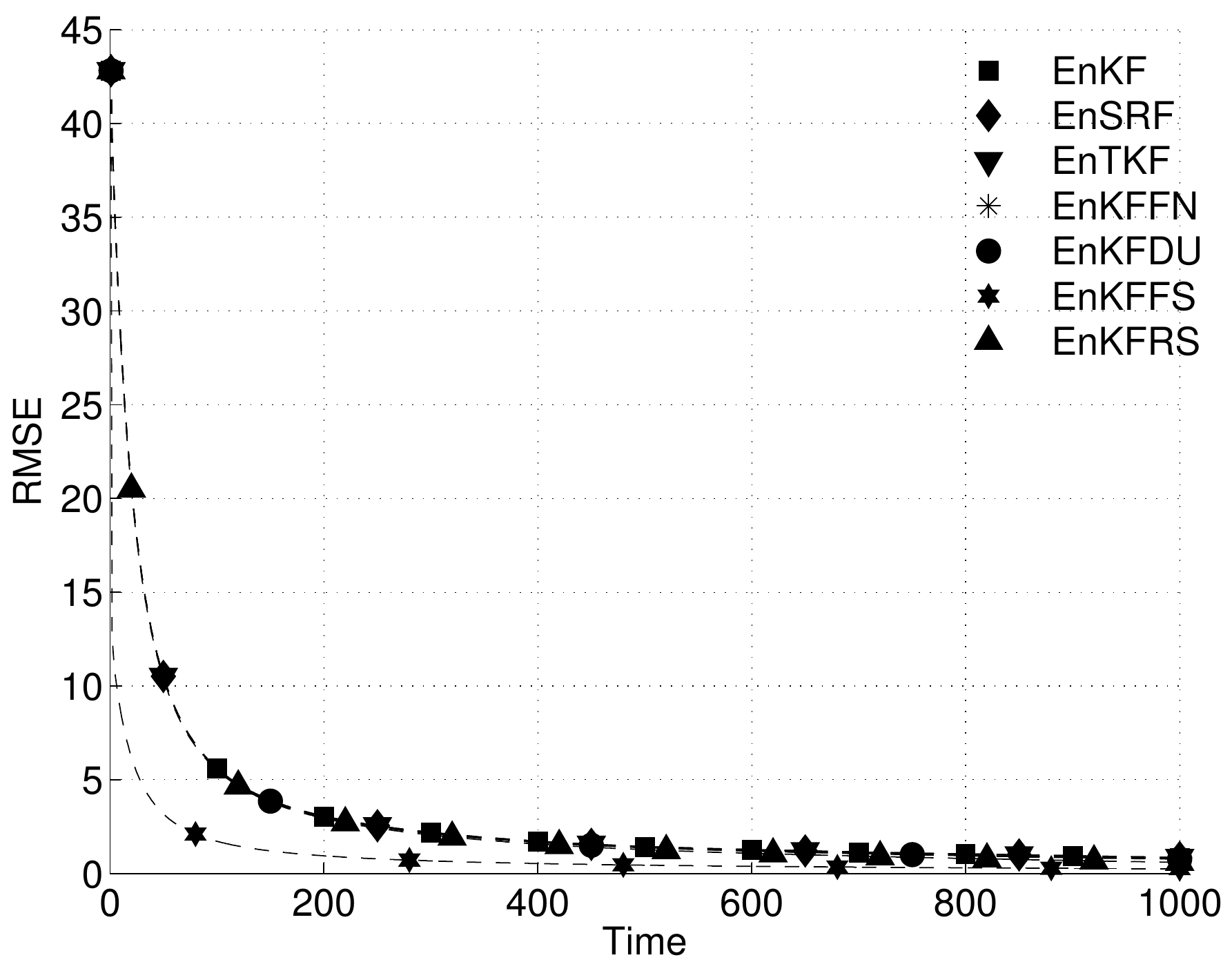}\label{fig-80-07-001-005-7}} %
          \subfloat[$\Nstate =  16129$ and $\sigma^B = 0.15$]{\includegraphics[width=0.4\textwidth]{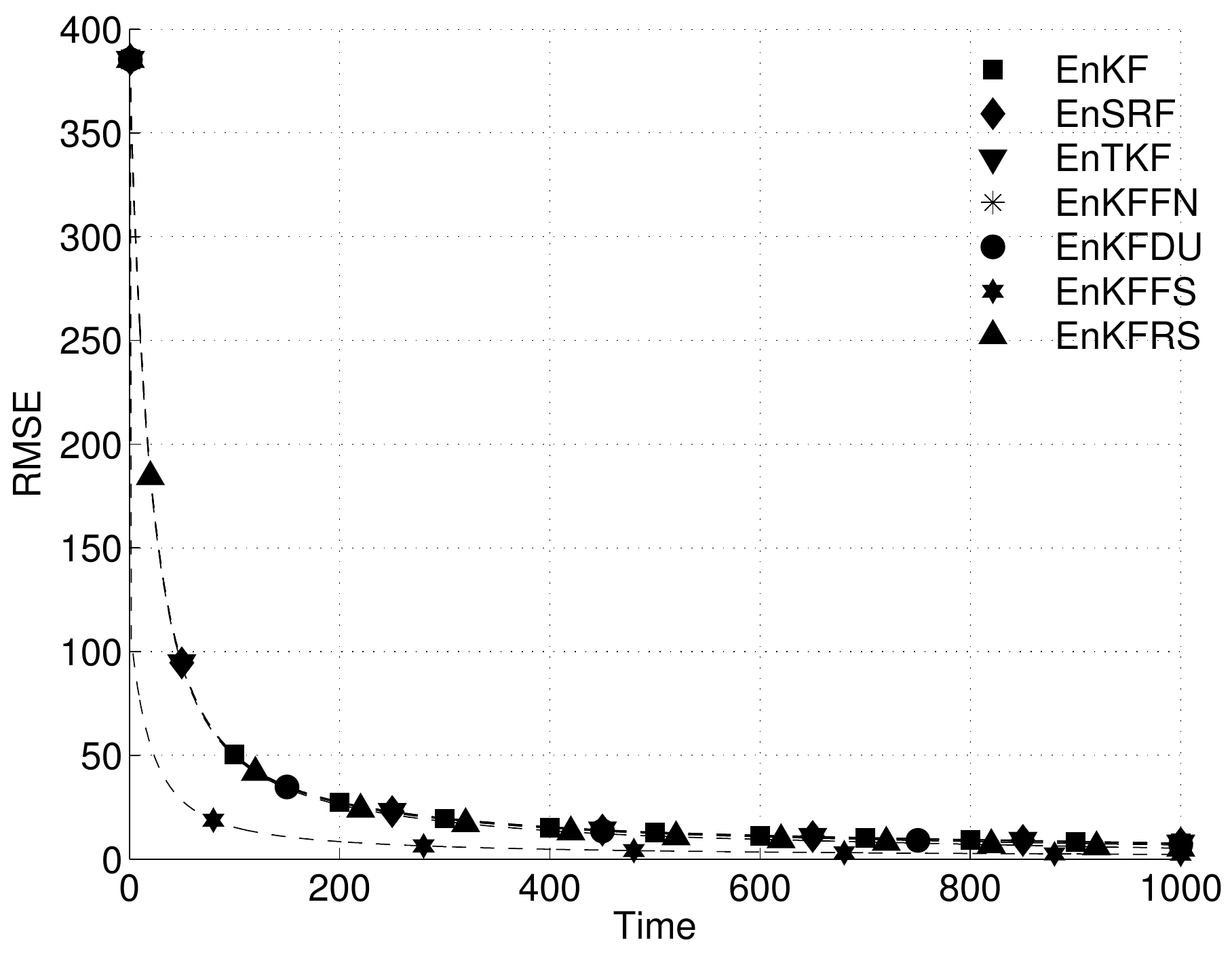}\label{fig-80-07-001-015-7}}
     \caption{Plots of RMSE values of the EnKF, EnSRF, EnTKF, EnKF-FN, EnKF-DU, EnKF-FS and EnKF-RS for $\Nens = 80$ and $\Pobs = 0.7$.}
     \label{fig:Comparison-Nens-80}
\end{figure}
%
\subsection{The impact of the number of ensemble members on performance}

The EnKF-RS and EnKF-FS implementations depend on the $K$ samples taken from the distribution \eqref{eq:artificial-member}. We now study how the performance of the proposed methods varies for different values of $\Nens$ (number of real members) and $K$ (number of artificial members). For this we let $K$ and $\Nens$ to be related by
\begin{eqnarray*}
K = C \cdot \Nens \,,
\end{eqnarray*}
where $C$ is a constant. Practical ensemble sizes range in $40 \le \Nens \le 80$ \cite{EnKF1657419}. For each ensemble size we use several values of $C$ between $0$ and $10$, e.g., $\Nens = 40$ and $C=10$ lead to $K = 400$. When $C=0$ no artificial members are added, but the error covariance matrix $\BEST$ is estimated. In the numerical experiments, the variances of the initial background error are set to $\sigma^{\B} = 0.15$. The analyses RMSE values and the compute times for the proposed implementations using the $\qga$ instance are reported in figures \ref{fig:RMSE-proposed-methods} and \ref{fig:TIME-proposed-methods}, respectively. The analysis times for both implementations are small. Moreover, as expected, EnKF-RS is sensitive to changes in any of the parameters $\Nens$ and $K$. The RMSE is decreased when the values of those parameters are high as shown in figures \ref{ERMSE_5_7_07} and \ref{ERMSE_5_7_09}. The RMSE of the EnKF-FS analysis decreases only with increasing $\Nens$ as can be seen in figures \ref{ERMSE_5_6_07} and \ref{ERMSE_5_6_09}.

An important question is how well do the proposed implementations perform with a small number of  real members. Hopefully the inexpensive addition of artificial members can compensate for a small number of real ones. To this end we consider a small number of real members $10 \le \Nens \le 30$ and a large number of artificial members with $10 \le C \le 60$. The results for the EnKF-RS are shown in figures \ref{ERMSE2_5_6_07} and \ref{ERMSE2_5_6_09} for $p$ equal to 0.7 and 0.9, respectively. The EnKF-RS implementation improves the estimated analysis state whenever $\Nens$ or $K$ are increased. Moreover, the quality of the analyses obtained with small real ensembles ($10 \le \Nens \le 30$) is comparable to those obtained with large real ensemble sizes ($40 \le \Nens \le 80$). This justifies the addition of inexpensive artificial members in order to increase the degrees of freedom of the ensemble. The EnKF-FS analysis improves only when the number of real members is increased. This can be seen in figures \ref{ERMSE2_5_6_07} and \ref{ERMSE2_5_6_09}. For the smallest ensemble size ($\Nens =10$) the results obtained by the EnKF-FS are better than those of any other implementation, including the traditional EnKF implementations with the large ensemble size $\Nens=80$.

Figures \ref{ERMSE2_5_7_07} and \ref{ERMSE2_5_7_09} show that the EnKF-RS analyses for a small number of observed components (70\%) and a large number of artificial members ($K$) are equivalent to those obtained with a large number of real members ($\Nens$) and many observed components (90\%). This is another indication of the positive impact obtained by increasing the number of degrees of freedom  with samples from the distribution \eqref{eq:artificial-member}. This is computationally less expensive than adding real members via running the model.

\begin{figure}[H]
     \centering
     \subfloat[EnKF-FS $\Pobs = 0.7$]{\includegraphics[width=0.49\textwidth]{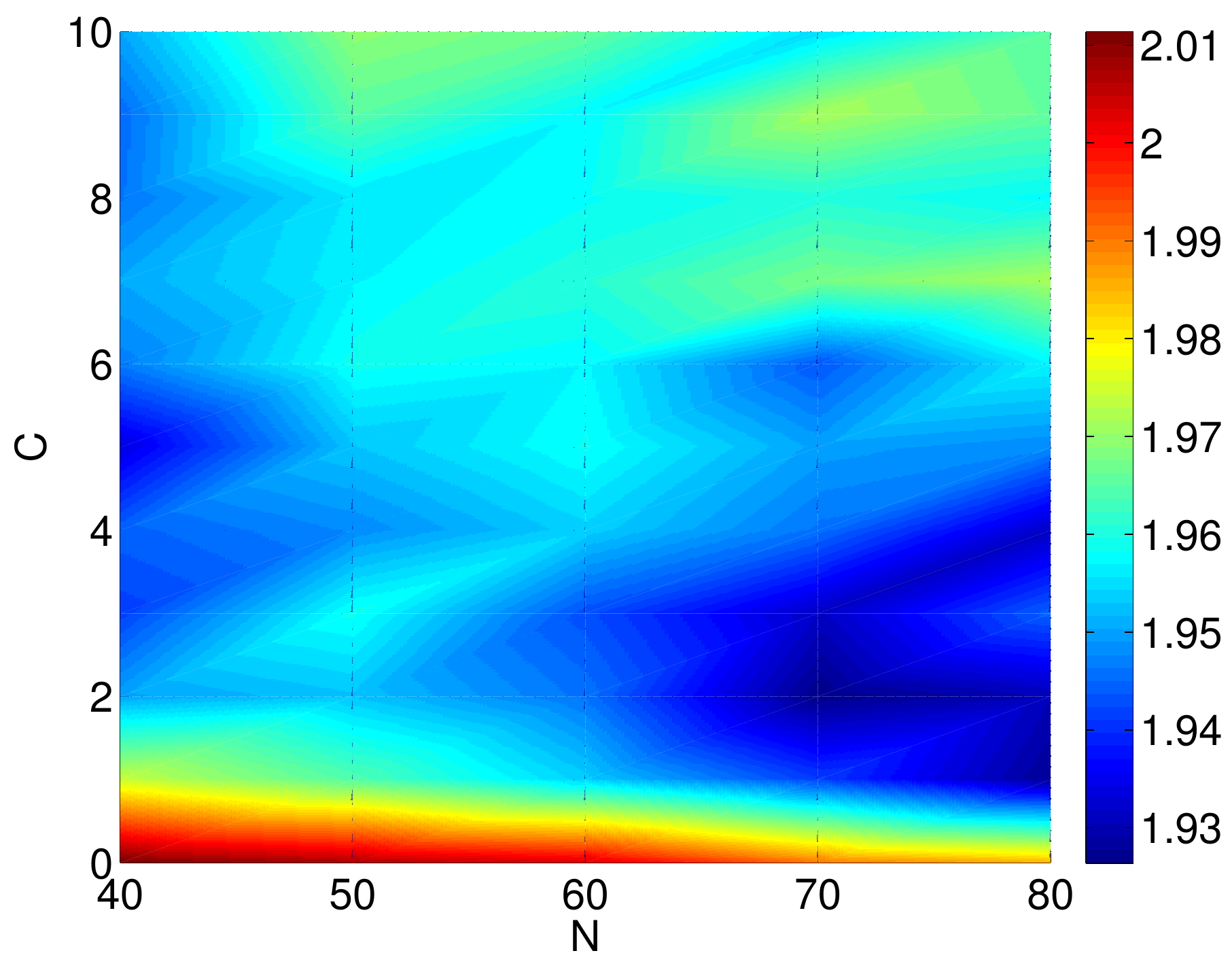}\label{ERMSE_5_6_07}} %
          \subfloat[EnKF-RS $\Pobs = 0.7$]{\includegraphics[width=0.48\textwidth]{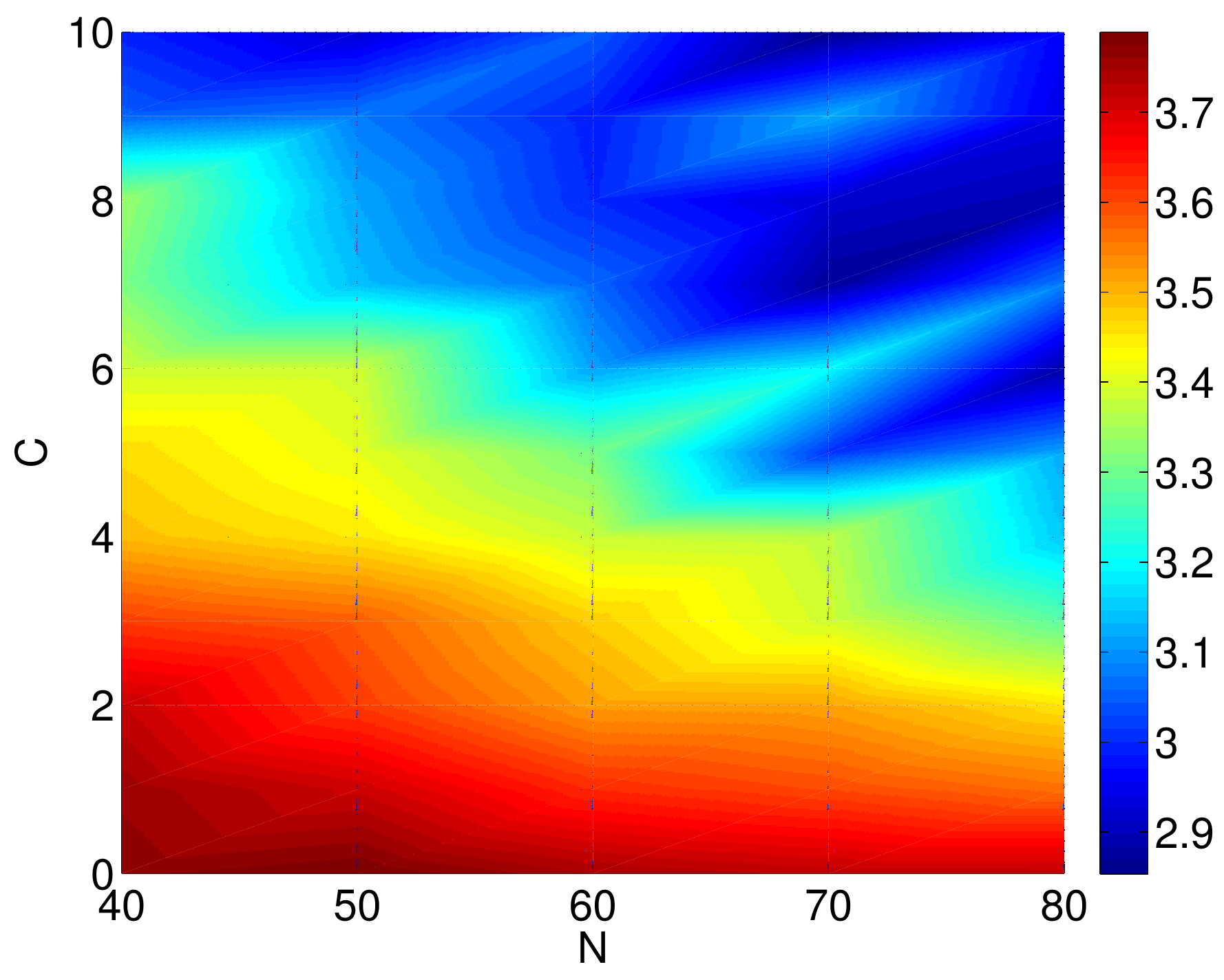}\label{ERMSE_5_7_07}}
          
     \subfloat[EnKF-FS $\Pobs = 0.9$]{\includegraphics[width=0.49\textwidth]{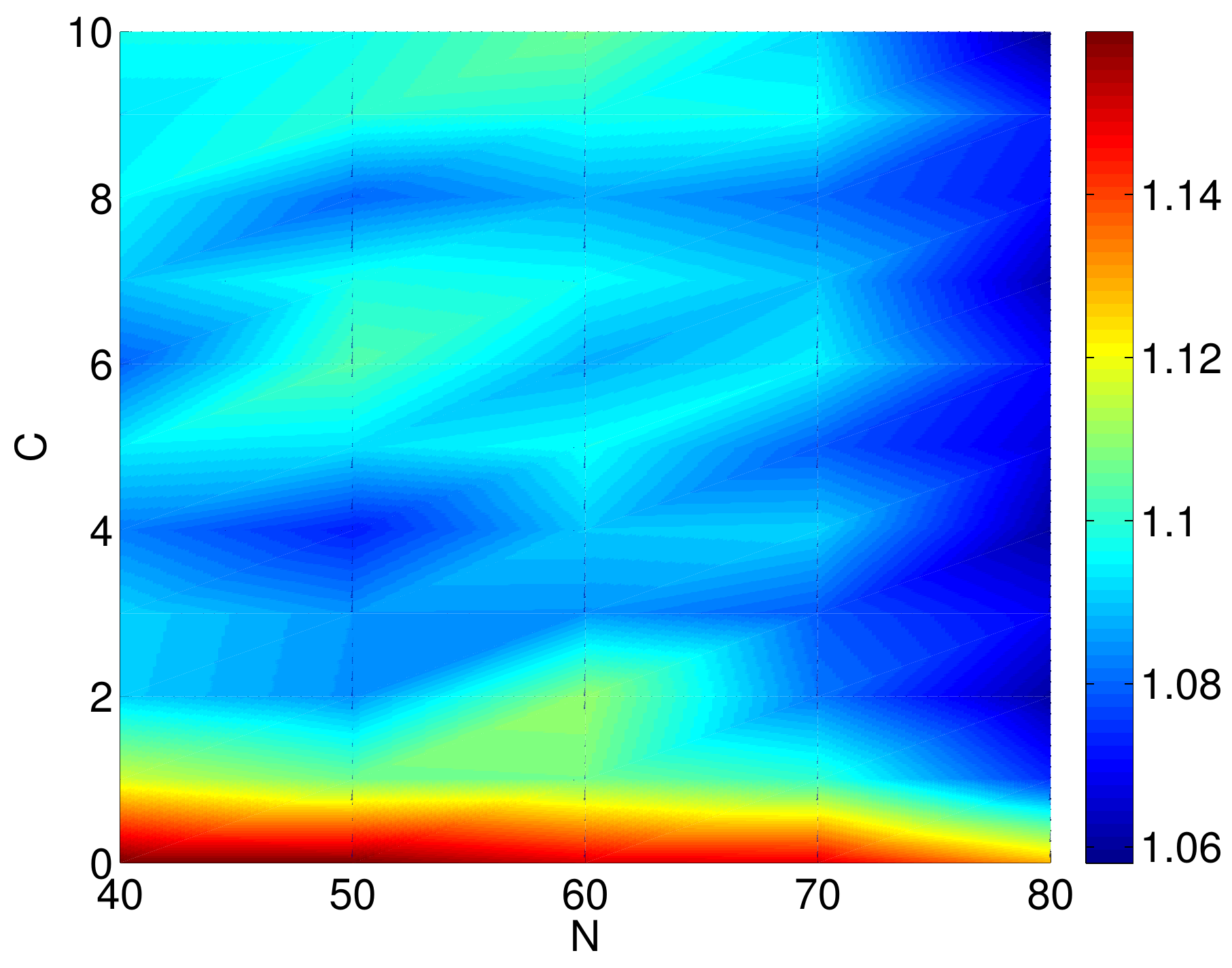}\label{ERMSE_5_6_09}} %
          \subfloat[EnKF-RS $\Pobs = 0.9$]{\includegraphics[width=0.48\textwidth]{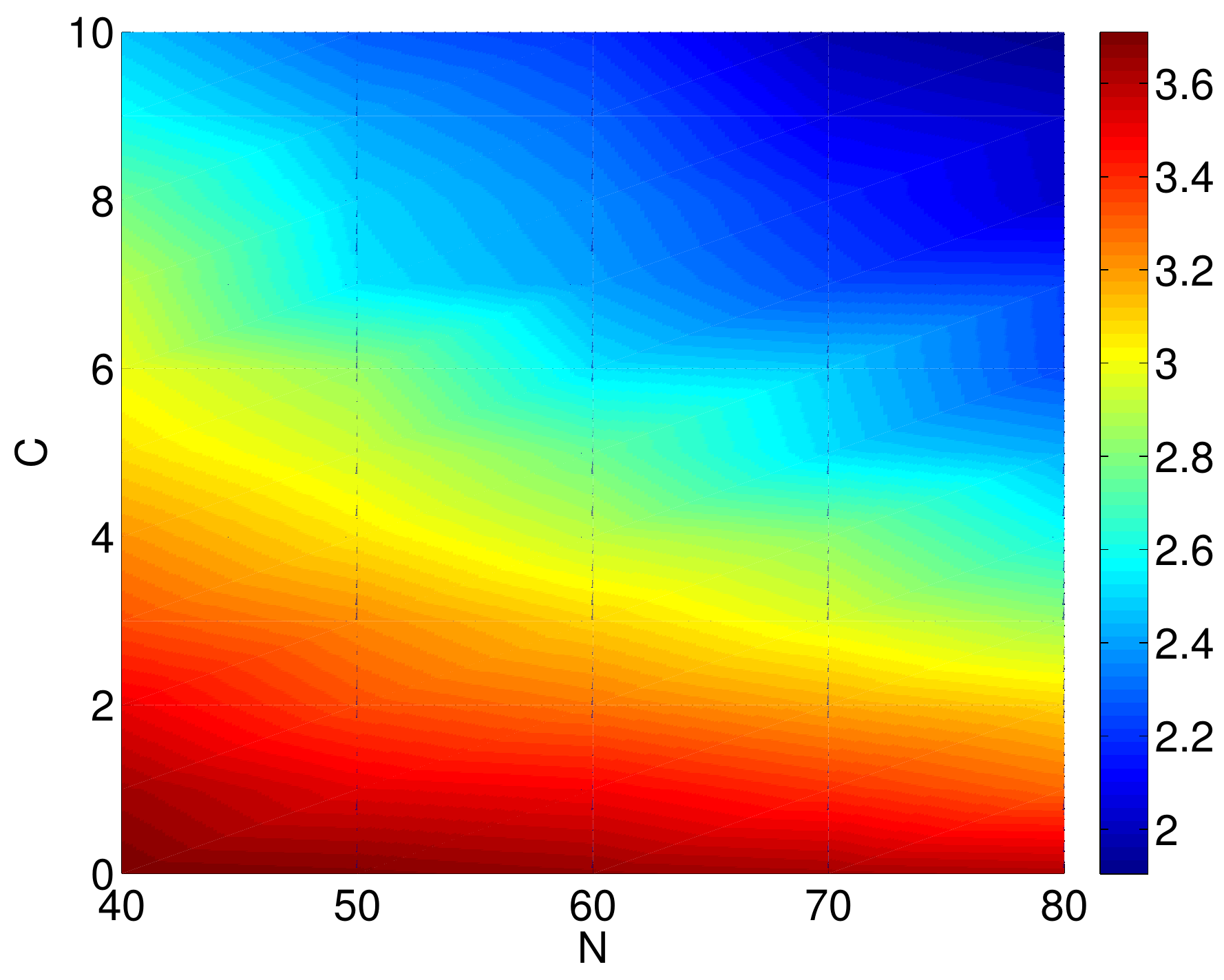}\label{ERMSE_5_7_09}}
         \caption{RMSE of the EnKF-FS and EnKF-RS implementations for different values of $0 \le C \le 10$ and $40 \le \Nens \le 80$ making use of the $\qga$ instance.}
     \label{fig:RMSE-proposed-methods}
\end{figure}
\begin{figure}[H]
     \centering
     \subfloat[EnKF-FS, $\Pobs = 0.7$]{\includegraphics[width=0.49\textwidth]{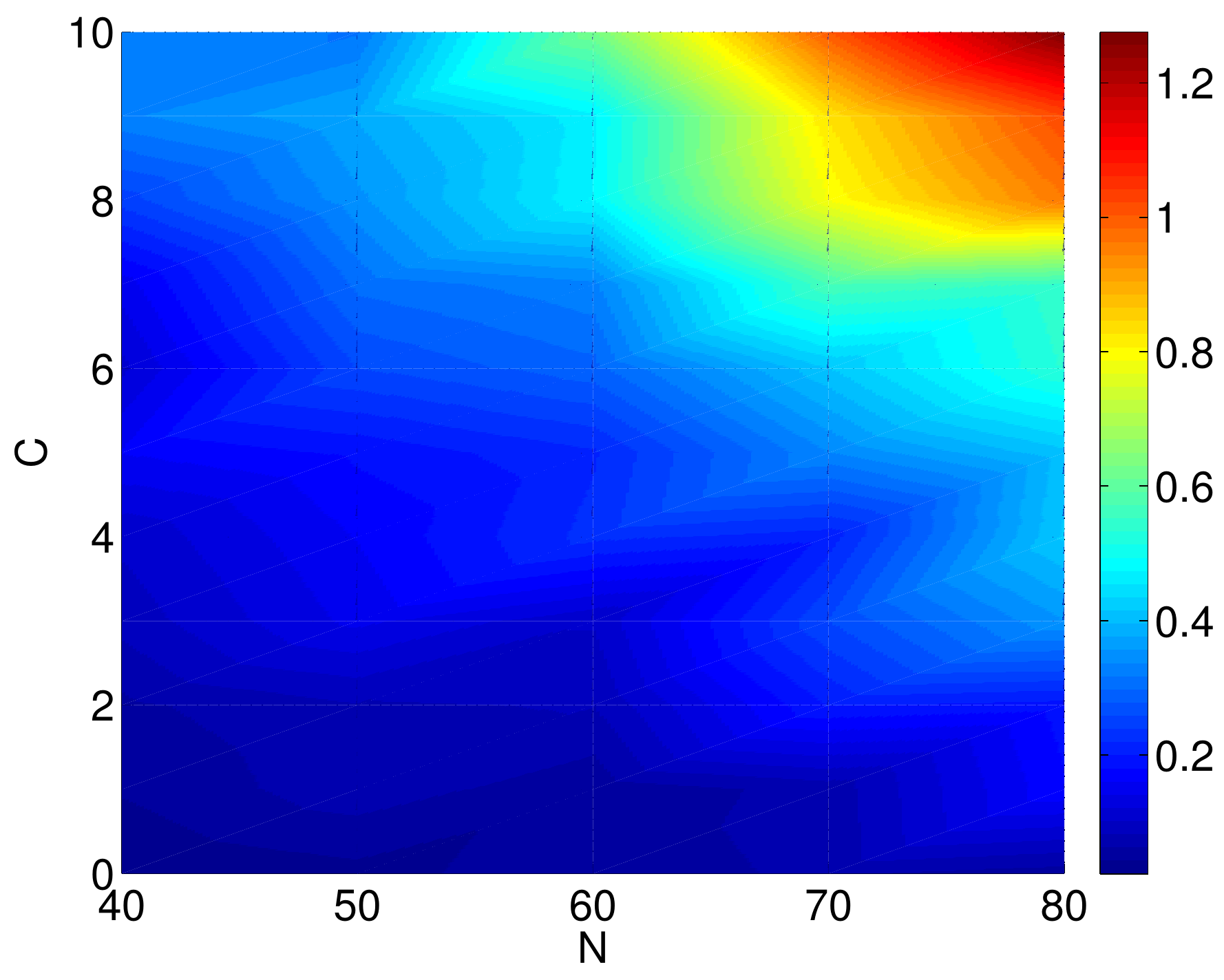}\label{ETIME_5_6_07}} %
          \subfloat[EnKF-RS, $\Pobs = 0.7$]{\includegraphics[width=0.49\textwidth]{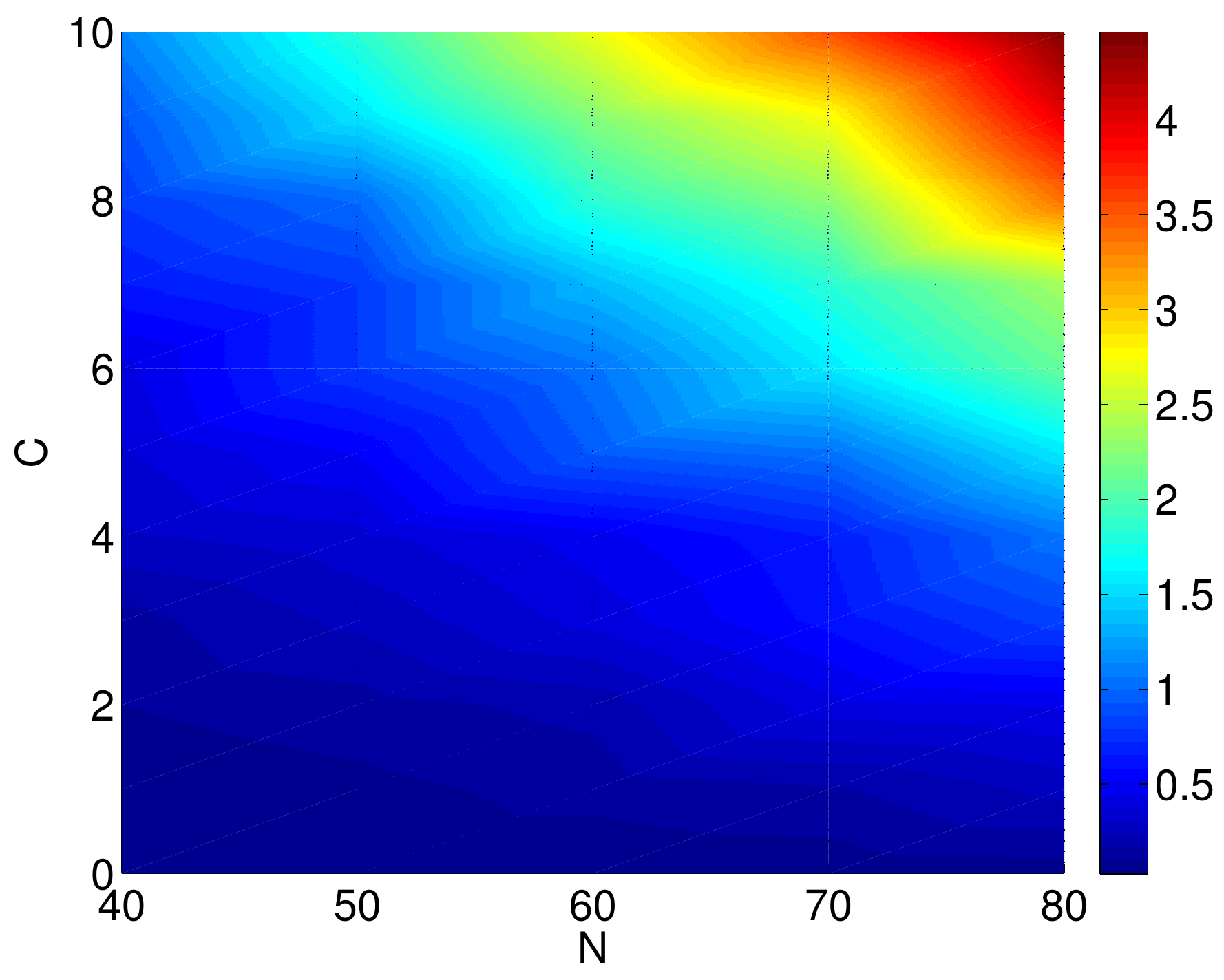}\label{ETIME_5_7_07}}
          
     \subfloat[EnKF-FS, $\Pobs = 0.9$]{\includegraphics[width=0.49\textwidth]{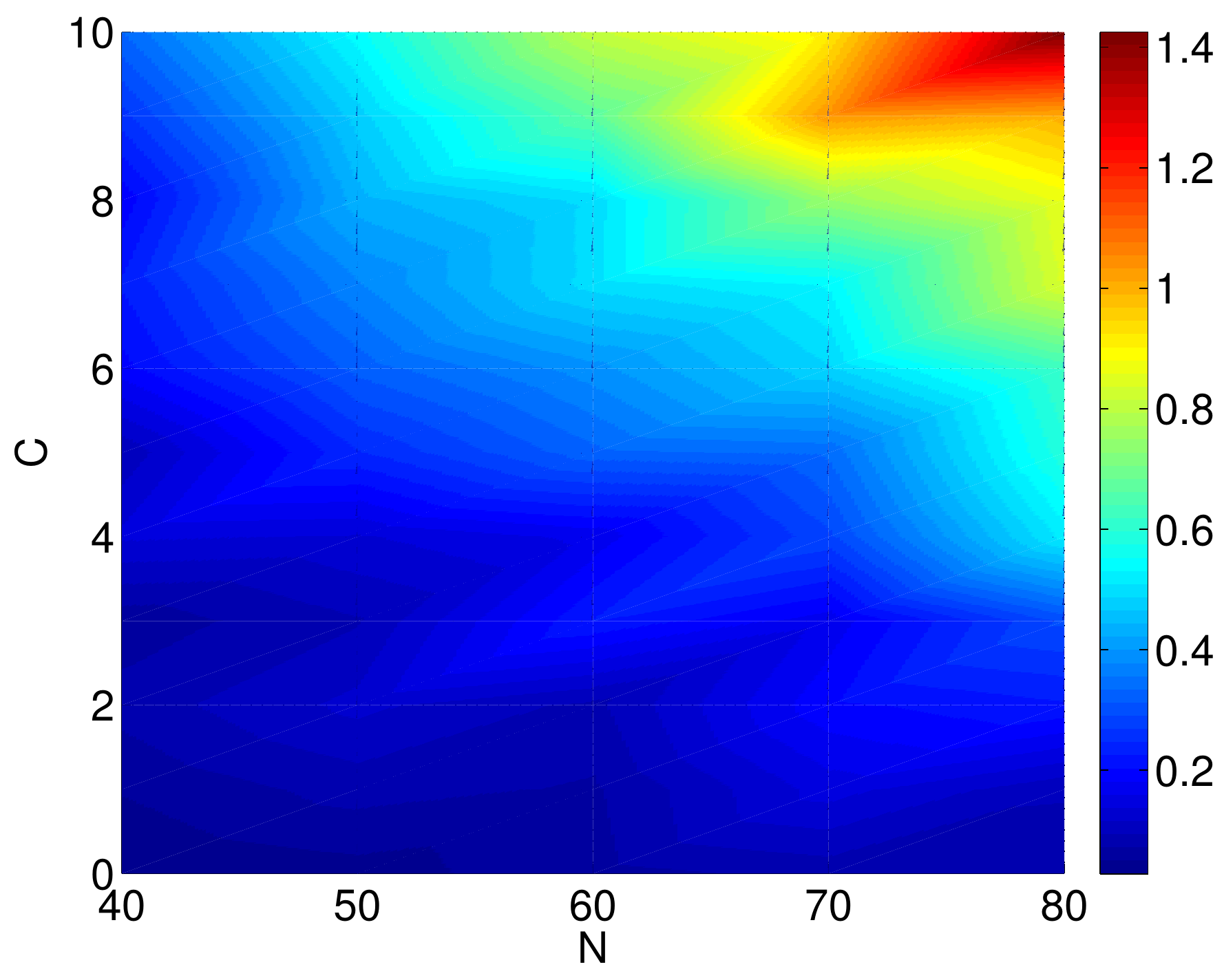}\label{ETIME_5_6_09}} %
          \subfloat[EnKF-RS, $\Pobs = 0.9$]{\includegraphics[width=0.48\textwidth]{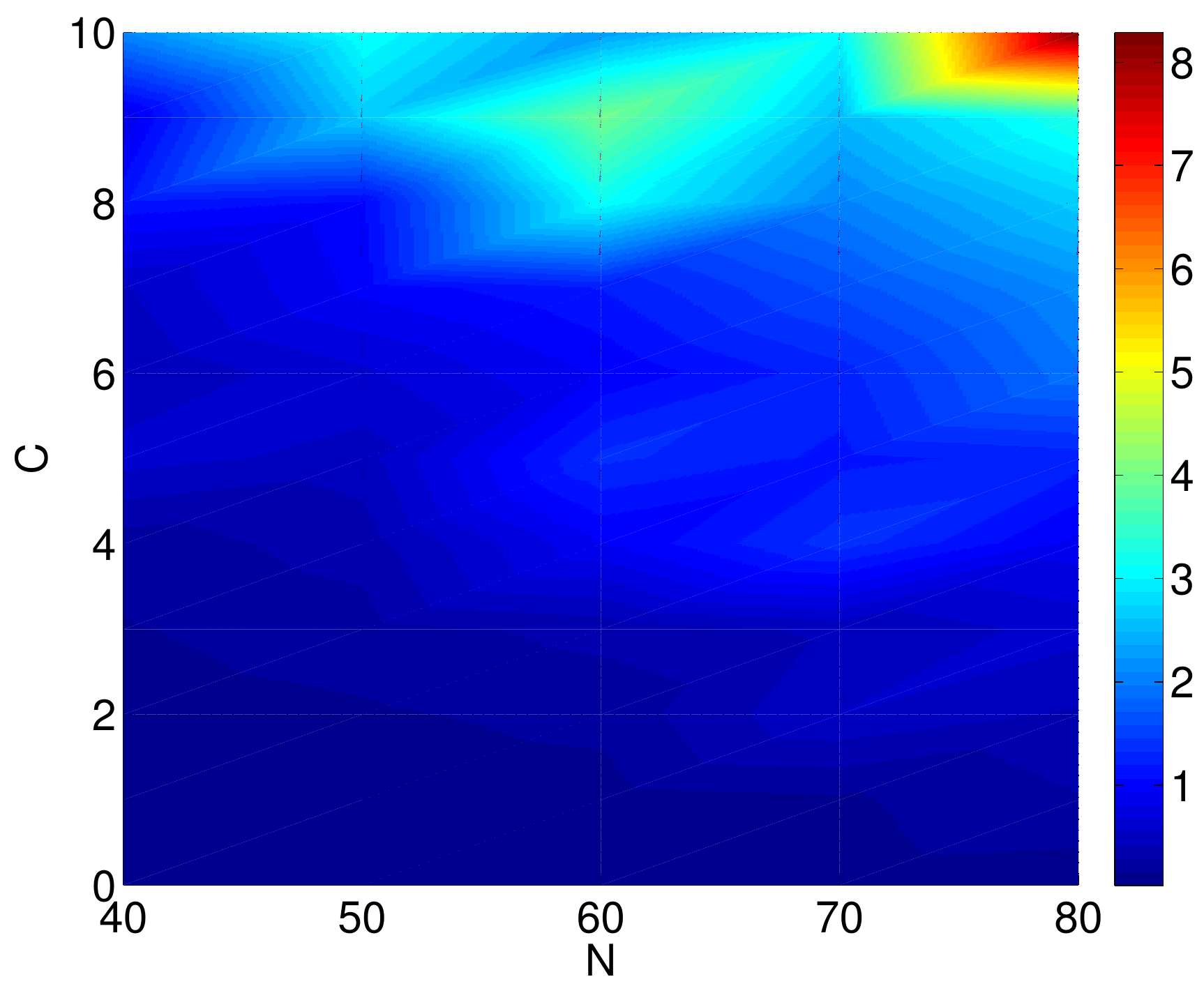}\label{ETIME_5_7_09}}
         \caption{Assimilation times of the EnKF-FS and EnKF-RS implementations for different values of $0 \le C \le 10$ and $40 \le \Nens \le 80$ making use of the $\qga$ instance}
     \label{fig:TIME-proposed-methods}
\end{figure}
\begin{figure}[H]
     \centering
     \subfloat[EnKF-FS, $\Pobs = 0.7$]{\includegraphics[width=0.49\textwidth]{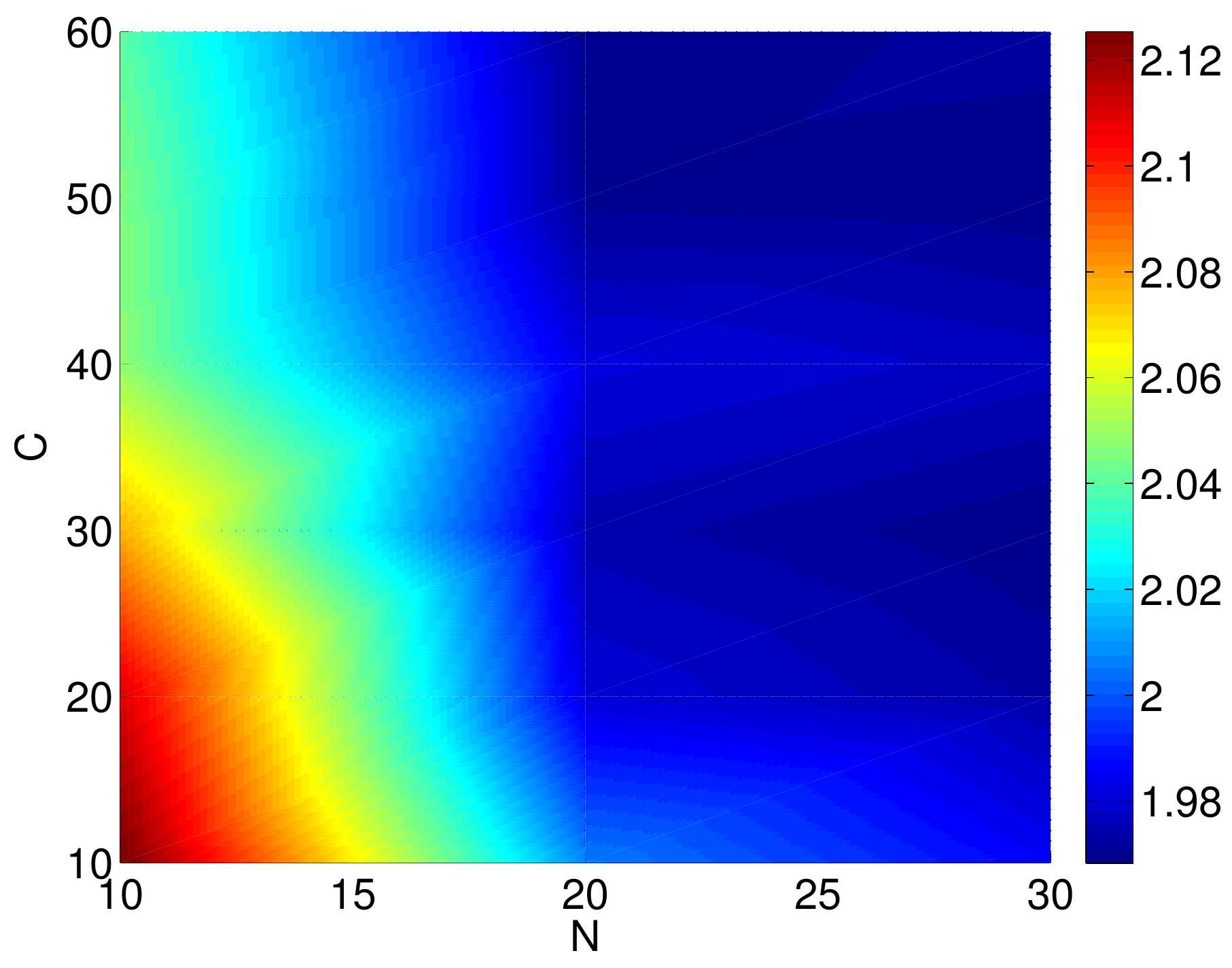}\label{ERMSE2_5_6_07}} %
          \subfloat[EnKF-RS, $\Pobs = 0.7$]{\includegraphics[width=0.48\textwidth]{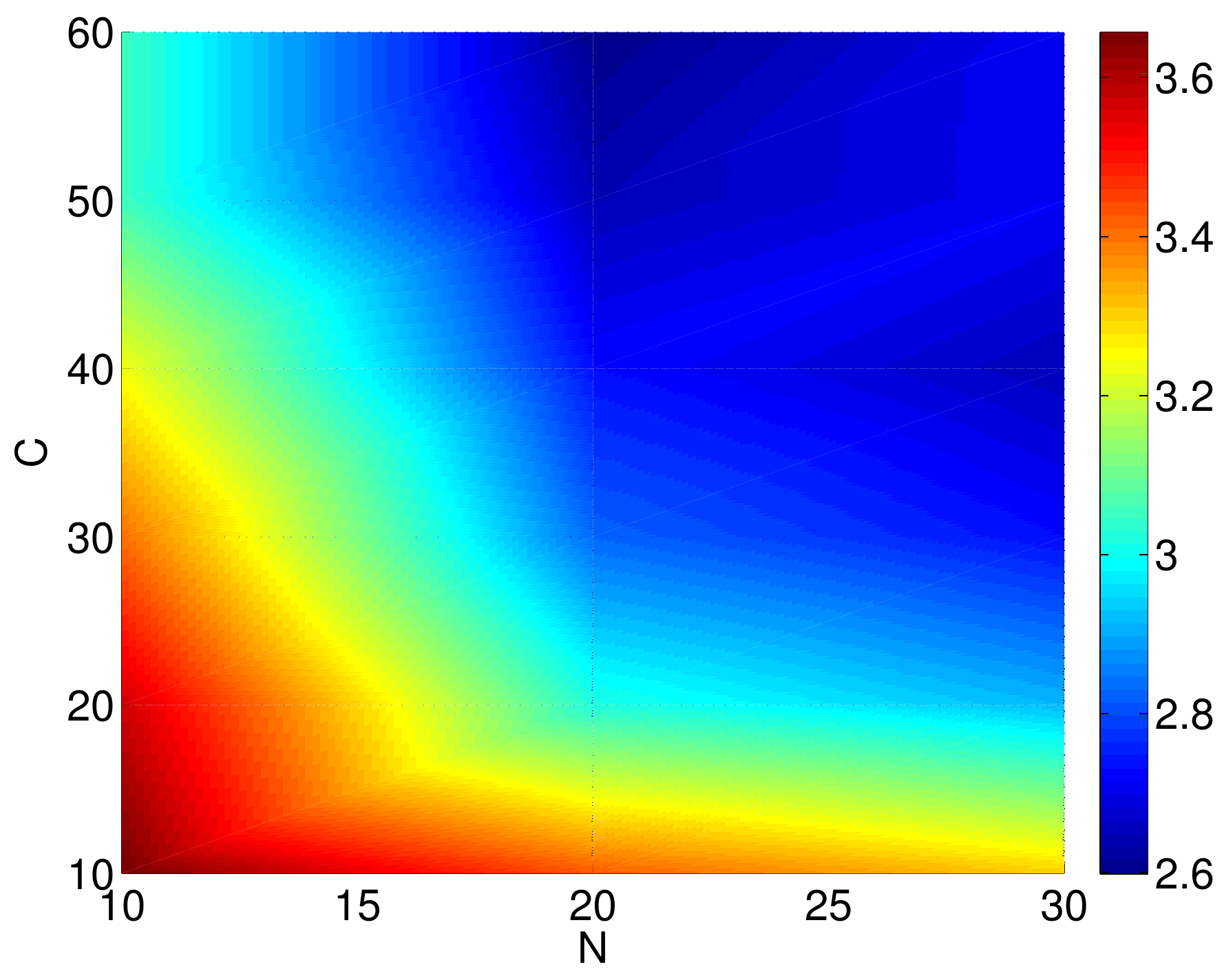}\label{ERMSE2_5_7_07}}
          
     \subfloat[EnKF-FS, $\Pobs = 0.9$]{\includegraphics[width=0.49\textwidth]{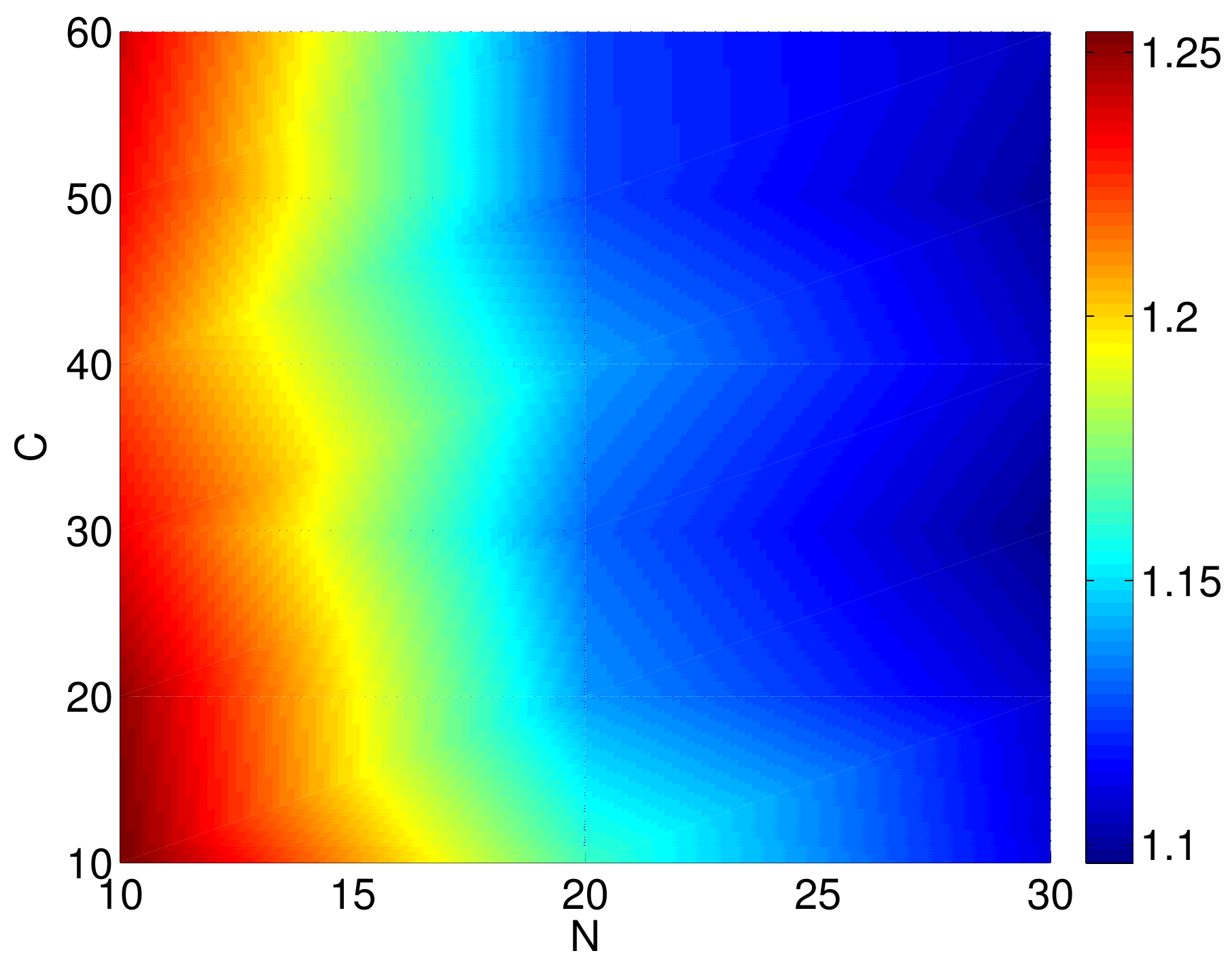}\label{ERMSE2_5_6_09}} %
          \subfloat[EnKF-RS, $\Pobs = 0.9$]{\includegraphics[width=0.48\textwidth]{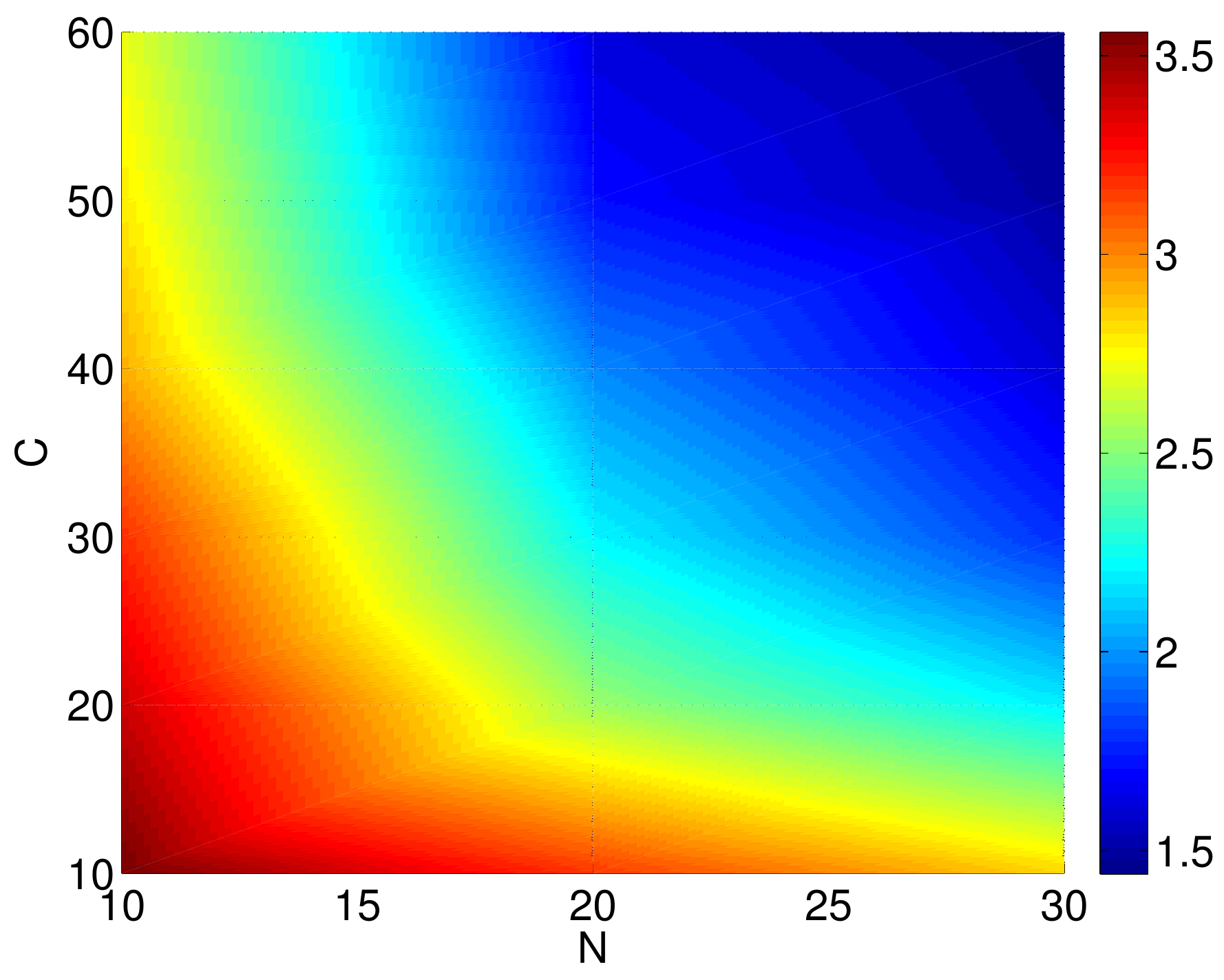}\label{ERMSE2_5_7_09}}
         \caption{RMSEs of the EnKF-FS and EnKF-RS implementations for different values of $10 \le C \le 60$ and $10 \le \Nens \le 30$ making use of the $\qga$ instance.}
     \label{fig:RMSE2-proposed-methods}
\end{figure}
\begin{figure}[H]
     \centering
     \subfloat[EnKF-FS, $\Pobs = 0.7$]{\includegraphics[width=0.49\textwidth]{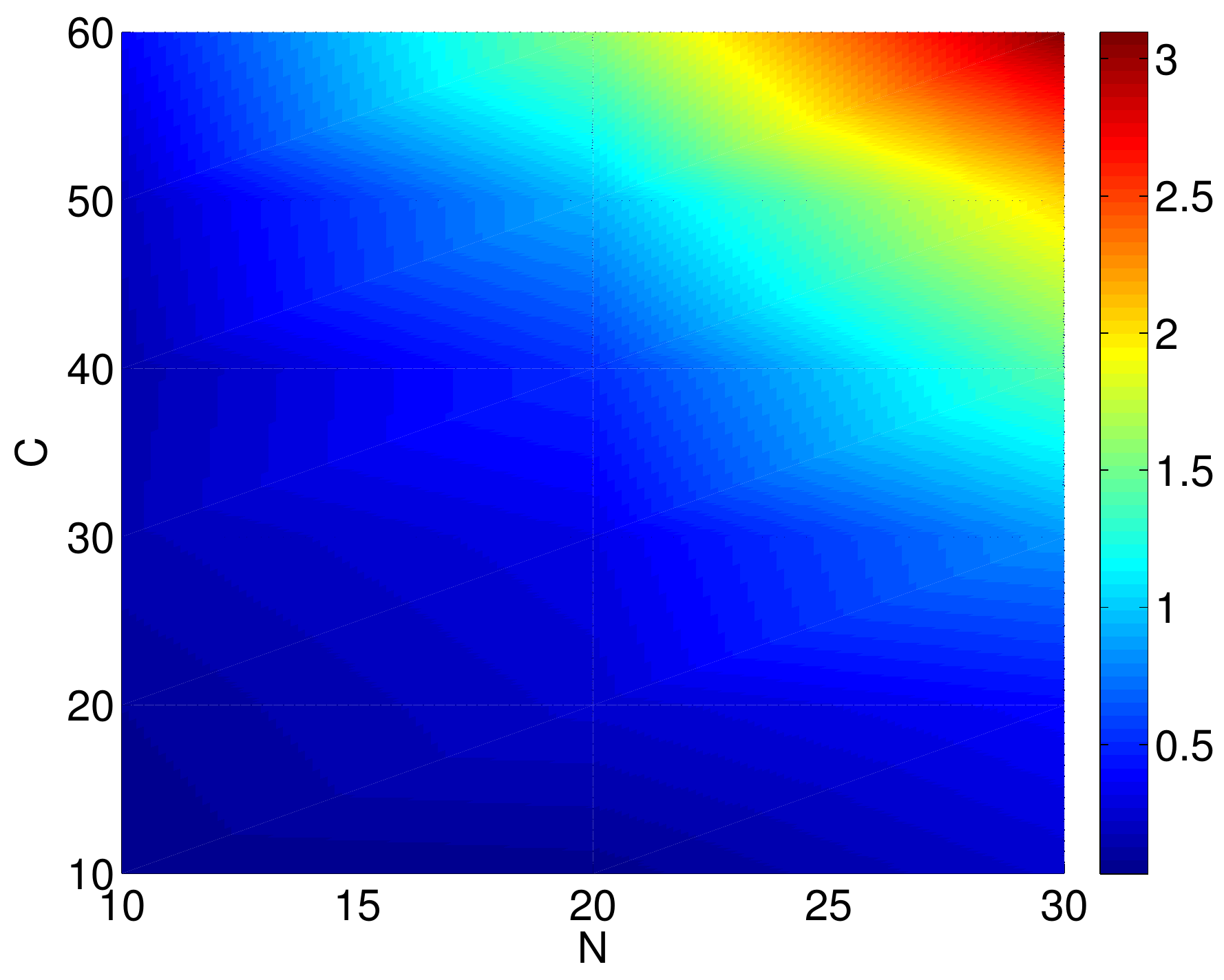}\label{ETIME2_5_6_07}} %
          \subfloat[EnKF-RS, $\Pobs = 0.7$]{\includegraphics[width=0.48\textwidth]{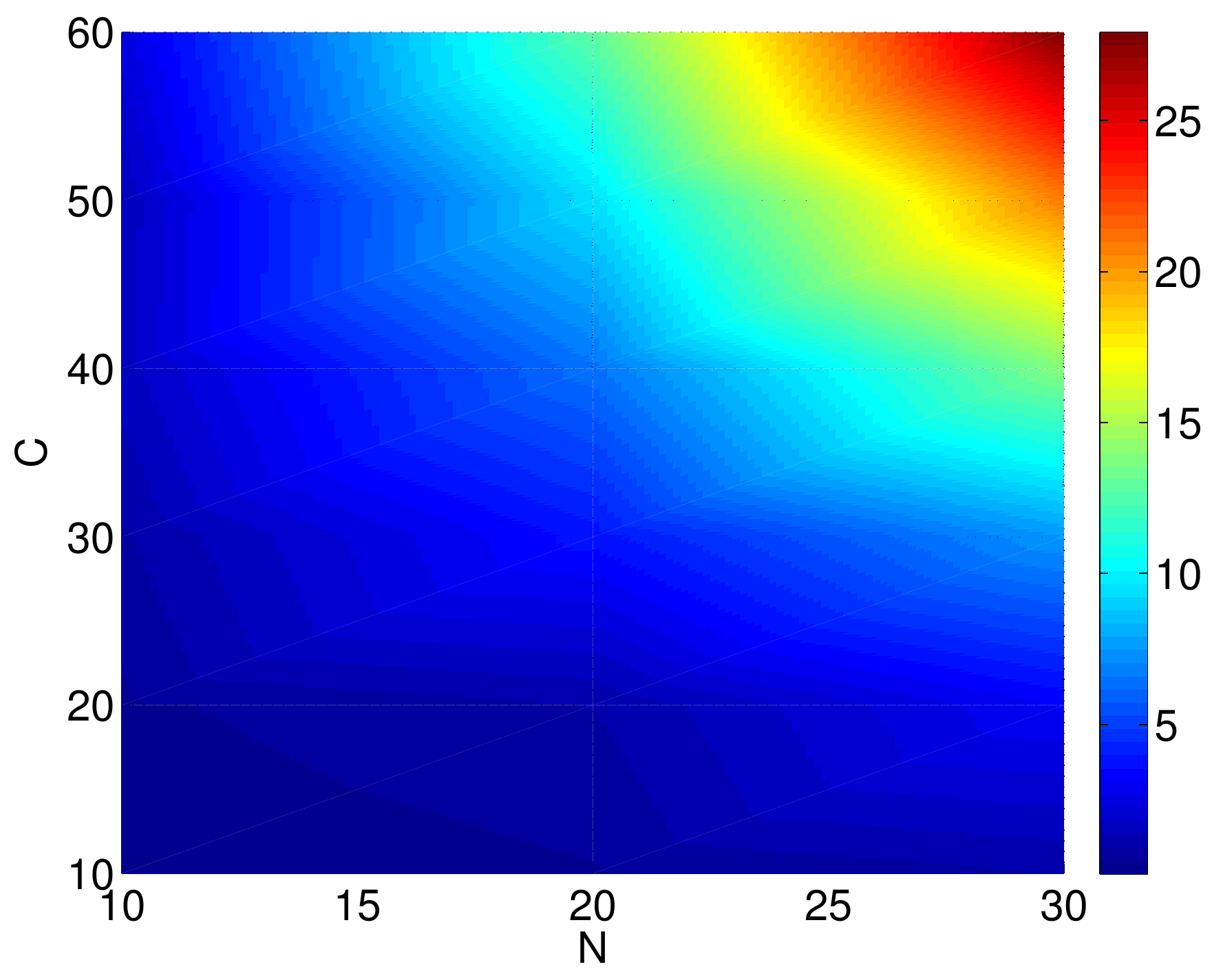}\label{ETIME2_5_7_07}}
          
     \subfloat[EnKF-FS, $\Pobs = 0.9$]{\includegraphics[width=0.49\textwidth]{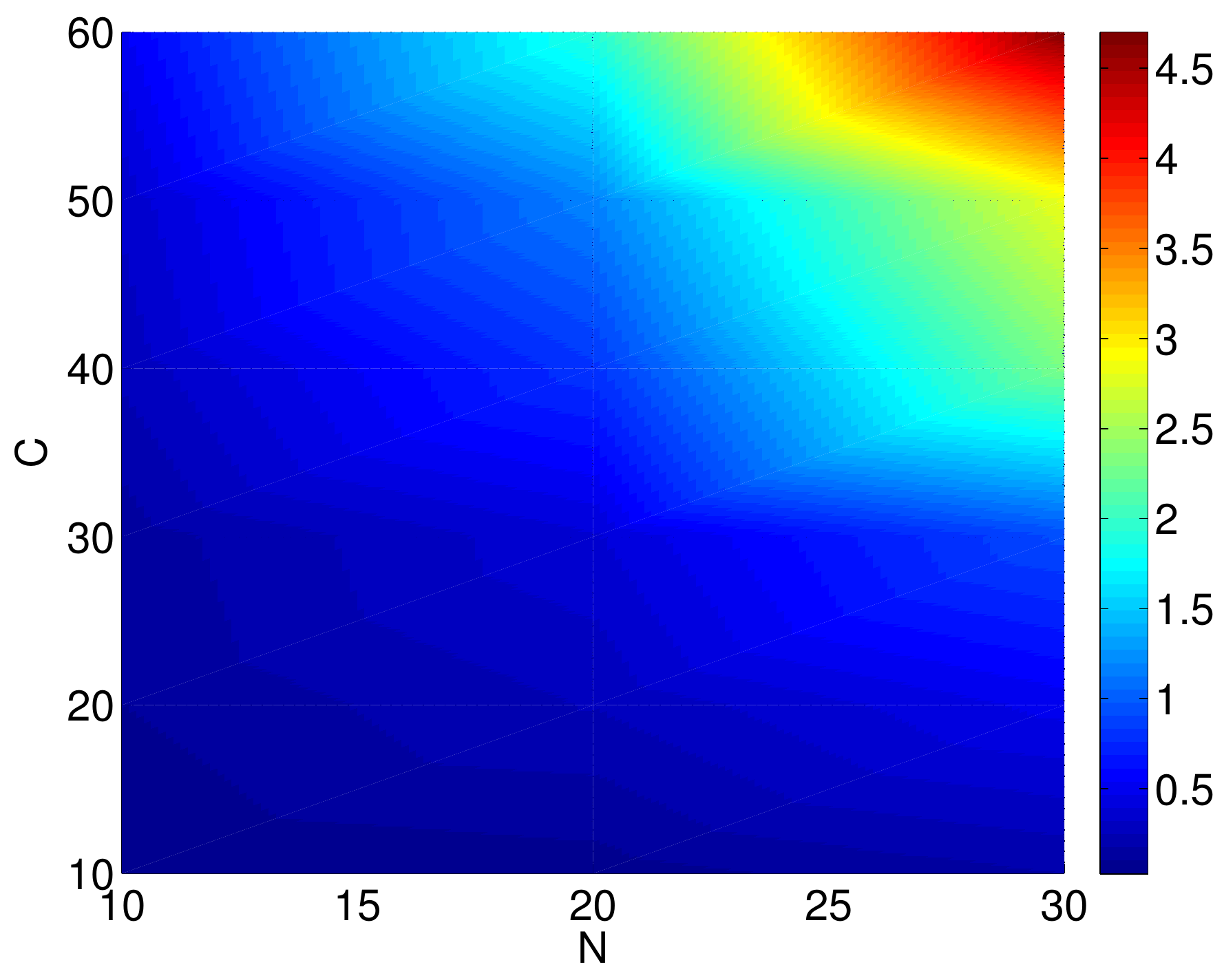}\label{ETIME2_5_6_09}} %
          \subfloat[EnKF-RS, $\Pobs = 0.9$]{\includegraphics[width=0.48\textwidth]{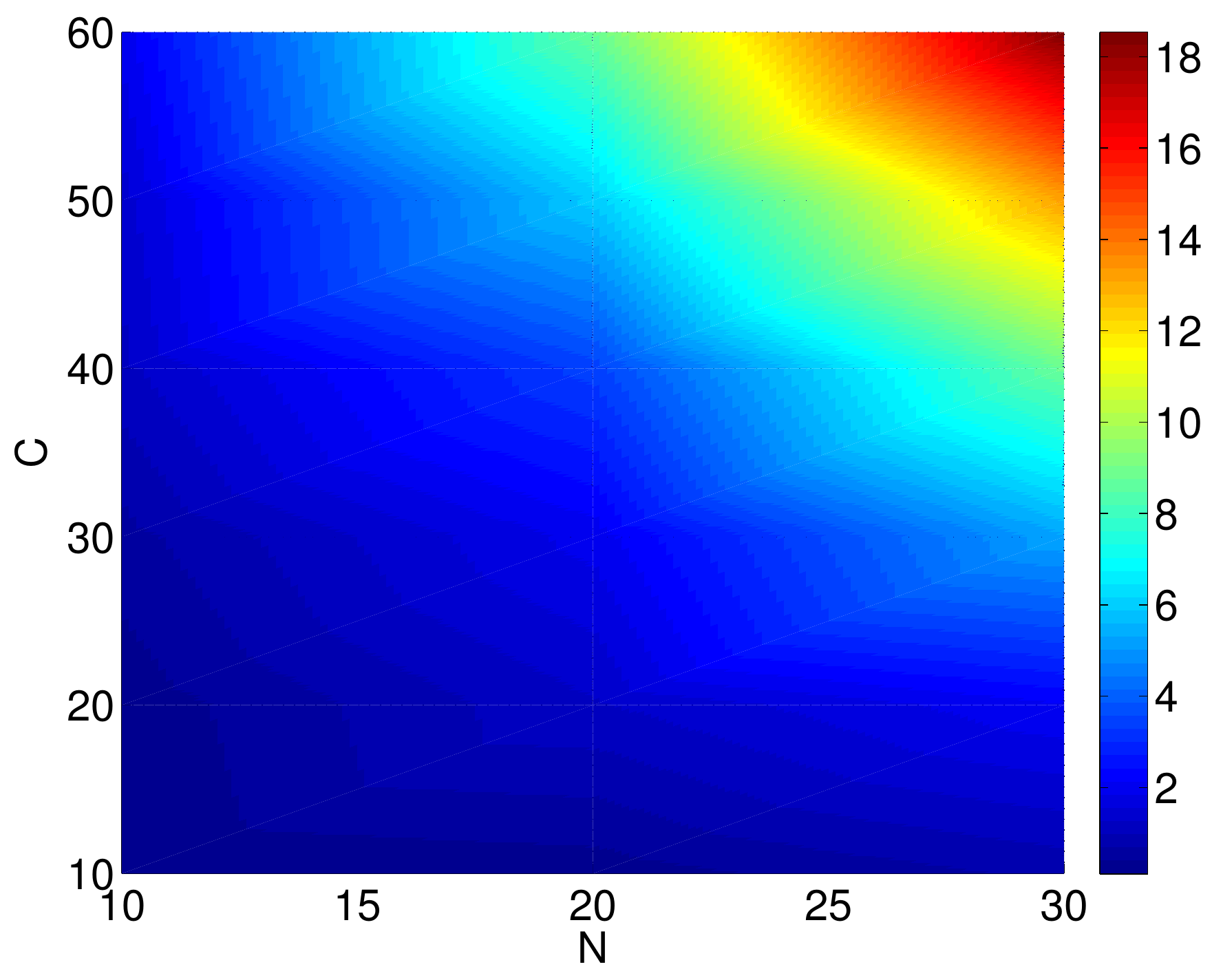}\label{ETIME2_5_7_09}}
         \caption{Assimilation times of the EnKF-FS and EnKF-RS implementations for different values of $10 \le C \le 60$ and $10 \le \Nens \le 30$ making use of the $\qga$ instance.}
     \label{fig:TIME2-proposed-methods}
\end{figure}

\section{Conclusions}
\label{sec:conclusions}

This paper develops two new implementations of the ensemble Kalman filter (EnKF) based on shrinkage covariance estimation. The background error covariance matrices used in analysis are obtained via the Rao-Blackwell Ledoit and Wolf estimator, which has been proved optimal in the estimation of high-dimensional covariance matrices from a small number of samples. This covariance matrix and the background state (ensemble mean) serve as parameters of the normal error distribution associated with the ensemble members. Samples from this distribution are taken in order to increase the number of ensemble members, and therefore, to decrease the sampling error in representing the background error distribution, and to increase the number of degrees of freedom in the assimilation. The two proposed implementations differ in the space where the assimilation process is performed:  EnKF Full-Space (EnKF-FS) performs the analysis in the model space, while EnKF Reduce-Space (EnKF-RS) computes the analysis state in the space spanned by the ensemble members. Numerical experiments are carried out using a quasi-geostrophic model. They show that the two new implementations perform better than current EnKF implementations such as the traditional EnKF, square root filters, and inflation-free EnKF methods. For all the scenarios and experimental settings, the EnKF-FS outperforms the other implementations by at least $60\%$ in terms of accuracy (root mean square error). Moreover, for a small number of ensemble members ($\sim 10$)  and a moderate percentage of observed components from the vector state ($\sim 70 \%$), the solutions obtained by the proposed methods are similar to those obtained by large ensemble sizes ($\sim 80$) and large percentage of observed components ($\sim 90\%$). The computational time for analysis of the proposed implementations is reasonably low. Since the total compute time is dominated by the multiple model runs considerable savings are expected from reducing the number of real ensemble members without deteriorating the quality of the results.

\section*{Acknowledgements}



\bibliographystyle{alpha}


\end{document}